\documentclass[a4paper,oneside,11pt,final]{amsart}

\usepackage{amsmath,amsthm,amssymb}
\usepackage[matrix,arrow]{xy}
\usepackage{showkeys}

\numberwithin{equation}{section}
\allowdisplaybreaks

\theoremstyle{definition}
 \newtheorem {thm}  {Theorem} [section]
 \newtheorem {prop} [thm]     {Proposition}
 
 \newtheorem {lem}  [thm]     {Lemma}
 \newtheorem {cor}  [thm]     {Corollary}
 \newtheorem {conj} [thm]     {Conjecture}
 \newtheorem*{thm*}           {Theorem}
 \newtheorem*{conj*}          {Conjecture}
\theoremstyle{definition}
 \newtheorem {eg}   [thm]     {Example}
 \newtheorem {fact} [thm]     {Fact}
 \newtheorem {dfn}  [thm]     {Definition}
 \newtheorem*{dfn*}           {Definition}
 \newtheorem*{ack}            {Acknowledgements}
 \newtheorem*{nota}           {Notation}
\theoremstyle{remark}
 \newtheorem {rmk}  [thm]     {Remark}

\newcommand{\ep} {\epsilon}

\newcommand{\bbC}{\mathbb{C}}
\newcommand{\bbP}{\mathbb{P}}
\newcommand{\bbQ}{\mathbb{Q}}
\newcommand{\bbR}{\mathbb{R}}
\newcommand{\bbZ}{\mathbb{Z}}

\newcommand{\bfD}{\mathbf{D}}
\newcommand{\bfE}{\mathbf{E}}
\newcommand{\bfF}{\mathbf{F}}
\newcommand{\bfG}{\mathbf{G}}
\newcommand{\bfP}{\mathbf{P}}
\newcommand{\bfR}{\mathbf{R}}
\newcommand{\bfL}{\mathbf{L}}

\newcommand{\calD}{\mathcal{D}}
\newcommand{\calE}{\mathcal{E}}
\newcommand{\calF}{\mathcal{F}}

\newcommand{\calO}{\mathcal{O}}
\newcommand{\calU}{\mathcal{U}}
\newcommand{\calV}{\mathcal{V}}
\newcommand{\calZ}{\mathcal{Z}}

\newcommand{\frako}{\mathfrak{o}}
\newcommand{\frakH}{\mathfrak{H}}
\newcommand{\frakM}{\mathfrak{M}}
\newcommand{\frakS}{\mathfrak{S}}

\newcommand{\ch}{\operatorname{ch}}
\newcommand{\id}{\operatorname{id}}
\newcommand{\rk}{\operatorname{rk}}
\newcommand{\tr}{\operatorname{tr}}

\newcommand{\Coker}{\operatorname{Coker}}
\newcommand{\Img}  {\operatorname{Im}}
\newcommand{\Ker}  {\operatorname{Ker}}

\newcommand{\End}{\operatorname{End}}
\newcommand{\Hom}{\operatorname{Hom}}
\newcommand{\Ext}{\operatorname{Ext}}

\newcommand{\Coh} {\operatorname{\rm Coh}}

\newcommand{\seteq} {\mathbin{:=}}

\newcommand{\simto} {\xrightarrow{\sim}}

\newcommand{\GL}   {\operatorname{GL}}
\newcommand{\SL}   {\operatorname{SL}}
\newcommand{\SO}   {\operatorname{SO}}
\newcommand{\PSL}  {\operatorname{PSL}}
\newcommand{\LieO} {\operatorname{O}}
\newcommand{\Liesl}{\operatorname{sl}}
\newcommand{\Ad}   {\operatorname{Ad}}

\newcommand{\Eq} {\operatorname{Eq}}
\newcommand{\FM} {\operatorname{FM}}
\newcommand{\NS} {\operatorname{NS}}
\newcommand{\Sym}{\operatorname{Sym}}
\newcommand{\WIT}{\operatorname{WIT}}
\newcommand{\alg}{{\rm alg}}
\newcommand{\ev} {{\rm ev}}
\newcommand{\op} {{\rm op}}
\newcommand{\sgn}{{\rm sgn}}

\newcommand{\Pic}    [1]{{\rm Pic}(#1)}
\newcommand{\Piczero}[1]{{\rm Pic}^0(#1)}
\newcommand{\Hilb}   [2]{\operatorname{\rm Hilb}^{#1}(#2)}
\newcommand{\mpr}    [1]{\langle #1\rangle}

\begin{document}

\title[Semi-homogeneous sheaves, FMT and moduli on abelian surfaces]
{Semi-homogeneous sheaves, Fourier-Mukai transforms and \\
moduli of stable sheaves on abelian surfaces}
\date{June 25, 2009}

\author{Shintarou Yanagida}
\address{
{\rm Shintarou Yanagida}
\newline
\indent
Department of Mathematics, Faculty of Science, Kobe University,
Kobe 657-8501, Japan}
\email{yanagida@math.kobe-u.ac.jp}
\thanks{The first author is supported by JSPS Research Fellowships for Young Scientists (No.\ 21$\cdot$2241).}

\author{K\={o}ta Yoshioka}
\address{
{\rm K\={o}ta Yoshioka}
\newline
\indent
Department of Mathematics, Faculty of Science, Kobe University,
Kobe 657-8501, Japan
\newline
\indent 
Max-Planck-Institut f\"{u}r Mathematik,
Vivatsgasse 7
53111 Bonn}
\email{yoshioka@math.kobe-u.ac.jp}
\thanks{The second author is supported by the Grant-in-aid for Scientific
Research (No.\ 18340010), JSPS}


\begin{abstract}
This paper studies stable sheaves on abelian surfaces of Picard number one. Our main tools are semi-homogeneous sheaves and Fourier-Mukai transforms. We introduce the notion of semi-homogeneous presentation and investigate the behavior of stable sheaves under Fourier-Mukai transforms. As a consequence, an affirmative proof is given to the conjecture proposed by Mukai in the 1980s. This paper also includes an explicit description of the birational correspondence between the moduli spaces of stable sheaves and the Hilbert schemes.
\end{abstract}

\maketitle
\tableofcontents

\section{Introduction}

This paper analyzes stable sheaves on abelian surfaces using two notions. The first is the {\it semi-homogeneous sheaves} and the second is the {\it Fourier-Mukai transforms}. Both of them were introduced by Mukai~\cite{Mukai:1978,Mukai:1981}.

Semi-homogeneous sheaves are well-behaved coherent sheaves on abelian surfaces. They are semi-stable for any ample divisor, and their cohomological properties and classifications were fully investigated in \cite{Mukai:1978}. Semi-homogeneous sheaves may be regarded as building blocks of more complicated stable sheaves.

We consider presenting stable sheaves as kernel or cokernel of morphisms between semi-homogeneous sheaves. These presentations will be called \emph{semi-homogeneous presentations}.
\begin{dfn*}[{Definition \ref{defn:semihom-presen}}]
For $E\in\Coh(X)$, a \emph{semi-homogeneous presentation} is an exact sequence
\begin{align*}
\begin{matrix}
0\to  E   \to  E_1 \to  E_2 \to  0 \quad \mbox{or}\quad
0\to  E_1 \to  E_2 \to  E \to  0
\end{matrix}
\end{align*}
with the condition $v(E_1)=\ell_1v_1$, $v(E_2)=\ell_2v_2$ and
\begin{align*}
\begin{split}
&\ell_1, \ell_2\in\bbZ_{>0},\quad
 (\ell_1-1)(\ell_2-1)=0,\\
&\mpr{v_1,v_2}=-1, \quad 
 v_1, v_2\colon \mbox{positive, primitive and isotropic Mukai vectors}.
\end{split}
\end{align*}
\end{dfn*}
The uniqueness of the presentation and the conditions for its existence  will be clarified in \S\,\ref{section:semihom-presen} under the condition that the Picard number of the surface is one. In the main Theorem~\ref{thm:presentation}, we show that a numerical condition on Chern characteristics (or Mukai vectors) of stable sheaves totally controls the existence of the semi-homogeneous presentations. The condition is described by an indefinite equation, called \emph{numerical equation} (see Definition \ref{defn:semihom-presen}). This equation plays an important role in this paper.

The proof of Theorem~\ref{thm:presentation} requires moduli space $\frakM^{+}(v_1,v_2,\ell_1,\ell_2)$ of simple two-term complexes associated to the solution $(v_1,v_2,\ell_1,\ell_2)$ of the numerical equation. We construct these spaces in \S\,\ref{sect:moduli}. These spaces are essentially the Hilbert schemes of points over the abelian surfaces and enable us to treat the family of two-term complexes including the desired semi-homogeneous presentations.

As an application of the semi-homogeneous presentation, we prove the Conjecture~\ref{conj:mukai1'}, originally proposed by Mukai~\cite{Mukai:1980}:
\begin{conj*}[{Conjecture~\ref{conj:mukai1'}}]
Let $X$ be an abelian surface with $\NS(X)=\bbZ H$ and $v$ be a Mukai vector with $\ell\seteq \mpr{v^2}/2>0$. Suppose that $v$ has at least one solution of the numerical equation \eqref{eq:numerical_equation}.
\begin{enumerate}
\item
Among the numerical solutions for $v$, take the solution $(v_1,v_2,\ell_1,\ell_2)$ such that $\rk(v_i)$ is minimum, where the index $i\in\{1,2\}$ is determined by $\ell_i=\ell$. Then for a general member $\bigl[E^{-1}\xrightarrow{f} E^0\bigr]$ of $\frakM^{+}(v_1,v_2,\ell_1,\ell_2)$, $f$ is either surjective or injective. 
\item
In the situation of (1), the kernel or cokernel of $f$ is stable.
\item
A general member of $M_X^H(v)$ has a semi-homogeneous presentation corresponding to the numerical solution $(v_1,v_2,\ell_1,\ell_2)$ of $v$ such that $\rk(v_i)$ is minimum, where the index $i\in\{1,2\}$ is determined by $\ell_i=\ell$.
\end{enumerate}
\end{conj*}
\begin{thm*}[{Theorem~\ref{thm:mukai-conj1'}}]
The above conjecture is true.
\end{thm*}
We also obtained some explicit conditions when sheaves are transformed into sheaves under Fourier-Mukai transforms. These will be explained in \S\,\ref{section:appli}.

In order to see how semi-homogeneous presentations behave under a Fourier-Mukai transform, we need to calculate Mukai vectors explicitly. One can describe Mukai vectors economically using quadratic matrices, which will be explained in \S\,\ref{section:matrix}. These matrices turn out to form an arithmetic group $G$, introduced in Definition~\ref{dfn:G}. During the study of this \emph{matrix description}, we are led to the study of cohomological representation of the equivalences among the derived categories of coherent sheaves on abelian surfaces. The result is shown in Theorem~\ref{thm:bij}.

In the subsequent \S\,\ref{section:complex}, we extensively use the semi-homogeneous presentation and the matrix description to analyze the behavior of stable sheaves under Fourier-Mukai transforms. The resulting image is in general a complex, and the analysis is too complicated. We circumvent this difficulty by constructing a sequence of \emph{tame complexes}, for which one can completely describe the cohomology sheaves and the behavior under Fourier-Mukai transforms. As a consequence we obtained Theorems~\ref{thm:csv:1} and \ref{thm:csv:2}, which give the numerical condition under which a stable sheaf is transformed into a sheaf by Fourier-Mukai transforms.

In \S\,\ref{sect:birational}, we present an explicit construction of the birational correspondence between the moduli spaces of stable sheaves and the Hilbert schemes. The main result is shown in Proposition~\ref{prop:biraional}. The birational correspondence requires the existence of a solution of the numerical solution and the resulting birational map is realized by a Fourier-Mukai transform. In the final \S\S\,\ref{subsec:pp}, we treat the case of the principally polarized abelian surfaces.  In this case the arithmetic group $G$ introduced in \S\,\ref{section:matrix} reduces to $\SL(2,\bbZ)$, and the condition above reduces to the existence of solutions of integral quadratic forms. Thus using the classical theory of quadratic forms we can obtain a clear description of the birational maps between moduli and the Hilbert schemes, as shown in 
\begin{thm*}[{Theorem~\ref{thm:class-number}.}]
Let $X$ be a principally polarized abelian surface with $\NS(X)=\bbZ H$.
Let $v=(r,d H,a)$ be a Mukai vector satisfying the following condition.
\begin{align*}
\begin{split}
&\ell\seteq\mpr{v^2}/2\in\bbZ_{>0} \mbox{ is not a square number, and }
\\
&\mbox{the class number of quadratic forms with discriminant $\ell$ is 1.}
\end{split}
\end{align*}
Then the birational morphism $M_X^H(v)\cdots\to X\times{\rm Hilb}^{\ell}X$ is given by the following description.
\begin{enumerate}
\item
Take a Mukai vector $v_1=(p_1^2,p_1q_1H,q_1^2)$ which is of minimum rank among those satisfying $2p_1q_1d-p_1^2a-q_1^2r=\epsilon$, where $\epsilon=1$ or $-1$.  Then 
\begin{align*}
\gamma\seteq\pm
\begin{bmatrix}
 q_1&\epsilon(d q_1-a p_1)\\-p_1&\epsilon(d p_1-r q_1)
\end{bmatrix}
\in\PSL(2,\bbZ)
\end{align*}
diagonalizes the matrix $Q_v=\begin{bmatrix}r&d\\d&a\end{bmatrix}$, that is, 
${}^t\gamma Q_v\gamma=-\epsilon\begin{bmatrix}1&0\\0&-\ell\end{bmatrix}$.

\item
Set $v_2\seteq(p_2^2,p_2q_2H,q_2^2)$, where $q_2=-\epsilon(d q_1-a p_1)$ and $p_2=\epsilon(d p_1-r q_1)$. Then a general member of $M_X^H(v)$ has a semi-homogeneous presentation. 
Moreover the FMT $\Phi\seteq\Phi_{X\to X}^{\bfE }$ 
such that $\theta(\Phi)=\gamma$ or the composition $\calD_X\circ\Phi$ 
gives the birational correspondence $M_X^H(v)\cdots\to X\times\Hilb{\ell}{X}$ 
up to shift.
\end{enumerate}
\end{thm*}

\begin{ack}
In the final stage for preparing this article,
the second named author stayed at Max Planck Institut f\"{u}r Mathematik.
He would like to thank the institution for
the hospitality and support very much.
\end{ack}

\begin{nota}
All schemes are of finite type over $\bbC$. For coherent sheaves $E$ and $F$ on a scheme $X$, we abbreviate $\Hom_{\calO_X}(E,F)$ and $\Ext^i_{\calO_X}(E,F)$ to $\Hom(E,F)$ and $\Ext^i(E,F)$.

The abelian category of coherent sheaves on a projective variety $X$ is denoted by $\Coh(X)$, and the bounded derived category of $\Coh(X)$ is denoted by $\bfD(X)$. An object of $\bfD(X)$, namely a complex of coherent sheaves on $X$, will be denoted like $E^\bullet$. We also use the symbol $[E^{-1}\to E^0]$ for a two-term complex. In this notation, the $(-1)$-th and the zeroth component of the complex are $E^{-1}$ and $E^0$, and the other components are the zero sheaf. For an object $E^\bullet\in\bfD(X)$, $H^i(E^\bullet)$ means the $i$-th cohomology sheaf of $E^\bullet$. We denote $\Hom_{\bfD(X)}(E^\bullet,F^\bullet)$ and $\Ext^i_{\bfD(X)}(E^\bullet,F^\bullet)$ by $\Hom(E^\bullet,F^\bullet)$ and $\Ext^i(E^\bullet,F^\bullet)$ if confusion does not arise. We also denote the opposite category of $\bfD(X)$ by $\bfD(X)_{\op}$ and the dualizing functor by
\begin{align*}
\calD_X\colon\bfD(X)\to\bfD(X)_{\op},\quad 
x\mapsto x^\vee\seteq\bfR\mathcal{H}om(x,\calO_X).
\end{align*}

Two smooth projective varieties $X$ and $Y$ are said to be Fourier-Mukai partners if there is an equivalence $\bfD(X)\simeq\bfD(Y)$. The set of isomorphism classes of Fourier-Mukai partners of $X$ is denoted by $\FM(X)$. The set of equivalences between $\bfD(X)$ and $\bfD(Y)$ is denoted by $\Eq(\bfD(X),\bfD(Y))$.

For an abelian variety $X$, 
its dual variety $\Piczero{X}$ is denoted as $\widehat{X}$. 
The Poincar\'{e} line bundle of $X$ is denoted by $\bfP_X$ 
and we sometimes abbreviate it to $\bfP$ if confusion does not occur.

Assume that an ample divisor $H$ is fixed on a projective variety $X$. For $E\in\Coh(X)$, we define its slope by $\mu(E)\seteq (c_1(E),H)/\rk(E)$ if $\rk(E)>0$ . Here $(\cdot,\cdot)$ is the intersection form on $H^2(X,\bbQ)$. If $\rk(E)=0$, then we define $\mu(E)\seteq+\infty$.

In this paper we use the word `stability' in the sense of Simpson (for detail, see \cite{HuybrechtsLehn:book}). $M_X^H(v)^{ss}$ denotes the coarse moduli scheme which parametrizes $S$-equivalent classes of semi-stable shaves $E$ on $X$ having $v$ as their Mukai vectors (see \S\S\,\ref{subsec:ml}  for the definition). We denote by $M_X^H(v)$ the open subscheme consisting of stable sheaves.
\end{nota}

\section{Preliminaries}\label{section:prelim}

\subsection{Mukai lattice and Mukai vector}\label{subsec:ml}

First we introduce the \emph{Mukai lattice} \cite{Mukai:1987:tata} for an abelian surface $X$. Consider the lattice $(H^{\rm ev}(X,\bbZ),\mpr{\cdot,\cdot})$ consisting of the even part of the cohomology group
\begin{align*}
H^{\rm ev}(X,\bbZ)\seteq\bigoplus_{i=0}^2 H^{2i}(X,\bbZ),
\end{align*}
and the \emph{Mukai pairing} $\mpr{\cdot,\cdot}$ defined as
\begin{align*}
\mpr{x,y}\seteq
&\int_X(x_1\cup y_1-x_0\cup y_2-x_2\cup y_0),
\\
&x=(x_0,x_1,x_2),y=(y_0,y_1,y_2),\ x_i,y_i\in H^{2i}(X,\bbZ).
\end{align*} 
This paring is symmetric and bilinear. We also have
$(H^{\rm ev}(X,\bbZ),\mpr{ \cdot,\cdot})\cong U^{\perp 4}$, where $U$ is the hyperbolic lattice and $\perp$ is the orthogonal direct sum. In particular, $(H^{\rm ev}(X,\bbZ),\mpr{ \cdot,\cdot})$ is a non-degenerate even lattice.

For an object $E^\bullet\in \bfD(X)$, the \emph{Mukai vector} of $E^\bullet$ is defined as its Chern character. 
\begin{align*}
v(E^\bullet)\seteq
\ch(E^\bullet)=(\rk(E^\bullet),c_1(E^\bullet),\chi(E^\bullet))
\in H^{\ev}(X,\bbZ).
\end{align*}
By the Grothendieck-Riemann-Roch formula we have 
\begin{align*}
\mpr{v(E^\bullet),v(F^\bullet)}=-\chi(E^\bullet,F^\bullet)
\end{align*}
for $E^\bullet, F^\bullet\in\bfD(X)$. For a class $e\in K(X)$ of the Grothendieck group of $X$, we similarly define $v(e)\seteq\ch(e)=(\rk(e),c_1(e),\chi(e))$.

We call an element $v\in H^{\ev}(X,\bbZ)$ a \emph{Mukai vector} if there exists an object $E^\bullet\in \bfD(X)$ such that $v=v(E^\bullet)$. The set of Mukai vectors is nothing but the submodule $H^{\ev}(X,\bbZ)_{\alg}$ of $H^{\ev}(X,\bbZ)$ defined by
\begin{align*}
H^{\ev}(X,\bbZ)_{\alg}
&\seteq
H^0(X,\bbZ)\oplus\NS(X)\oplus H^4(X,\bbZ).
\end{align*}
Here $\NS(X)$ denotes the N\'{e}ron-Severi group of $X$. $(H^{\ev}(X,\bbZ)_{\alg},\mpr{\cdot,\cdot})$ is called the \emph{Mukai lattice} of $X$.

We list several definitions concerning Mukai vectors.

\begin{dfn}
Let $v=(r,\xi,a)$ be a Mukai vector. $v$ is called \emph{isotropic} if $\mpr{v^2}\seteq\mpr{v,v}=0$. $v$ is called \emph{positive} and denoted by $v>0$ if $r>0$, or if $r=0$ and $\xi$ is effective, or if $r=0$ and $\xi=0$ and $a>0$. $v$ is called \emph{primitive} if $H^{\ev}(X,\bbZ)/(\bbZ v)$ is torsion free. For a non-zero Mukai vector $v$, its \emph{multiplicity} $\ell$ is defined to be the non-negative integer defined by the relation $v=\ell w$, where $w$ is a primitive Mukai vector. We also define $\rk(v)\seteq r$, $c_1(v)\seteq\xi$ and $\chi(v)\seteq a$. When an ample divisor $H$ is fixed, we define $\mu(v)\seteq(\xi,H)/r$ if $r>0$ and $\mu(v)\seteq +\infty$ if $r=0$. At last, we define $v^\vee\seteq(r,-\xi,a)$.
\end{dfn}

In this article we sometimes denote a Mukai vector $v=(r,\xi,a)$ as $v=r+\xi+a\rho_X$ additively. Here $\rho_X$ is the element of $H^4(X,\bbZ)$ satisfying $\int_X\rho_X=1$. Also we will use the exponential form of the Chern characteristics. For example, for $d\in \bbZ$ and an ample divisor $H$, $e^{d H}$ means $(1,d H, d^2(H^2)/2)$.

\subsection{Moduli spaces of stable sheaves}\label{subsec:ms}

We mention only a few properties of the moduli of stable sheaves, although there are considerable amount of study.

\begin{dfn}\label{defn:divisor:general}
Let $v$ be a Mukai vector. An ample divisor $H$ is \emph{general} with respect to $v$, if the following condition holds. For any $\mu$-semi-stable sheaf $E$ with $v(E)=v$, if $F\subset E$ satisfies $(c_1(F),H)/\rk(F)=(c_1(E),H)/\rk(E)$, then $c_1(F)/\rk(F)=c_1(E)/\rk(E)$.
\end{dfn}

\begin{fact}\label{fact:moduli}
Let $v=(r,\xi,a)$ be a Mukai vector of $r>0$ and $\mpr{v^2}>0$, and $H$ be an ample divisor which is  general  with respect to $v$. Then $M_X^H(v)$ is smooth, projective and irreducible.
\end{fact}
\begin{proof}
Smoothness and projectivity are proved in \cite{Mukai:1984}. 
Irreducibility is proved in \cite[Theorem 3.18]{Yoshioka:2003:twistI}.
\end{proof}

\subsection{Fourier-Mukai transforms on abelian surfaces}

The next ingredient is the Fourier-Mukai transform initiated by Mukai in \cite{Mukai:1981}. In this subsection we summarize its main properties, focusing on the case of abelian surfaces. For a general treatment of this functor we recommend the reader to see \cite{Huybrechts:book}.

\begin{dfn}
Let $X$ and $Y$ be smooth projective varieties. For an object $\bfE^\bullet\in \bfD(X\times Y)$, the integral functor $\Phi^{\bfE^\bullet}_{X\to Y}$ is defined as
\begin{align*}
\begin{array}{c c c c}
\Phi^{\bfE^\bullet}_{X\to Y}\colon
&\bfD(X)&\to &\bfD(Y)\\
&x      &\mapsto 
        &\bfR p_{Y*}(\bfE^\bullet\stackrel{\bfL}{\otimes} p_X^*(x)),
\end{array}
\end{align*}
where $p_X\colon X\times Y\to X$ and $p_Y\colon X\times Y\to Y$ are the natural projections. We call $\bfE^\bullet$ its \emph{kernel}. If $\Phi^{\bfE^\bullet}_{X\to Y}$ is an exact equivalence, it is called a \emph{Fourier-Mukai transform}. In this paper  we abbreviate it to FMT.
\end{dfn}

In order to describe FMT in detail, it is convenient to introduce the next terminology.

\begin{dfn}
Let $X$ and $Y$ be smooth projective varieties, $E$ be a coherent sheaf over $X$ and $\Phi\colon\bfD(X)\to\bfD(Y)$ be an equivalence. We say that $E$ satisfies  \emph{the weak index theorem} with respect to $\Phi$ and that \emph{its index is equal to} $i$ if $H^j(\Phi(E))=0$ for every $j\neq i$. If this condition holds, then we also say that ${\rm WIT}_i$ holds for $E$ with respect to $\Phi$.
\end{dfn}

We will extensively utilize the following classification of FMTs on abelian surfaces due to Orlov \cite{Orlov:2002}. For the proof, see \cite[Theorem 1.4]{Yoshioka:2009:ann}.

\begin{fact}\label{fact:FMT-abelian-surface}

\begin{enumerate}
\item
Let $X$ be an abelian surface and $Y$ be a smooth projective variety. If $\Phi^{\bfE^\bullet}_{X\to Y}$ is a FMT, then $Y$ is an abelian surface and the inverse functor is given by $\Phi_{Y\to X}^{\bfE^{\bullet \vee}[2]}$.
\item
Let $X$ be an abelian surface and $v=(r,\xi,a)$ be a primitive isotropic Mukai vector. Suppose that there exists a universal family $\bfE$ on $M_X^H(v)\times X$. Then $\Phi_{X\to M_X^H(v)}^\bfE$ is a FMT.
\item
Let $X$ and $Y$ be abelian surfaces. Suppose that there exists an equivalence $\Phi\colon \bfD(X)\to \bfD(Y)$. Then there exists a primitive isotropic Mukai vector $v$ such that $Y=M_X^H(v)$. Moreover there exists a universal family $\bfE$ on $Y\times X$ and an integer $k\in\bbZ$ such that $\Phi=\Phi_{X\to Y}^{\bfE[k]}$.
\end{enumerate}
\end{fact}

Let $\Phi=\Phi_{X\to Y}^{\bfE^\bullet}$ be a FMT 
between abelian surfaces $X$ and $Y$. 
This functor induces an isomorphism 
$\Phi^K\colon K(X)\to K(Y)$ 
on the Grothendieck groups and an isomorphism 
$\Phi^H\colon H^{\ev}(X,\bbQ)\to H^{\ev}(Y,\bbQ)$ 
on the cohomology groups. 
Moreover, by the Grothendieck-Riemann-Roch formula, 
the following diagram is commutative. 
(For the proof, see \cite[\S\S\,5.2]{Huybrechts:book}.)
\begin{align}\label{diag:comm1}
\xymatrix{
 \ar@{}[rd]|{\circlearrowright}
 \bfD(X) 
 \ar[r]^{\Phi}    
 \ar[d]_{[\ ]} 
&\bfD(Y)
 \ar[d]^{[\ ]} 
\\
 \ar@{}[rd]|{\circlearrowright}
 K(X)    
 \ar[r]^{\Phi^K}
 \ar[d]_{v}
&K(Y)
 \ar[d]^{v}
\\
 H^{\ev}(X,\bbQ)
 \ar[r]_{\Phi^H}
&H^{\ev}(Y,\bbQ)
}
\end{align}
Here the symbol $[\ ]$ means taking the class of the object and $v$ means taking the Mukai vector.

Moreover, one can prove that 
$\Phi^H(x)\in H^{\ev}(Y,\bbZ)$ 
for $x\in H^{\ev}(X,\bbZ)$ and that 
$\Phi^H$ induces isometries 
$(H^{\ev}(X,\bbZ),\mpr{\cdot,\cdot})\to(H^{\ev}(Y,\bbZ),\mpr{\cdot,\cdot})$ 
and 
$(H^{\ev}(X,\bbZ)_\alg,\mpr{\cdot,\cdot}) \to
 (H^{\ev}(Y,\bbZ)_\alg,\mpr{\cdot,\cdot})$. 
(For the proof, see \cite[Lemma 2.2]{Yoshioka:2001:ann}.)
$\Phi^H$ will be called the \emph{cohomological FMT} 
induced by $\Phi$. 
We often suppress the symbol $H$ of the cohomological FMT 
and simply denote it by $\Phi$, using the same notation of the original FMT.

\begin{lem}\label{lem:picard}
Let $X$ be an abelian surface with $\NS(X)=\bbZ H$, $(H^2)=2n$ and let us take $Y\in\FM(X)$. Then $\NS(Y)$ is of rank one. Moreover if we put $\NS(Y)=\bbZ \widehat{H}$, where $\widehat{H}$ is the ample generator of $\NS(Y)$, then we have $(\widehat{H}^2)=2n$.
\end{lem}

\begin{proof}
A FMT $\bfD(X)\simto\bfD(Y)$ induces an isometry $H^{\ev}(X,\bbZ)_{\alg}\simto H^{\ev}(Y,\bbZ)_\alg$. Since in our case $H^{\ev}(X,\bbZ)_\alg\cong \NS(X) \perp U$ and $H^{\ev}(Y,\bbZ)_\alg\cong \NS(Y) \perp U$, where $U$ is the hyperbolic lattice, we find that $\rho(X)+2=\rk(H^{\ev}(X,\bbZ)_\alg)=\rk(H^{\ev}(Y,\bbZ)_\alg)=\rho(Y)+2$ and $(H^2)=-\det H^{\ev}(X,\bbZ)_\alg=-\det H^{\ev}(Y,\bbZ)_\alg=(\widehat{H}^2)$. Therefore the conclusion holds.
\end{proof}

Finally, the Serre duality yields the following lemma.

\begin{lem}\label{lem:Serre}
For a FMT $\Phi_{X \to Y}^{\bfE^\bullet}\colon\bfD(X)\to\bfD(Y)$ 
between abelian surfaces $X$ and $Y$, we have 
\begin{align*}
\calD_{Y}\Phi_{X \to Y}^{{\bfE}^{\bullet\vee}}
=\Phi_{X \to Y}^{\bfE^\bullet[2]} \calD_{X}.
\end{align*}
\end{lem}

\subsection{Semi-homogeneous sheaves}\label{subsection:semi-homog}

Here we recall the notion of semi-homogeneous sheaves over abelian surfaces, 
which was introduced by Mukai \cite{Mukai:1978}. 
Fix an abelian surface $X$.

\begin{dfn}
For any coherent sheaf $E$ on $X$, the subset 
\begin{align*}
S(X)\seteq \{(x,\hat{x})\in X\times \widehat{X}
 \mid T_x^*(E)\otimes\bfP_{\hat{x}}\cong E\}
\end{align*}
is an abelian subvariety of $X\times\widehat{X}$ with $\dim S(X)\le 2$ 
(see \cite[Proposition 3.3]{Mukai:1978}). 
If $\dim S(X)=2$, then $E$ is called \emph{semi-homogeneous}.
\end{dfn}

Semi-homogeneous sheaves are classified as below. 
For the proof, see \cite{Mukai:1978} and \cite[\S\,4]{Yoshioka:2009:ann}. 

\begin{fact}\label{fact:semihom}
Let $H$ be an ample divisor on $X$ and 
$v=(r,\xi,a)$ be a positive isotropic Mukai vector.
\begin{enumerate}
\item
Let $E$ be a coherent sheaf with $v(E)=v$. Fix an ample divisor $H$ on $X$.
\begin{enumerate}
\item
$E$ is simple $\iff$ $E$ is stable $\iff$ $v$ is primitive.
The stability of $E$ is independent of $H$.

\item
If $v$ is primitive and $E$ is semi-stable with respect to $H$ and $v(E)=n v$, 
then $E$ is $S$-equivalent to $\bigoplus_{i=1}^n E_i$, $E_i\in M_X^H(v)$. 
The semi-stability of $E$ is independent of $H$.

\item
Every semi-stable sheaf with isotropic Mukai vector is semi-homogeneous.

\end{enumerate}
\item
Stable semi-homogeneous sheaves are classified as follows.

\begin{enumerate}
\item
Semi-homogeneous vector bundle.

For any $\delta\in\NS(X)\otimes_\bbZ\bbQ$, there is a simple semi-homogeneous vector bundle $E$ satisfying $c_1(E)/\rk(E)=\delta$.

If $E$ and $E'$ are simple semi-homogeneous vector bundles with $c_1(E)/\rk(E)=c_1(E')/\rk(E')$, then $E\cong E'\otimes M$ for some $M\in{\rm Pic}^0(X)$.

For any semi-homogeneous vector bundle $E$, there exist an isogeny $\pi\colon X'\to X$ of abelian surfaces and a line bundle $L$ on $X'$ such that $\pi_*(L)\cong E$.

\item
Torsion semi-homogeneous sheaves.

${\rm Supp}(E)$ is an elliptic curve $C$ and 
$E$ is a locally free semi-stable sheaf over $C$, 

 or ${\rm Supp}(E)$ is one point and $E$ is a skyscraper sheaf.
\end{enumerate}
\end{enumerate}
\end{fact}

\begin{rmk}\label{rmk:semihom-ns-rank1}
Notice that if $\NS(X)=\bbZ H$ then the first case in (2) (b) does not occur. Therefore non-locally-free semi-homogeneous sheaves are nothing but coherent sheaves with finite support.
\end{rmk}

The cohomology groups of semi-homogeneous sheaves were also investigated by Mukai~\cite{Mukai:1978}. The results are as follows. For the proof, see Mukai's paper or \cite[Proposition 4.4]{Yoshioka:2009:ann}.

\begin{fact}\label{fact:semihom-ext}
Let $E$ and $F$ be semi-homogeneous sheaves.
\begin{enumerate}
\item
Suppose that $\rk(E)>0$ and $\rk(F)>0$.
\begin{enumerate}
\item
If $\mpr{v(E),v(F)}>0$, then $\Hom(E,F)=\Ext^2(E,F)=0$.

\item
If $\mpr{v(E),v(F)}<0$, then $\mu(E)\neq\mu(F)$, $\Ext^1(E,F)=0$, and
\begin{align*}
\begin{cases}
\Hom(E,F)=0   &\mbox{if}\quad \mu(E)>\mu(F),\\ 
\Ext^2(E,F)=0 &\mbox{if}\quad \mu(E)<\mu(F),\end{cases}
\end{align*}
\end{enumerate}
\item
Suppose that $\rk(E)>0$ and $F$ is a torsion sheaf.
\begin{enumerate}
\item
If $\mpr{v(E),v(F)}>0$, then $\Hom(E,F)=\Ext^2(E,F)=0$.

\item
If $\mpr{v(E),v(F)}<0$, then $\Ext^1(E,F)=\Ext^2(E,F)=0$.
\end{enumerate}
\item
Suppose $E$ and $F$ are torsion sheaves. Then we have $\mpr{v(E),v(F)}\ge0$. 
If $\mpr{v(E),v(F)}>0$, then $\Hom(E,F)=\Ext^2(E,F)=0$.
\end{enumerate}
\end{fact}

\begin{lem}\label{lem:rho=1}
Let $E$ and $F$ are semi-homogeneous sheaves on $X$.
Assume that $\NS(X)= \bbZ H$.
Then $\mpr{v(E),v(F)} \le 0$.
If the equality holds, then $\bbQ v(E)=\bbQ v(F)$.
\end{lem}

\begin{proof}
If $\rk(E)>0$ and $\rk(F)>0$, 
then one can put 
\begin{align}
\label{eq:mv}
v(E)=r_1 e^{(d_1 /r_1)H}
    =\bigl(r_1,d_1 H, d_1^2(H^2)/(2 r_1)\bigr),
\quad  
v(F)=r_2 e^{(d_2/r_2)H}
    =\bigl(r_2,d_2 H,d_2^2(H^2)/(2 r_2)\bigr),
\end{align}
with $r_1,r_2\in\bbZ_{>0}$ and $d_1,d_2\in\bbZ$. 
The Mukai paring becomes 
$\mpr{v(E),v(F)}=-(r_1 d_2-r_2 d_1)^2(H^2)/(2 r_1 r_2)\le 0$, 
which yields the first statement. 
The equality holds if and only if $d_1/r_1=d_2/r_2$. 
This is equivalent to $\bbQ v(E)=\bbQ v(F)$ by \eqref{eq:mv}.
Thus we have the second statement.

Assume that $\rk(E)=0$. 
By Remark~\ref{rmk:semihom-ns-rank1} 
one can put $v(E)=(0,0,a)$ with $a\in\bbZ_{>0}$. 
Then $\mpr{v(E),v(F)}=-a \rk(F)\le 0$. 
The equality holds if and only if $\rk(F)=0$ if and only if $v(F)=(0,0,b)$ for some $b\in\bbZ_{>0}$. 
Therefore we have the conclusion.

The case $\rk(F)=0$ is the same as the case $\rk(E)=0$.
\end{proof}

\begin{prop}\label{prop:semihom-fmt}
Assume that $\NS(X)=\bbZ H$. 
Let $v$ be a positive, primitive and isotropic Mukai vector 
and let $Y\seteq M_X^H(v)$ be the moduli space. 
Suppose that there exists a Mukai vector $w$ such that $\mpr{v,w}=-1$. 
Then, there exists a universal family $\bfE$ on $Y\times X$ 
by \cite[Proposition A.6]{Mukai:1987:tata}. 
Now set $\Phi\seteq\Phi_{X\to Y}^{\bfE^\vee}$ and $\mu\seteq\mu(v)$. 
Under these assumptions, the next statements hold.

\begin{enumerate}
\item
For any semi-homogeneous sheaf $E$ on $X$, $\WIT_i$ holds w.r.t. $\Phi$ 
and the index $i$ is 
\begin{align*}
i=\begin{cases} 0&\mu(E)>\mu,\\ 2&\mu(E)\le\mu. \end{cases}
\end{align*}
Moreover the image $\Phi^i(E)$ is a semi-homogeneous sheaf. 

\item
For a semi-homogeneous sheaf $E$ with $v(E)=v$, we have $\Phi(E)=\bbC_y[-2]$. 
Here $y\in Y$ is a closed point corresponding to $E$. 

\item
Suppose that the Mukai vector $w$ in the assumption is isotropic. 
Then for a semi-homogeneous sheaf $F$ with $v(F)=w$, we have $\Phi(F)=L[-i]$, 
where $L$ is a line bundle on $Y$ and the index $i$ is $0$ or $2$ 
according to $\mu(F)>\mu$ or $\mu(F)\le\mu$.

\item
Under the additional condition $\mpr{w^2}=0$, 
let $E$ and $F$ be coherent sheaves with $v(E)=v$ and $v(F)=w$ 
as in (2) and (3). 
Assume further that $F$ is stable. 
Then there is a coherent sheaf $\bfF\in\Coh(X\times Y)$ 
such that $\Phi'\seteq\Phi_{X\to Y}^{\bfF^\vee}$ satisfies 
$\Phi'(E)=\bbC_y[-2]$ and $\Phi'(F)=\calO_Y[-i]$.
\end{enumerate}
\end{prop}

\begin{proof}
(1)
Since $\Phi(E)=\bfR\Hom_{p_Y}(\bfE,p_X^*(E))$, 
we have $H^i(\Phi(E))=\Ext^i_{p_Y}(\bfE,p_X^*(E))$. 
Then from the base change theorem for the relative extension 
and the cohomological result in Fact~\ref{fact:semihom-ext}, 
WIT holds for $E$ and the index coincides with the one given in the statement. 
Thus the first statement is verified.

Next, we recall that $\Phi$ induces an isometry of the Mukai lattices. 
Hence $\mpr{v(\Phi^i(E))^2}=\mpr{v^2}=0$. 
If $E$ is simple, 
then from this equality and Fact~\ref{fact:semihom} (1) (a) 
we find that $\Phi^i(E)$ is a stable semi-homogeneous sheaf. 
For a proper semi-stable semi-homogeneous sheaf $E$, 
$\Phi^i(E)$ is a also semi-homogeneous sheaf, 
since it is an extension of the stable semi-homogeneous sheaves 
with the same Mukai vectors.

(2)
It is an immediate consequence of (1).

(3)
{}From the second statement of (1), 
we know that $\Phi^i(F)$ is a semi-homogeneous sheaf. 
Its rank is calculated as follows.
\begin{align*}
 \rk(\Phi^{i}(F))
=\chi(\Phi^{i}(F),\bbC_y)
=\chi(F,\Psi^{2-i}(\bbC_y))
=\chi(F,\bfE|_{\{y\}\times X})
=-\mpr{w,v}
=1.
\end{align*}
Here we used the notation $\Psi\seteq\Phi_{Y\to X}^{\bfE}$ and the fact that $\Psi[2]$ is the inverse of $\Phi$. Thus $\Phi^{i}(F)$ is a semi-homogeneous sheaf of rank $1$, that is, a line bundle.

(4)
We use the same notation as in the proof of (3). If $L\cong\calO_Y$, then it is enough to set $\bfF\seteq\bfE$. Otherwise, by setting $\bfF\seteq\bfE\otimes_{\calO_{Y\times X}}p_Y^*L$, we have the desired conclusion.
\end{proof}

The next propositions are essential for our analysis. In fact several results of this paper, including Theorems \ref{thm:presentation}, \ref{thm:mukai-conj1'} and \ref{thm:class-number}, depend on these results.

\begin{fact}[{\cite[Theorem 3.1]{Yoshioka:2009:ann}}]
\label{fact:yoshioka-thm2-1}
Let $\Phi\seteq\Phi_{Y\to X}^\bfE\colon\bfD(Y)\to\bfD(X)$ be a FMT. We denote by $\Psi$ this functor $\Phi$ or the composition $\calD_X\circ\Phi$. Let $w$ be a primitive Mukai vector on $Y$ with $\mpr{w^2}>0$. Let $H'$ be a general ample divisor on $Y$ with respect to $w$. Set $v\seteq\pm\Psi(w)$. Let $H$ be a general ample divisor on $X$ with respect to $v$. 

If $\Psi(E)$ is not a sheaf up to shift for all $E\in M_Y^{H'}(w)$, then there is an integer $k$ such that for a general $E\in M_Y^{H'}(w)$, $\Psi(E)$ fits in an exact triangle
\begin{align*}
A_0\to\Psi(E)\to A_1[-1]\to A_0[1]
\end{align*}
in $\bfD(X)$ or $\bfD(X)_{\op}$. Here $A_0$ and $A_1$ are semi-homogeneous sheaves with 
\begin{align*}
v(A_i)=\ell_i v_i,\ \mpr{v_0^2}=\mpr{v_1^2}=0,\ (\ell_0-1)(\ell_1-1)=0.
\end{align*}
\end{fact}

\begin{rmk}
In \cite[Theorem 3.1]{Yoshioka:2009:ann} only the statement for $\Phi$ is proved, but the proof works in this general situation. 
\end{rmk}

\begin{fact}[{\cite[Lemma 3.2]{Yoshioka:2009:ann}}]
\label{fact:yoshioka-lem2-2}
Keep notations in Fact~\ref{fact:yoshioka-thm2-1}. 
\begin{enumerate}
\item
If $\Psi(E_0) \in \Coh(X)$ for an element $E_0\in M_Y^{H'}(w)$, 
then $\Psi(E)$ is stable with respect to $H$ 
for a general $E\in M_Y^{H'}(w)$.
\item
If $V^\bullet\seteq \Psi(E)$ is not a sheaf 
for all $E \in M_Y^{H'}(w)$, 
then there is an equivalence 
$\calF\colon\bfD(X)\to \bfD(X)$ or $\bfD(X)\to \bfD(X)_{\op}$ 
such that $\calF(v)=v$ 
and for a general $E\in M_Y^{H'}(w)$ 
one of the following three conditions holds.
\begin{enumerate}
\item
$\rk(H^0(\calF(V^\bullet)))+\rk(H^1(\calF(V^\bullet))) 
<\rk(H^0(V^\bullet))+\rk(H^1(V^\bullet))$, or
\item
$\deg H^1(\calF(V^\bullet))<\deg H^1(V^\bullet)$ 
if $\rk(H^1(V^\bullet))=0$, or
\item
$H^1(\calF(V^\bullet))=0$ is of 0-dimensional.
\end{enumerate}
\end{enumerate}
\end{fact}

\section{Semi-homogeneous presentation of sheaves}
\label{section:semihom-presen}

We introduce the notion of \emph{semi-homogeneous presentations} and discuss their properties. We show the uniqueness of this presentation in Proposition~\ref{prop:semihom-presen-unique} and state the numerical criteria for the existence in Theorem~\ref{thm:presentation}. Its proof requires certain moduli of simple complexes, so we defer it to the next \S\,\ref{sect:moduli}. 

In this section, $X$ always means an abelian surface.

\subsection{Definitions and basic properties}

\begin{dfn}\label{defn:semihom-presen}
For $E\in\Coh(X)$, a \emph{semi-homogeneous presentation} is an exact sequence
\begin{align*}
\begin{matrix}
0\to  E   \to  E_1 \to  E_2 \to  0 \quad \mbox{or}\quad
0\to  E_1 \to  E_2 \to  E \to  0
\end{matrix}
\end{align*}
with the condition $v(E_1)=\ell_1v_1$, $v(E_2)=\ell_2v_2$ and
\begin{align}
\begin{split}
\label{cond:num_eq}
&\ell_1, \ell_2\in\bbZ_{>0},\quad
 (\ell_1-1)(\ell_2-1)=0,\\
&\mpr{v_1,v_2}=-1, \quad 
 v_1, v_2\colon \mbox{positive, primitive and isotropic Mukai vectors}.
\end{split}
\end{align}
 The former sequence is called a \emph{kernel presentation}, and the latter is called a \emph{cokernel presentation}. 

For a Mukai vector $v$, the equation
\begin{align}
\label{eq:numerical_equation}
\begin{split}
&v=\pm(\ell_2 v_2-\ell_1 v_1),
\\
&\ell_1, \ell_2 \in {\bbZ}_{>0},\quad
 (\ell_1-1)(\ell_2-1)=0,\\
&v_1, v_2\colon \mbox{positive, primitive and isotropic Mukai vectors with}
\\
&\langle v_1,v_2 \rangle=-1, \quad  \mu(v_1)<\mu(v_2)
\end{split}
\end{align}
is called the \emph{numerical equation}. 
A solution $(v_1,v_2,\ell_1,\ell_2)$ 
of this equation \eqref{eq:numerical_equation} 
is called \emph{numerical solution} of $v$.
\end{dfn}

Existence of semi-homogeneous presentation is not obvious \emph{a priori}. 
We will show in Theorem \ref{thm:presentation} 
that numerical conditions on sheaves control the existence. 
The next proposition shows 
that a semi-homogeneous presentation for a stable sheaf $E$ 
guarantees the existence of a good FMT under which $E$ is transformed 
into the ideal sheaf of 0-dimensional subscheme.

\begin{prop}\label{prop:sp}
Assume $\NS(X)=\bbZ H$. Let $\ell\in\bbZ_{>0}$ and $E$ be a stable sheaf on $X$ with $\mpr{v(E)^2}=2\ell$ and write the four possible semi-homogeneous presentations for $E$ as
\begin{align*}
\begin{array}{l l l l l l l l l l}
(1)\qquad&0&\to&E_1&\to&E_2&\to&E  &\to&0,\\
(2)\qquad&0&\to&E_2&\to&E_1&\to&E  &\to&0,\\
(3)\qquad&0&\to&E  &\to&E_2&\to&E_1&\to&0,\\
(4)\qquad&0&\to&E  &\to&E_1&\to&E_2&\to&0.
\end{array}
\end{align*}
Here $E_1$ is a simple semi-homogeneous sheaf with $v(E_1)=v_1$ and $E_2$ is a semi-homogeneous sheaf with $v(E_2)=\ell v_2$, where $v_1$ and $v_2$ are primitive Mukai vectors.

Let $\bfE$ be a universal family of semi-homogeneous sheaves with Mukai vector $v_2$ and set $Y\seteq M_X^H(v_2)$, $\Phi_1\seteq\Phi^{\bfE}_{X\to Y}$, $\Phi_2\seteq\Phi^{\bfE^\vee}_{X\to Y}$. Finally let $\calD_X\colon\bfD(X)\to\bfD(X)_{\op}$ be the dualizing functor and $\widehat{\Phi}_1\seteq\Phi_1\circ\calD_X$ be the composition.

Consider the conditions
\begin{align*}
\tag{\#1$i$}
{\rm WIT}_i \text{ holds for }E \text{ w.r.t. } \widehat{\Phi}_1 
\text{ and } \widehat{\Phi}_1^i(E)\cong L\otimes I_Z,
\\[4pt]
\tag{\#2$i$}
{\rm WIT}_i\ \mbox{holds for}\  E\ \mbox{w.r.t.}\  \Phi_2
\text{ and } \Phi_2^i(E)\cong L\otimes I_Z,
\end{align*}
where $Z$ is a 0-dimensional subscheme of length $\ell$, $I_Z$ is its ideal sheaf and $L$ is a line bundle on  $Y$. Then the existence of each presentation is equivalent to the following condition.
\begin{center}
(1)$\iff$ (\#21),
\quad
(2)$\iff$ (\#12),
\quad
(3)$\iff$ (\#11),
\quad
(4)$\iff$ (\#22).
\end{center}
Moreover, the condition that a stable sheaf on $X$ has a semi-homogeneous presentation is an open condition.
\end{prop}

\begin{proof}
This is the consequence of Lemma~\ref{lem:simple-complex2},
which will be given in \S\,\ref{sect:moduli}.
\end{proof}

\begin{prop}\label{prop:semihom-presen-unique}
Suppose $\NS(X)=\bbZ H$. 
Let $E$ be a coherent sheaf on $X$ with $\mpr{v(E)^2}>0$. 
Then a kernel (or cokernel) presentation of $E$ is unique, if it exists. 
\end{prop}

More generally, we have the next statements.

\begin{lem}
Suppose $\NS(X)=\bbZ H$. 
Let $E$ be a coherent sheaf on $X$ with $\mpr{v(E)^2}>0$. 
\begin{enumerate}
\item
Suppose that we are given two exact sequences containing $E$.
\begin{align}
\xymatrix{
 0   \ar[r]    
&E   \ar[r] \ar@{=}[d] 
&E_1 \ar[r]^{\phi}
&E_2 \ar[r]
&0\ 
\\
 0   \ar[r]
&E   \ar[r]
&E_3 \ar[r]_{\psi}
&E_4 \ar[r]
&0.}
\end{align}
Here $E_i$'s are semi-homogeneous sheaves 
and neither $\phi$ nor $\psi$ is trivial. 

Then $E_1\cong E_3$, $E_2\cong E_4$ and 
the diagram extends to the next commuting diagram.
\begin{align*}
\xymatrix{
 0   \ar[r]    
&\ar@{}[rd]|{\circlearrowright}
 E   \ar[r] \ar@{=}[d] 
&\ar@{}[rd]|{\circlearrowright}
 E_1 \ar[r]^{\phi} \ar[d]^{\wr}
&E_2 \ar[r] \ar[d]^{\wr}
&0\ 
\\
 0   \ar[r]
&E   \ar[r]
&E_3 \ar[r]_{\psi}
&E_4 \ar[r]
&0.}
\end{align*}

\item
Suppose that we are given two exact sequences containing $E$.
\begin{align*}
\xymatrix{
 0   \ar[r]    
&E_1 \ar[r]^{\phi}
&E_2 \ar[r]
&E   \ar[r]  \ar@{=}[d]
&0\ 
\\
 0   \ar[r]
&E_3 \ar[r]_{\psi}
&E_4 \ar[r]
&E   \ar[r]
&0,}
\end{align*}
where $E_i$ are semi-homogeneous sheaves 
and neither $\phi$ nor $\psi$ is trivial. 
Then we have $E_1\cong E_3$ and $E_2\cong E_4$. 
We also have the next commuting diagram.
\begin{align*}
\xymatrix{
 0   \ar[r]    
&\ar@{}[rd]|{\circlearrowright}
 E_1 \ar[r]^{\phi} \ar[d]_{\wr}
&\ar@{}[rd]|{\circlearrowright}
 E_2 \ar[r] \ar[d]_{\wr} 
&E   \ar[r]  \ar@{=}[d]
&0\ 
\\
 0   \ar[r]
&E_3 \ar[r]_{\psi}
&E_4 \ar[r]
&E   \ar[r]
&0.
}
\end{align*}
\end{enumerate}
\end{lem}
\begin{proof}
This is the corollary of the next Lemma~\ref{lem:unique}.
\end{proof}

\begin{lem}
\label{lem:unique}
Assume $\NS(X)=\bbZ H$.
Let $E^{\bullet}$ be a bounded complex on $X$ such that
$\langle v(E^{\bullet})^2 \rangle>0$.  
Assume that $E^{\bullet}$ has two expressions
\begin{equation}
E^{\bullet} \cong [E_1 \overset{\phi}{\to} E_2] 
\cong [F_1 \overset{\psi}{\to} F_2],
\end{equation}
where $E_1, E_2, F_1, F_2$ are semi-homogeneous sheaves
and $\phi, \psi$ are non-trivial.
Then there are isomorphisms $f_i:E_i \cong F_i$ 
which induce the following commutative diagram in $\bfD(X)$.
\begin{align}
\label{diag:unique}
\xymatrix{
 E_1 \ar[r]^{\phi} \ar[d]^{\wr}_{f_1}
&E_2 \ar[r]^{p}    \ar[d]^{\wr}_{f_2}
&E^\bullet   \ar[r]  \ar@{=}[d]
&E_1[1]\ar[d]_{\wr}^{f_1[1]}
\\
 F_1 \ar[r]_{\psi}
&F_2 \ar[r]_{q}
&E^\bullet   \ar[r]
&F_1[1].
}
\end{align} 
\end{lem}

\begin{proof}
The proof is divided into several steps.

{\bf Step 1}.
$\mu(E_1)<\mu(E_2)$, $\mu(F_1)<\mu(F_2)$ and 
$\Ext^1(F_2,E_1)=\Ext^1(E_2,F_1)= 0$.

Note that
\begin{equation}\label{eq:unique}
v(E^{\bullet})=v(E_2)-v(E_1)=v(F_2)-v(F_1).
\end{equation}
Then we see that $0<\mpr{v(E^{\bullet})^2}=
-2\mpr{v(E_1),v(E_2)}=-2\mpr{v(F_1),v(F_2)}$.
Thus $\mpr{v(E_1),v(E_2)}=\mpr{v(F_1),v(F_2)}<0$.
We also find from \eqref{eq:unique} that
\begin{align*}
\mpr{v(E_1), v(F_2)}=
\mpr{v(E_1),v(E_2)+v(F_1)-v(E_1)}
=\mpr{v(E_1),v(E_2)}+\mpr{v(E_1),v(F_1)}.
\end{align*}
{}From Lemma \ref{lem:rho=1} and the inequality $\mpr{v(E_1),v(E_2)}<0$, 
we find that $\mpr{v(E_1), v(F_2)}<0$.
In the same way, we have $\mpr{v(F_1), v(E_2)}<0$. 
Then Fact~\ref{fact:semihom-ext} shows that  $\Ext^1(F_2,E_1)=\Ext^1(E_2,F_1)=0$.

Lemma \ref{lem:rho=1} and $\mpr{v(E_1),v(E_2)}<0$ also implies $\bbQ v(E_1)\neq \bbQ v(E_2)$. Then since $\phi \ne 0$, we have $\mu(E_1)<\mu(E_2)$.
Similarly we have $\mu(F_1)<\mu(F_2)$.

{\bf Step 2}.
There exist morphisms $f_2 \colon E_2\to F_2$ and $g_2 \colon F_2\to E_2$ which make the next diagram commutative.
\begin{align}
\label{diag:step3}
\xymatrix{
 E_1 \ar[r]^{\phi} 
&E_2 \ar[r]^{p} \ar@<-0.3ex>[d]_{f_2}
&E^\bullet  \ar[r]  \ar@{=}[d]
&E_1[1]\ 
\\
 F_1 \ar[r]_{\psi}
&F_2 \ar[r]_{q} \ar@<-0.3ex>[u]_{g_2}
&E^\bullet   \ar[r]
&F_1[1].
}
\end{align}

By assumption we have the next exact sequence.
\begin{align*}
\Hom(E_2,F_2) \to \Hom(E_2,E^\bullet) \to \Ext^1(E_2,F_1).
\end{align*}
By Step 1 the last term vanishes. 
Therefore any morphism $F_2\to E^\bullet$ lifts to $E_2\to F_2$. 
Since we are given a non-zero morphism $q \colon F_2\to E^\bullet$, 
there exists a non-zero morphism $f_2 \colon E_2\to F_2$. 
It is easy to verify that the resulting diagram is commutative. 
By the symmetry of the conditions, 
there exists a non-zero morphism $g_2 \colon F_2\to E_2$ 
which makes the resulting diagram commutative. 
Thus we have the diagram \eqref{diag:step3}.

{\bf Step 3}. $f_2$ is an isomorphism. 
We also have an isomorphism $E_1 \simto F_1$ 
and the resulting diagram \eqref{diag:unique} is commutative.

Let $h\seteq g_2\circ f_2-\id_{E_2}\in \End(E_2)$. 
Then the commutativity of the diagram \eqref{diag:step3} yields $p \circ h=0$, 
so that $h$ induces $\widetilde{h}\colon E_2\to E_1$. 
Hence from $\mu(E_1)<\mu(E_2)$ (Step 1) and the semi-stability of $E_i$, 
we have $\Hom(E_2,E_1)=0$. 
Therefore $\widetilde{h}=0$, so that $h=0$. 
Then $g_2 \circ f_2={\rm id}_{E_2}$. 
By the symmetry of the condition, 
we also have $f_2 \circ g_2=\id_{F_2}$.
Thus $f_2$ is an isomorphism.

The second statement is shown by routine procedures.
\end{proof}

\subsection{Numerical criterion of semi-homogeneous presentation}
\label{subsec:uniquness}

Now we state the numerical criterion for the existence of semi-homogeneous presentation. Fix an abelian surface $X$.

\begin{thm}\label{thm:presentation}
Let $\NS(X)=\bbZ H$ and $v$ be a Mukai vector with $\mpr{v^2}>0$. 
\begin{enumerate}
\item
If there is at least two numerical solutions for $v$, 
then a general member of $M_X^H(v)$ has 
both kernel presentation and cokernel presentation. 
Each presentation is unique.
\item
If there is only one numerical solution for $v$, 
then a general member of $M_X^H(v)$ has 
either kernel presentation or cokernel presentation. 
Such a presentation is unique.
\end{enumerate}
\end{thm}

The proof will be presented in \S\S\,\ref{subsect:prf}.

\section{Moduli of simple complexes of numerical invariant $(v_1,v_2,\ell_1,\ell_2)$}
\label{sect:moduli}

In this section we construct fine moduli spaces of simple complexes, which will be used in the proof of Theorem~\ref{thm:presentation}. We fix an abelian surface $X$ with $\NS(X)=\bbZ H$.

\subsection{Complexes associated to numerical solutions}

In this subsection we fix a positive Mukai vector $v\in H^{\ev}(X,\bbZ)_\alg$ with $\ell\seteq\mpr{v^2}/2\in\bbZ_{>0}$. 

Let $(v_1,v_2,\ell_1,\ell_2)$ be a numerical solution of $v$.
For $i=1,2$, let $Y_i\seteq M_X^H(v_i)$ and ${\bfE}_i$ be a universal family
such that $v(\Phi_{X \to Y_i}^{{\bfE}_i^{\vee}}(v_j))=(1,0,0)$, $j \ne i$. 
By Lemma~\ref{lem:picard}, $\NS(Y_i)$ is of rank one. 
We define $H_i$ to be the ample generator of $\NS(Y_i)$.
We also set ${\bfP}\seteq 
{\bfR}p_{Y_1 \times Y_2 *}({\bfE}_1^{\vee} \otimes_{{\calO}_X} {\bfE}_2)$,
where $p_{Y_1 \times Y_2}:Y_1 \times X \times Y_2 \to Y_1 \times Y_2$
is the projection.
Then ${\bfP}$ is a line bundle on $Y_1 \times Y_2$ with
${\bfP}|_{\{ y_1\} \times Y_2} \in \Piczero{Y_2}$
and ${\bfP}|_{Y_1 \times \{ y_2 \}} \in \Piczero{Y_1}$.
It is easy to see that $Y_2 \cong \widehat{Y}_1$ and 
${\bfP}$ is the Poincar\'{e} line bundle on $Y_1 \times Y_2$.

For each numerical solution,
we will construct families of complexes with Mukai vector
$\ell_2 v_2-\ell_1 v_1$ parametrized by
the moduli spaces of rank 1 sheaves $Y_j \times \Hilb{\ell}{Y_i}$.
Here the indices $\{i,j\}=\{1,2\}$ are defined to be $(\ell_i,\ell_j)=(\ell,1)$. 
We start with characterizations of our complexes.

\begin{lem}\label{lem:simple-complex1}
Let $V^{\bullet}$ be a bounded complex on $X$.
The following three conditions are equivalent.
\begin{enumerate}
\item
$V^{\bullet}$ is simple, that is, $\Hom(V^{\bullet},V^{\bullet}) \cong {\bbC}$,  and
$H^i(V^{\bullet})=0$ for $i \ne -1,0$,
$H^{-1}(V^{\bullet}) \in M_X^H(\ell_1 v_1)^{ss}$ and
$H^0(V^{\bullet}) \in M_X^H(\ell_2 v_2)^{ss}$.
\item
$\Phi_{X \to Y_1}^{{\bfE}_1^{\vee}}(V^{\bullet}) \in 
M_{Y_1}^{H_1}(\ell_2,0,-\ell_1)$.
\item
$\Phi_{X \to Y_2}^{{\bfE}_2}\bigl((V^{\bullet})^{\vee}[1]\bigr) \in 
M_{Y_2}^{H_2}(\ell_1,0,-\ell_2)$.
\end{enumerate}
\end{lem}

\begin{proof}  
(a)
We first prove the equivalence of (2) and (3).
From $\Phi_{Y_2 \to Y_1}^{{\bfP}}=
\Phi_{X \to Y_1}^{{\bfE}_1^{\vee}}\Phi_{Y_2 \to X}^{{\bfE}_2}$,
we find that
$\Phi_{X \to Y_2}^{{\bfE}_2^{\vee}[2]}=
\Phi_{Y_1 \to Y_2}^{{\bfP}^{\vee}[2]}
\Phi_{X \to Y_1}^{{\bfE}_1^{\vee}}$.
Then Lemma~\ref{lem:Serre} yields that
\begin{align*}
\Phi_{X \to Y_2}^{{\bfE}_2}\calD_X
=\calD_{Y_2}\,
\Phi_{X \to Y_2}^{{\bfE}_2^{\vee}[2]}
=\calD_{Y_2}\,
\Phi_{Y_1 \to Y_2}^{{\bfP}^{\vee}[2]}\,
\Phi_{X \to Y_1}^{{\bfE}_1^{\vee}}
=\Phi_{Y_1 \to Y_2}^{{\bfP}}\,
\calD_{Y_1}\,
\Phi_{X \to Y_1}^{{\bfE}_1^{\vee}}.
\end{align*}
Hence 
\begin{align*}
\Phi_{Y_1 \to Y_2}^{{\bfP}[1]}\calD_{Y_1}
\bigl(\Phi_{X \to Y_1}^{{\bfE}_1^{\vee}}(V^{\bullet})\bigr)=
\Phi_{X \to Y_2}^{{\bfE}_2}\bigl((V^{\bullet})^{\vee}[1]\bigr).
\end{align*}
Then the equivalence follows from Lemma~\ref{lem:dual} below.

(b)
Next we prove the statement (1) $\Rightarrow$ (2) (or (3)). 
Assume that $V^{\bullet}$ satisfies (1).
Then we have an exact triangle 
\begin{align*}
H^{-1}(V^{\bullet})[1] \to V^{\bullet} \to
H^0(V^{\bullet}) \to H^{-1}(V^{\bullet})[2].
\end{align*} 
Since $V^{\bullet}$ is simple, 
we have $\Hom\left(H^0(V^{\bullet}),H^{-1}(V^{\bullet})[2]\right) \ne 0$.
Since $\Phi_{X \to Y_1}^{{\bfE}_1^{\vee}}\left(H^{-1}(V^{\bullet})[2]\right) 
\in M_{Y_1}^{H_1}(0,0,\ell_1)^{ss}$  
and 
$\Phi_{X \to Y_1}^{{\bfE}_1^{\vee}}
(H^0(V^{\bullet})) \in M_{Y_1}^{H_1}(\ell_2,0,0)^{ss}$,
$H^i\bigl(\Phi_{X \to Y_1}^{{\bfE}_1^{\vee}}(V^{\bullet})\bigr)=0$ for
$i \ne 0,1$ and we have an exact sequence 
\begin{align}
\label{eq:cpx1}
0 \to 
H^0\bigl(\Phi_{X \to Y_{1}}^{{\bfE}_{1}^{\vee}}(V^{\bullet})\bigr) \to 
\Phi_{X \to Y_1}^{{\bfE}_1^{\vee}}\bigl(H^0(V^{\bullet})\bigr) \to
\Phi_{X \to Y_1}^{{\bfE}_1^{\vee}}\bigl(H^{-1}(V^{\bullet})[2]\bigr) \to
H^1\bigl(\Phi_{X \to Y_1}^{{\bfE}_1^{\vee}}(V^{\bullet})\bigr) 
\to 0.
\end{align}

We will prove that 
$H^1\bigl(\Phi_{X \to Y_1}^{{\bfE}_1^{\vee}}(V^{\bullet})\bigr)=0$
by showing 
$\Hom\bigl({\bfE}_{1}|_{\{ y \} \times X},V^{\bullet}[1]\bigr)=0$
for all $y \in Y_1$.
Assume that $\Hom\bigl({\bfE}_1|_{\{ y \} \times X},V^{\bullet}[1]\bigr) \ne 0$
for a point $y \in Y_1$.
Then the Serre duality implies that
$\Hom\bigl(V^{\bullet},{\bfE}_1|_{\{ y \} \times X}[1]\bigr) \ne 0$.
Since $\mpr{ v(H^0(V^{\bullet})),v_1 }=-1$,
Fact~\ref{fact:semihom-ext} implies that
$\Hom\bigl(H^0(V^{\bullet}),{\bfE}_1|_{\{ y \} \times X}[1]\bigr)=0$.
Hence we have a morphism 
$\varphi:V^{\bullet} \to {\bfE}_1|_{\{ y \} \times X}[1]$ 
such that 
$H^{-1}(V^{\bullet}) \to H^{-1}\bigl({\bfE}_1|_{\{ y \} \times X}[1]\bigr)$ 
is surjective in $\Coh(X)$.
By Lemma~\ref{lem:iso} below and the semi-stability,
there is an injective homomorphism
${\bfE}_1|_{\{ y \} \times X} \to H^{-1}(V^{\bullet})$. 
Hence we have a morphism $V^{\bullet} \to 
{\bfE}_1|_{\{ y \} \times X}[1]
\to H^{-1}(V^{\bullet})[1] \to V^{\bullet}$ which induces
a non-zero homomorphism $H^{-1}(V^{\bullet}) \to 
H^{-1}({\bfE}_1|_{\{ y \} \times X}[1]) \to
H^{-1}(V^{\bullet})$ and the zero map
$H^0(V^{\bullet}) \to H^0(V^{\bullet})$.
Since $V^{\bullet}$ is simple, this is impossible.
Therefore $\Hom({\bfE}_1|_{\{ y \} \times X},V^{\bullet}[1])= 0$
for all $y \in Y_{1}$.

By \eqref{eq:cpx1}, $\Phi_{X \to Y_1}^{{\bfE}_1^{\vee}}(V^{\bullet})$
is a $\mu$-semi-stable sheaf with Mukai vector $(\ell_2,0,-\ell_1)$.
We will prove that $\Phi_{X \to Y_1}^{{\bfE}_1^{\vee}}(V^{\bullet})$
is stable.
The claim for $\ell_2=1$ is trivial. Hence we assume that
$(\ell_1,\ell_2)=(1,\ell)$.
Assume that there is a stable subsheaf $F$ of 
$\Phi_{X \to Y_1}^{{\bfE}_1^{\vee}}(V^{\bullet})$
such that  
\begin{align*}
\rk(F)<\ell,\quad
c_1(F)=0,\quad
\chi(F)/\rk F \ge -1/\ell.
\end{align*}
Set $v(F)=(r,0,a)$. Then since
$\mpr{v(F)^2}/2=-r a \geq 0$, we have $a=0$, which yields
$F \in M_{Y_1}^{H_1}(1,0,0)$.

Since $F$ is a subsheaf of a semi-homogeneous sheaf
$\Phi_{X \to Y_1}^{{\bf E}_1^{\vee}}(H^0(V^{\bullet}))$
with $v(F)=(\ell,0,0)$,
we find from Lemma~\ref{lem:iso} below that
$\Hom\bigl(\Phi_{X \to Y_1}^{{\bf E}_1^{\vee}}(H^0(V^{\bullet})),F\bigr)\ne0$.
Then we see that 
$\Phi_{X \to Y_1}^{{\bfE}_1^{\vee}}(V^{\bullet})$
is not simple.
Therefore $\Phi_{X \to Y_1}^{{\bfE}_1^{\vee}}(V^{\bullet})$ is stable. 
Thus (2) holds.

(c)
Finally we prove the statement (2) $\Rightarrow$ (1).
Assume that (2) (and hence (3)) holds.  

We first assume that $\ell_2=1$.
Then $\Phi_{X \to Y_1}^{{\bfE}_1^{\vee}}(V^{\bullet}) \cong
L \otimes I_Z$ has the kernel presentation
\begin{align*}
0 \to L \otimes I_Z \to L \to L \otimes {\calO}_Z \to 0,
\end{align*}
where $L \in \Piczero{Y_1}$ and $I_Z \in \Hilb{\ell}{Y_1}$.
Applying $\Phi_{Y_1 \to X}^{{\bfE}_1}$ to this exact sequence,
we get 
\begin{align*}
H^i(V^{\bullet})=0\; (i \ne -1,0), \quad
H^{-1}(V^{\bullet})=\Phi_{Y_1 \to X}^{{\bfE}_1[1]}(L \otimes {\calO}_Z),\quad
H^0(V^{\bullet})=\Phi_{Y_1 \to X}^{{\bfE}_1[2]}(L).
\end{align*}
Thus (1) holds.

If $\ell_1=1$, then 
$\Phi_{X \to Y_2}^{{\bfE}_2}\bigl((V^{\bullet})^{\vee}[1]\bigr) \cong
L \otimes I_Z$ has the kernel presentation
\begin{align*}
0 \to L \otimes I_Z \to L \to L \otimes {\calO}_Z \to 0,
\end{align*}
where $L \in \Piczero{Y_2}$ and $I_Z \in \Hilb{\ell}{Y_2}$.
Applying $\calD_X\, \Phi_{Y_2 \to X}^{{\bfE}_2^{\vee}}
= \Phi_{Y_2 \to X}^{{\bfE}_2[2]}\, \calD_{Y_2}$ 
(see Lemma~\ref{lem:Serre})
to this exact sequence,
we get 
\begin{align*}
H^i(V^{\bullet})=0\; (i \ne -1,0), \quad
H^{-1}(V^{\bullet})=
\Phi_{Y_2 \to X}^{{\bfE}_2[2]}\bigl((L \otimes {\calO}_Z)^{\vee}\bigr),\quad
H^0(V^{\bullet})=
\Phi_{Y_2 \to X}^{{\bfE}_2[2]}(L^{\vee}).
\end{align*}
Thus (1) holds again in this case.
\end{proof}

\begin{lem}
\label{lem:iso}
Let $E$ be a simple semi-homogeneous sheaf and
$F$ a semi-homogeneous sheaf with $v(F)=\ell v(E)$, $\ell \in {\bbZ}_{>0}$.
Then $\Hom(E,F) \ne 0$ if and only if $\Hom(F,E) \ne 0$. 
\end{lem}

\begin{proof}
See \cite[Proposition 6.10]{Mukai:1978}.
\end{proof}

\begin{lem}
\label{lem:dual}
The functors $\Phi_{Y_1 \to Y_2}^{{\bfP}[1]}\calD_{Y_1}$ and 
$\Phi_{Y_1 \to Y_2}^{{\bfP}^{\vee}[1]}\calD_{Y_1}$ induce 
an isomorphism
$M_{Y_1}^{H_1}(r,0,-a) \cong M_{Y_2}^{H_2}(a,0,-r)$.
\end{lem}

\begin{proof}
See \cite[Theorem 8.4]{Yoshioka:2001:ann} or 
\cite[Corollary 4.5]{Mukai:1987:ASPM} for $a=1$.
\end{proof}

\begin{rmk}
$\Phi_{Y_2 \to Y_1}^{{\bfP}[1]}\,\calD_{Y_2}\,
\Phi_{Y_1 \to Y_2}^{{\bfP}[1]}\,\calD_{Y_1} \cong 1_{Y_1}$.
Hence $\Phi_{Y_2 \to Y_1}^{{\bfP}[1]}\calD_{Y_2}$ gives the inverse
of $\Phi_{Y_1 \to Y_2}^{{\bfP}[1]}\calD_{Y_1}$.
\end{rmk}

\begin{lem}\label{lem:simple-complex2}
Let $V^{\bullet}$ be a bounded complex on $X$.
The following three conditions are equivalent.
\begin{enumerate}
\item
$V^{\bullet}$ is simple and quasi-isomorphic to a complex 
$\bigl[U^{-1} \to U^0\bigr]$ with
$U^{-1} \in M_X^H(\ell_1 v_1)^{ss}$ and 
$U^0 \in M_X^H(\ell_2 v_2)^{ss}$. 
\item
$\Phi_{X \to Y_2}^{{\bfE}_2^{\vee}}(V^{\bullet}[1])
\in M_X^H(\ell_1,0,-\ell_2)$.
\item
$\calD_{Y_1}\Phi_{X \to Y_1}^{{\bf E}_1^{\vee}}(V^{\bullet})=
\Phi_{X \to Y_1}^{{\bfE}_1}\bigl((V^{\bullet})^{\vee}[2]\bigr)
\in M_X^H(\ell_1,0,-\ell_2)$. 
\end{enumerate}
\end{lem}

\begin{proof}
$\Phi_{Y_2 \to Y_1}^{{\bfP}}=
\Phi_{X \to Y_1}^{{\bfE}_1^{\vee}}\,\Phi_{Y_2 \to X}^{{\bfE}_2}$
implies that
$\Phi_{Y_1 \to Y_2}^{{\bfP}^{\vee}}=
\Phi_{X \to Y_2}^{{\bfE}_2^{\vee}}\,\Phi_{Y_1 \to X}^{{\bfE}_1[2]}$.
Using Lemma~\ref{lem:Serre} we have 
\begin{align*}
(\Phi_{Y_1 \to Y_2}^{{\bfP}^{\vee}}\,\calD_{Y_1})\,
\Phi_{X \to Y_1}^{{\bfE}_1[2]}\,\calD_X
=
\Phi_{Y_1 \to Y_2}^{{\bfP}^{\vee}}\,
\Phi_{X \to Y_1}^{{\bfE}_1^{\vee}}
=\Phi_{X \to Y_2}^{{\bfE}_2^{\vee}}.
\end{align*}
Therefore Lemma~\ref{lem:dual} implies that (2) and (3) are equivalent.

The proof of the equivalence between (1) and (2) (or (3)) 
is similar to that in Lemma~\ref{lem:simple-complex1}. 
Hence we only prove the direction (1) $\Rightarrow$ (2).
Assume that (1) holds. 

We will show that 
$\Hom\bigl({\bfE}_2|_{\{ y \} \times X},V^{\bullet}[2]\bigr)=0$
for all $y \in Y_2$.
If 
$\Hom\bigl(V^{\bullet},{\bfE}_2|_{\{ y \} \times X}\bigr)
=\Hom\bigl({\bfE}_2|_{\{ y \} \times X},V^{\bullet}[2]\bigr)^{\vee} \ne 0$,
then we have a non-zero homomorphism 
$V^0 \to {\bfE}_2|_{\{ y \} \times X}$.
By our choice of $v_2$, it is surjective.
Then by Lemma~\ref{lem:iso},
we also have an injective homomorphism
${\bfE}_2|_{\{ y \} \times X} \to V^0$, which gives a morphism
${\bfE}_2|_{\{ y \} \times X} \to V^{\bullet}$.
Hence we have a non-zero morphism 
$V^{\bullet} \to {\bfE}_2|_{\{y\}\times X} \to V^{\bullet}$.
Obviously it is not isomorphic.
By the simpleness of $V^{\bullet}$, it is impossible.
Therefore 
$\Hom\bigl({\bfE}_2|_{\{y\} \times X},V^{\bullet}[2]\bigr)=0$ for
all $y \in Y_2$.

Then $\Phi_{X \to Y_2}^{{\bfE}_2^\vee}\bigl((V^{\bullet})^{\vee}[1]\bigr)$ is
a simple sheaf.
Then the same arguments as in 
Lemma~\ref{lem:simple-complex1} show the stability.
\end{proof}

\subsection{Existence and some properties of the moduli spaces}

As in the last subsection, we fix a positive Mukai vector $v\in H^{\ev}(X,\bbZ)_\alg$ with $\ell\seteq\mpr{v^2}/2\in\bbZ_{>0}$. For a numerical solution $(v_1,v_2,\ell_1,\ell_2)$ of $v$, we continue to use the symbols $Y_i$, ${\bfE}_i$, $H_i$, $\bfP$.

Let ${\calZ}_1 \subset \Hilb{\ell}{Y_1} \times Y_1$
and ${\calZ}_2 \subset \Hilb{\ell}{Y_2} \times Y_2$ 
be the universal family of $\ell$-points on $Y_1$ and $Y_2$ respectively.
We set ${\frakH}_1\seteq Y_2 \times \Hilb{\ell}{Y_1}$ and
${\frakH}_2\seteq Y_1 \times \Hilb{\ell}{Y_2}$.
Then we have an exact sequence on ${\frakH}_i \times Y_i$:
\begin{align}
\label{eq:kernel-presentation}
0 \to {\bfP} \otimes I_{{\calZ}_i} \to {\bfP}
  \to {\bfP} \otimes {\calO}_{{\calZ}_i} \to 0,
\end{align}
where $i=1,2$ and we simply denote the pull-backs of
${\bfP}$, $I_{{\calZ}_i}$ and ${\calO}_{{\calZ}_i}$
by the same symbol.
For $i=1,2$, we set $c(1)\seteq 2$ and $c(2)\seteq 1$.
Applying $\Phi_{Y_i \to Y_j}^{{\bfP}[1]}$ to the dual
of \eqref{eq:kernel-presentation}, 
we have an exact sequence on ${\frakH}_i \times Y_{c(i)}$:
\begin{align}
\label{seq:09}
0 \to 
\Phi_{Y_i \to Y_{c(i)}}^{{\bfP}[1]}
 \bigl(({\bfP} \otimes I_{{\calZ}_i})^{\vee}\bigr) 
\to 
\Phi_{Y_i \to Y_{c(i)}}^{{\bfP}[2]}
 \bigl(({\bfP} \otimes {\calO}_{{\calZ}_i})^{\vee}\bigr) 
\to
\Phi_{Y_i \to Y_{c(i)}}^{{\bfP}[2]}({\bfP}^{\vee}) 
\to 0.
\end{align}
Then for $i=1,2$,
${\bfP} \otimes I_{{\calZ}_i}$ and 
$\Phi_{Y_{c(i)} \to Y_i}^{{\bfP}[1]}
(({\bfP} \otimes I_{{\calZ}_{c(i)}})^{\vee})$ 
are universal families of stable sheaves on $Y_i$ with
Mukai vectors $(1,0,-\ell)$ and $(\ell,0,-1)$ respectively.
$\Phi_{Y_{c(i)} \to Y_i}^{{\bfP}^{\vee}[1]}
\bigl(({\bfP} \otimes I_{{\calZ}_{c(i)}})^{\vee}\bigr)$ is also
a universal family of stable sheaves on $Y_i$ with Mukai vector $(\ell,0,-1)$. 

\begin{rmk}
$\Phi_{Y_{c(i)} \to Y_i}^{{\bfP}[2]}
({\bfP}^{\vee}) \cong {\calO}_{{\Delta}}$,
where $\Delta \subset Y_i \times Y_i$ is the diagonal.
\end{rmk}

Applying $\Phi_{Y_i \to X}^{{\bfE}_i}$ 
to the exact sequences \eqref{eq:kernel-presentation} and \eqref{seq:09}, 
we have the following families of simple complexes.

\begin{prop}\label{prop:univ-cpx}
For a numerical solution $(v_1,v_2,\ell_1,\ell_2)$ of $v$, we set the indices $i,j\in\{1,2\}$ by $(\ell_i,\ell_j)=(\ell,1)$ and set ${\frakS}\seteq{\frakH}_i$.

\begin{enumerate}
\item 
We have families of simple complexes 
${\calU}^{\bullet}=\bigl[{\calU}^{-1} \to {\calU}^0\bigr]$ of 
$v({\calU}_s^{\bullet})=\ell_2 v_2-\ell_1 v_1$, 
$s \in {\frakS}$ with the following properties. 
\begin{enumerate}
\item
If $\ell_1=\ell$, we set
${\calU}^{\bullet}\seteq
\Phi_{Y_1 \to X}^{{\bfE}_1[2]}({\bfP} \otimes I_{{\calZ}_1})$.
Then
$H^{-1}({\calU}^{\bullet})=
\Phi_{Y_1 \to X}^{{\bfE}_1}({\bfP} \otimes {\calO}_{{\calZ}_1})$
and $H^0({\calU}^{\bullet})=
\Phi_{Y_1 \to X}^{{\bfE}_1[2]}({\bfP})$ 
are semi-homogeneous sheaves with Mukai vectors $\ell v_1$ and $ v_2$
respectively.
\item
If $\ell_2=\ell$, we set
${\calU}^{\bullet}\seteq
\Phi_{Y_1 \to X}^{{\bfE}_1[2]}\,\Phi_{Y_2 \to Y_1}^{{\bfP}[1]}
 \bigl(({\bfP} \otimes I_{{\calZ}_2})^{\vee}\bigr)
=
\Phi_{Y_2 \to X}^{{\bfE}_2[1]}
 \bigl(({\bfP} \otimes I_{{\calZ}_2})^{\vee}\bigr)$.
Then 
$H^{-1}({\calU}^{\bullet})=\Phi_{Y_2 \to X}^{{\bfE}_2}({\bfP}^{\vee})$
and 
$H^0({\calU}^{\bullet})
=
\Phi_{Y_2 \to X}^{{\bfE}_2[2]}
 \bigl(({\bfP} \otimes {\calO}_{{\calZ}_2})^{\vee}\bigr)$ 
are semi-homogeneous sheaves with Mukai vectors $v_1$ and $\ell v_2$
respectively.
\end{enumerate}
\item
We have families of simple complexes 
${\calV}^{\bullet}:{\calV}^{-1} \to {\calV}^0$ of 
$v({\calV}_s^{\bullet})=\ell_2 v_2-\ell_1 v_1$,
$s \in {\frakS}$ with the following properties. 
\begin{enumerate}
\item
If $\ell_2=\ell$, we set
${\calV}^{\bullet}\seteq
\Phi_{Y_2 \to X}^{{\bfE}_2[1]}({\bfP} \otimes I_{{\calZ}_2})$.
Then 
${\calV}^{-1}=\Phi_{Y_2 \to X}^{{\bfE}_2}({\bfP})$ 
and 
${\calV}^0=\Phi_{Y_2 \to X}^{{\bfE}_2}
({\bfP} \otimes {\calO}_{{\calZ}_2})$ 
are families of the semi-homogeneous sheaves
with Mukai vectors $v_1$ and $\ell v_2$ respectively.
\item
If $\ell_1=\ell$, we set
${\calV}^{\bullet}=
\Phi_{Y_2 \to X}^{{\bfE}_2[1]}
\Phi_{Y_1 \to Y_2}^{{\bfP}^{\vee}[1]}
\bigl(({\bfP} \otimes I_{{\calZ}_1})^{\vee}\bigr)
=
\Phi_{Y_1 \to X}^{{\bfE}_1[2]}
 \bigl(({\bfP} \otimes I_{{\calZ}_1})^{\vee}\bigr)$.
Then
${\calV}^{-1}=
\Phi_{Y_1 \to X}^{{\bfE}_1}
\bigl(({\bfP} \otimes {\calO}_{{\calZ}_1}[2])^{\vee}\bigr)$ 
and 
${\calV}^0=
\Phi_{Y_1 \to X}^{{\bfE}_1}({\bfP}^{\vee})$ 
are families of the semi-homogeneous sheaves
with Mukai vectors $\ell v_1$ and $v_2$ respectively.
\end{enumerate}
\end{enumerate}
\end{prop}

\begin{prop}\label{prop:open}
For a scheme $S$, let
$V^{\bullet}$ be a family of simple complexes on $X$
such that each $V^i$ is flat over $S$.
\begin{enumerate}
\item
We set
\begin{align*}
S' \seteq 
\left\{
 s \in S\, 
 \left|\,
 \begin{aligned}
  &H^i(V^{\bullet}_s)=0\, (i \ne -1,0),\\
  &H^{-1}(V^{\bullet}_s) \in M_X^H(\ell_1 v_1)^{ss},\;
   H^0(V^{\bullet}_s) \in M_X^H(\ell_2 v_2)^{ss}
 \end{aligned}
 \right. 
\right\}.
\end{align*}
Then $S'$ is an open subset of $S$ and 
there is a unique morphism
$f:S' \to {\frakS}$ such that
$(f \times 1_X)^*({\calU}) \cong V^{\bullet}|_{S' \times X} 
\otimes p_{S'}^*(N)$,
where $N \in \Pic{S'}$ and $p_{S'}:S' \times X \to S'$ is the
projection.
\item
We set
\begin{align*}
S' \seteq 
 \left\{ s \in S \mid 
  V^{\bullet}_s \cong \bigl[W_s^{-1} \to W_s^0\bigr],\; 
  W_s^{-1} \in M_X^H(\ell_1 v_1)^{ss},\;
  W_s^0 \in M_X^H(\ell_2 v_2)^{ss} 
 \right\}.
\end{align*}
Then $S'$ is an open subset of $S$ and
we have a unique morphism $f:S' \to {\frakS}$ such that
$(f \times 1_X)^*({\calV}) \cong 
V^{\bullet}|_{S' \times X} \otimes p_{S'}^*(N)$, 
where $N \in \Pic{S'}$ and $p_{S'}:S' \times X \to S'$ is the
projection.
\end{enumerate}
\end{prop}

\begin{proof}
(1)
By Lemma~\ref{lem:simple-complex1},
$S'$ is the set of $s \in S$ such that
$\Phi_{X \to Y_1}^{{\bfE}_1^{\vee}}(V^{\bullet})_s$
is a stable sheaf with 
$v(\Phi_{X \to Y_1}^{{\bfE}_1^{\vee}}(V^{\bullet})_s)=
(\ell_2,0,-\ell_1)$.
By the openness of stability \cite[Proposition 2.3.1]{HuybrechtsLehn:book}, 
$S'$ is an open subset of $S$.
Then we have a unique morphism $f \colon S' \to {\frakS}$ 
and a line bundle $N \in \Pic{S'}$ such that
$\Phi_{X \to Y_1}^{{\bfE}_1^{\vee}}(V^{\bullet})|_{S' \times Y_1}
\otimes p_{S'}^*(N)$ is the pull-back 
of the universal family on $ {\frakH} \times Y_1$.  
Hence the claim holds by Proposition~\ref{prop:univ-cpx}.
The proof of (2) is similar. 
\end{proof}

By Proposition \ref{prop:open}, 
the universal families ${\calU}$ and ${\calV}$ in Proposition 
\ref{prop:univ-cpx} give moduli spaces of simple complexes.

\begin{thm}
\label{thm:moduli}
Let $v$ be a positive Mukai vector with $\mpr{v^2}>0$ and 
$(v_1,v_2,\ell_1,\ell_2)$ be a numerical solution of $v$.
\begin{enumerate}
\item
We have the fine moduli space ${\frakM}^{-}(v_1,v_2,\ell_1,\ell_2)$
of simple complexes $V^{\bullet}$ such that
$H^i(V^{\bullet})=0, i \ne -1,0$,
$H^{-1}(V^{\bullet}) \in M_X^H(\ell_1 v_1)^{ss}$ and 
$H^0(V^{\bullet}_s) \in M_X^H(\ell_2 v_2)^{ss}$.
\item
We have the fine moduli space ${\frakM}^+(v_1,v_2,\ell_1,\ell_2)$
of simple complexes $V^{\bullet}$ such that
$V^{\bullet} \cong [W^{-1} \to W^0]$,  
$W^{-1} \in M_X^H(\ell_1 v_1)^{ss}$,
$W^0 \in M_X^H(\ell_2 v_2)^{ss}$.
\end{enumerate}
\end{thm}

The following assertions show that 
${\frakM}^{\pm}(v_1,v_2,\ell_1,\ell_2)$ behaves very well
under the Fourier-Mukai transforms.

\begin{prop}
\label{prop:FMT_M_pm}
For a FMT $\Phi_{X \to Y}^{{\bf E}^{\vee}}$,
${\bf E} \in \Coh(X \times Y)$,
we set $v_i' \seteq \Phi_{X \to Y}^{{\bf E}^{\vee}}(v_i)$  ($i=1,2$).
Then the following assertions hold.
\begin{enumerate}
\item
If $\mu({\bfE}|_{\{ x\} \times Y})<\mu(v_1)$
or $\mu(v_2) \leq \mu({\bfE}|_{\{x\} \times Y})$, then
we have an isomorphism
\begin{align*}
{\frakM}^{\pm}(v_1,v_2,\ell_1,\ell_2) \simto
{\frakM}^{\pm}(v_1',v_2',\ell_1,\ell_2).
\end{align*}
\item
If $\mu(v_1) \le \mu({\bfE}_{|\{ x\} \times Y})
<\mu(v_2)$, then we have an isomorphism
\begin{align*}
{\frakM}^{\pm}(v_1,v_2,  \ell_1,\ell_2) \simto
{\frakM}^{\mp}(v_2',v_1',\ell_2,\ell_1).
\end{align*}
\end{enumerate}
\end{prop}

\begin{proof}
Let $E_i$ ($i=1,2$) be semi-homogeneous sheaves on $X$
with Mukai vectors $\ell_i v_i$.
Then we have
\begin{align*}
\begin{cases}
\Phi_{X \to Y}^{{\bfE}^{\vee}}(E_1),\ 
\Phi_{X \to Y}^{{\bfE}^{\vee}}(E_2) \in \Coh(Y)
& 
\mu({\bf E}|_{\{x\} \times Y})<\mu(v_1),
\\   
\Phi_{X \to Y}^{{\bfE}^{\vee}[2]}(E_1),\ 
\Phi_{X \to Y}^{{\bfE}^{\vee}}(E_2) \in \Coh(Y)
& 
\mu(v_1) \leq \mu({\bfE}|_{\{ x\} \times Y})
<\mu(v_2),
\\
\Phi_{X \to Y}^{{\bfE}^{\vee}[2]}(E_1),\  
\Phi_{X \to Y}^{{\bfE}^{\vee}[2]}(E_2) \in \Coh(Y)
& 
\mu(v_2) \leq \mu({\bfE}_{|\{ x\} \times Y}).
\end{cases}
\end{align*}
{}From these results, we can easily prove the claim. 
\end{proof}

For $V^{\bullet} \in {\frakM}^{+}(v_1,v_2,\ell_1,\ell_2)$,
we set
\begin{align*}
\Psi(V^{\bullet})\seteq
\begin{cases}
\Phi_{Y_1 \to X}^{{\bfE}_1[2]}\,\calD_{Y_1}\,
\Phi_{X \to Y_1}^{{\bfE}_1^{\vee}}(V^{\bullet})
& \ell_1=\ell, 
\\
\Phi_{Y_2 \to X}^{{\bfE}_2[1]}\,\calD_{Y_2}\,
\Phi_{X \to Y_2}^{{\bfE}_2^{\vee}[1]}(V^{\bullet})
& \ell_2=\ell. 
\end{cases}
\end{align*}
Then Proposition~\ref{prop:univ-cpx} immediately implies the following.

\begin{prop}
\label{prop:dual}
$\Psi$ induces an isomorphism
\begin{align*}
{\frakM}^{+}(v_1,v_2,\ell_1,\ell_2) \simto
{\frakM}^{-}(v_1,v_2,\ell_1,\ell_2). 
\end{align*}
\end{prop}

\begin{prop}\label{prop:cpx-stable}
For a numerical solution $(v_1,v_2,\ell_1,\ell_2)$ of $v$,
we set
\begin{align*}
k\seteq
\begin{cases}
0 & \ell_2 v_2-\ell_1 v_1>0,\\
1 & \ell_1 v_1-\ell_2 v_2>0,
\end{cases}
\end{align*}
and set
\begin{align*}
{\frakM}^{+}(v_1,v_2,\ell_1,\ell_2)_0
\seteq
&
\{ V^{\bullet} \in {\frakM}^{+}(v_1,v_2,\ell_1,\ell_2) 
   \mid
   V^{\bullet}[k] \in M_X^H(v) \},\\
M_X^H(v)_0
\seteq
&
\{E \in M_X^H(v) \mid E[-k] \in {\frakM}^+(v_1,v_2,\ell_1,\ell_2) \}.
\end{align*}
Then we have the following assertions.
\begin{enumerate}
\item
${\frakM}^{+}(v_1,v_2,\ell_1,\ell_2)_0$ and $M_X^H(v)_0$
are open subschemes of 
${\frakM}^{+}(v_1,v_2,\ell_1,\ell_2)$ and $M_X^H(v)$ respectively, and
we have an isomorphism
\begin{align*}
\begin{matrix}
\phi\colon 
 & {\frakM}^{+}(v_1,v_2,\ell_1,\ell_2)_0 
 & \to & M_X^H(v)_0
\\
 & V^{\bullet} 
 & \mapsto & V^{\bullet}[k] .
\end{matrix}
\end{align*}
\item
If there is a complex $V^{\bullet} \in {\frakM}^{+}(v_1,v_2,\ell_1,\ell_2)$
with $V^{\bullet}[k] \in \Coh(X)$ for some $k\in\bbZ$, then
${\frakM}^{+}(v_1,v_2,\ell_1,\ell_2)_0 \ne \emptyset$.
\end{enumerate}
\end{prop}

\begin{proof}
(1) 
The openness of ${\frakM}^{+}(v_1,v_2,\ell_1,\ell_2)_0$ is well-known 
(for example, see \cite[Theorem 1.6 (1)]{Mukai:1987:ASPM})
and the existence of $\phi$ is the result of the universal property
of $M_X^{H}(v)$.
The openness of $M_X^{H}(v)_0$ and the existence of $\phi^{-1}$
follows from Proposition~\ref{prop:open}~(2).

(2) follows from Fact~\ref{fact:yoshioka-lem2-2}.
\end{proof}

\subsection{Proof of Theorem~\ref{thm:presentation}}
\label{subsect:prf}

Before starting the proof, we prepare some lemmas.

\begin{lem}
\label{lem:sheaf-criterion}
Let $E^{\bullet}$ be a bounded complex on $X$.  
Assume that $E^{\bullet}$ has two expressions
\begin{align*}
E^{\bullet} \cong \bigl[E_1 \to E_2\bigr] \cong \bigl[F_1 \to F_2\bigr][-1],
\end{align*}
where $E_1, E_2, F_1, F_2 \in \Coh(X)$.
Then $E^{\bullet} \in \Coh(X)$ and 
$E\seteq H^0(E^{\bullet})$ fits in the following two exact sequences in $\Coh(X)$.
\begin{align*}
0 \to E_1 \to E_2 \to E \to 0,\quad
0 \to E \to F_1 \to F_2 \to 0. 
\end{align*} 
\end{lem}

\begin{proof}
We have two exact triangles:
\begin{align*}
E_1 \to E_2 \to E^{\bullet} \to E_1[1],\quad
E^{\bullet} \to F_1 \to F_2 \to E^{\bullet}[1].
\end{align*}
By taking the cohomology groups of these triangles,
we get the claims.
\end{proof}

\begin{lem}
\label{lem:FM-sheaf-criterion}
Assume that a complex $E^{\bullet}$ on $X$ has two expressions
\begin{align*}
E^{\bullet} \cong \bigl[E_1 \overset{\phi}{\to} E_2\bigr] \cong 
\bigl[F_2 \overset{\psi}{\to} F_1[2]\bigr][-1],
\end{align*}
where $E_1, E_2, F_1, F_2$ are semi-homogeneous sheaves
and $\phi, \psi$ are non-trivial. Namely, we have an exact 
sequence in $\Coh(X)$:
\begin{align*}
0 \to F_1 \to E_1 \to E_2 \to F_2 \to 0.
\end{align*} 
Let $\Phi_{X \to Y}^{{\bfE}^{\vee}}$ be a FMT. Then the following holds.
\begin{enumerate}
\item
\begin{enumerate}
\item
If $\mu(F_1) \leq \mu({\bfE}|_{X \times \{ y \} }) <\mu(E_1)$, then
$\Phi_{X \to Y}^{{\bfE}^{\vee}}(E^{\bullet})\in\Coh(Y)$.

\item
If $\mu(E_2) \leq \mu({\bfE}|_{X \times \{ y \} }) <\mu(F_2)$, then
$\Phi_{X \to Y}^{{\bfE}^{\vee}[1]}(E^{\bullet})\in\Coh(Y)$.

\item
For the other cases, 
$\Phi_{X \to Y}^{{\bfE}^{\vee}[k]}(E^{\bullet})$ 
is not a coherent sheaf on $Y$ for any $k$.
\end{enumerate}
\item
\begin{enumerate}
\item 
If $\mu(F_1) < \mu({\bfE}|_{X \times \{ y \} }) \leq \mu(E_1)$, then
$\calD_Y \Phi_{X \to Y}^{{\bfE}^{\vee}}(E^{\bullet})\in\Coh(Y)$.

\item
If $\mu(E_2) < \mu({\bfE}|_{X \times \{ y \} })  \leq \mu(F_2)$, then
$\calD_Y \Phi_{X \to Y}^{{\bfE}^{\vee}[1]}(E^{\bullet})\in\Coh(Y)$.

\item
For the other cases, 
$\calD_Y \Phi_{X \to Y}^{{\bfE}^{\vee}[k]}(E^{\bullet})$ 
is not a coherent sheaf on $Y$ for any $k$.
\end{enumerate}
\end{enumerate}
\end{lem}

\begin{proof}
(1)
We have an isomorphisms
\begin{align}
\label{eq:two_reps}
\begin{split}
\Phi_{X \to Y}^{{\bfE}^{\vee}}(E^{\bullet})
 \cong  
\bigl[\Phi_{X \to Y}^{{\bfE}^{\vee}}(E_1) 
\to 
\Phi_{X \to Y}^{{\bfE}^{\vee}}(E_2)\bigr]
 \cong 
\bigl[\Phi_{X \to Y}^{{\bfE}^{\vee}}(F_2)
 \to 
\Phi_{X \to Y}^{{\bfE}^{\vee}}(F_1)[2]\bigr][-1].
\end{split}
\end{align}
If $\mu(F_1) \leq \mu({\bfE}|_{X \times \{ y \} }) <\mu(E_1)$, then
\begin{align*}
\Phi_{X \to Y}^{{\bfE}^{\vee}}(F_1[2]),\,
\Phi_{X \to Y}^{{\bfE}^{\vee}}(E_1),\,
\Phi_{X \to Y}^{{\bfE}^{\vee}}(E_2),\,
\Phi_{X \to Y}^{{\bfE}^{\vee}}(F_2) \in \Coh(Y).
\end{align*}
Hence the claim (a) holds by Lemma \ref{lem:sheaf-criterion}.

If $\mu(E_2) \le \mu({\bfE}|_{X \times \{ y \} }) <\mu(F_2)$, then
\begin{align*}
\Phi_{X \to Y}^{{\bfE}^{\vee}}(F_1[2]),\,
\Phi_{X \to Y}^{{\bfE}^{\vee}}(E_1[2]),\,
\Phi_{X \to Y}^{{\bfE}^{\vee}}(E_2[2]),\,
\Phi_{X \to Y}^{{\bfE}^{\vee}}(F_2) \in \Coh(Y).
\end{align*}
Hence the claim (b) holds by Lemma \ref{lem:sheaf-criterion}.

For the proof of (c), we check the following three cases.

If $\mu({\bfE}|_{X \times \{y\} }) < \mu(F_1)$, then
\begin{align*}
\Phi_{X \to Y}^{{\bfE}^{\vee}}(F_1),\,
\Phi_{X \to Y}^{{\bfE}^{\vee}}(E_1),\,
\Phi_{X \to Y}^{{\bfE}^{\vee}}(E_2),\,
\Phi_{X \to Y}^{{\bfE}^{\vee}}(F_2) \in \Coh(Y).
\end{align*}
Thus by \eqref{eq:two_reps} 
$\Phi_{X \to Y}^{{\bfE}^{\vee}[k]}(E^{\bullet})$ cannot be a sheaf for any $k\in\bbZ$.

If $\mu(E_1)\le \mu({\bfE}|_{X \times \{y\} }) <\mu(E_2)$, then
\begin{align*}
\Phi_{X \to Y}^{{\bfE}^{\vee}}(F_1[2]),\,
\Phi_{X \to Y}^{{\bfE}^{\vee}}(E_1[2]),\,
\Phi_{X \to Y}^{{\bfE}^{\vee}}(E_2),\,
\Phi_{X \to Y}^{{\bfE}^{\vee}}(F_2) \in \Coh(Y).
\end{align*}
By \eqref{eq:two_reps} 
$\Phi_{X \to Y}^{{\bfE}^{\vee}[k]}(E^{\bullet})$ cannot be a sheaf for any $k\in\bbZ$.

If $\mu(F_2)\le \mu({\bfE}|_{X \times \{y\} })$, then
\begin{align*}
\Phi_{X \to Y}^{{\bfE}^{\vee}}(F_1[2]),\,
\Phi_{X \to Y}^{{\bfE}^{\vee}}(E_1[2]),\,
\Phi_{X \to Y}^{{\bfE}^{\vee}}(E_2[2]),\,
\Phi_{X \to Y}^{{\bfE}^{\vee}}(F_2[2]) \in \Coh(Y).
\end{align*}
By \eqref{eq:two_reps} 
$\Phi_{X \to Y}^{{\bfE}^{\vee}[k]}(E^{\bullet})$ cannot be a sheaf for any $k\in\bbZ$.

(2)
We have an isomorphisms
\begin{align*}
\Phi_{X \to Y}^{{\bfE}^{\vee}}(E^{\bullet})^{\vee}
\cong 
\bigl[\Phi_{X \to Y}^{{\bfE}^{\vee}}(E_2)^{\vee} \to
      \Phi_{X \to Y}^{{\bfE}^{\vee}}(E_1)^{\vee}\bigr][-1]
\cong
\bigl[\Phi_{X \to Y}^{{\bfE}^{\vee}}(F_1[2])^{\vee} \to 
      \Phi_{X \to Y}^{{\bfE}^{\vee}}(F_2)^{\vee}\bigr].
\end{align*}
If $\mu(F_1) < \mu({\bfE}|_{X \times \{ y \} })  \leq \mu(E_1)$, then
\begin{align*}
\Phi_{X \to Y}^{{\bfE}^{\vee}}(E_1)^{\vee},\,
\Phi_{X \to Y}^{{\bfE}^{\vee}}(E_2)^{\vee},\,
\Phi_{X \to Y}^{{\bfE}^{\vee}}(F_1)^{\vee},\,
\Phi_{X \to Y}^{{\bfE}^{\vee}}(F_1[2])^{\vee} \in \Coh(Y).
\end{align*}
Hence the claim (a) holds by Lemma \ref{lem:sheaf-criterion}.
The other claims are verified similarly as (1).
\end{proof}

Now we can show Theorem~\ref{thm:presentation}.

\begin{proof}[{Proof of Theorem~\ref{thm:presentation}}]
The uniqueness statements follow from 
Proposition~\ref{prop:semihom-presen-unique}. 
In the following we show the existence statements.

(1)
Assume that the numerical equation of $v$ has two solutions 
$(v_1,v_2,\ell_1,\ell_2)$ and $(v'_1,v'_2,\ell'_1,\ell'_2)$. 
By Proposition \ref{prop:univ-cpx},
we have two families of simple complexes
${\calV}^{\bullet}=\bigl[{\calV}^{-1} \to {\calV}^0\bigr]$
and ${{\calV}'}^{\bullet}=\bigl[{{\calV}'}^{-1} \to {{\calV}'}^0\bigr]$
associated to the two numerical solutions
$(v_1,v_2,\ell_1,\ell_2)$ and
$(v'_1,v'_2,\ell'_1,\ell'_2)$.
We set $S \seteq {\frakM}^+(v_1,v_2,\ell_1,\ell_2)$ and 
$S' \seteq {\frakM}^+(v_1',v_2',\ell_1',\ell_2')$.

(a)
We first assume that
${\calV}_s^{\bullet}$ and ${{\calV}'}_{s'}^{\bullet}$
are sheaves up to shift for general $s\in S$ and $s'\in S'$.
By Proposition \ref{prop:cpx-stable},
there are non-empty open subschemes $S_0\subset S$ and $S_0'\subset S'$, 
and we have open immersions 
$f\colon S_0 \to M_X^H(v)$ and $f' \colon S_0' \to M_X^H(v)$.
Since $M_X^H(v)$ is irreducible (Fact \ref{fact:moduli}), $f(S_0) \cap f'(S_0') \ne \emptyset$, which implies that a general member $E \in M_X^H(v)$ has two semi-homogeneous presentations.
By the uniqueness of the presentations (Proposition~\ref{prop:semihom-presen-unique}), one of these two presentations must be a kernel presentation and the other should be a cokernel presentation. Thus we obtain the result.

(b)
Next we assume that ${\calV}_s^{\bullet}$ is not a sheaf for all $s \in S$. 
Then Fact~\ref{fact:yoshioka-thm2-1} implies that 
there is a numerical solution $(u_1,u_2,k_1,k_2)$ such that 
for a general $s \in S$,
$H^{-1}({\calV}_s^{\bullet})$ and
$H^0({\calV}_s^{\bullet})$ are semi-homogeneous sheaves
with $v(H^{-1}({\calV}_s^{\bullet}))=k_1 u_1$
and $v(H^0({\calV}_s^{\bullet}))=k_2 u_2$ respectively. 
Applying Fact~\ref{fact:yoshioka-lem2-2} to the family of complexes
${\calV}^{\bullet}$ successively,
we find an equivalence ${\calF}\colon\bfD(X)\to \bfD(X)$ 
or $\bfD(X)\to \bfD(X)_\op$ such that ${\calF}({\calV}_s^\bullet)\cong G[k]$, 
where $G$ is a stable sheaf, $k\in\bbZ$ and $v(G[k])=v$. 
Then by the proof of Lemma \ref{lem:FM-sheaf-criterion}, 
we have the kernel and the cokernel presentation of $G$.

(2)
Assume that the numerical equation of $v$ has a unique solution  
$(v_1,v_2,\ell_1,\ell_2)$. 
We have a complex 
${\calV}^{\bullet}=\bigl[{\calV}^{-1} \to {\calV}^0\bigr]$
associated to this unique solution.
This is a semi-homogeneous presentation, 
since otherwise we would have another numerical solution 
by Fact~\ref{fact:yoshioka-thm2-1}.
\end{proof}

\section{Applications of the semi-homogeneous presentation}\label{section:appli}

In this section we give some applications of the semi-homogeneous presentations. In particular Mukai's Conjecture~\ref{conj:mukai1'} is proved, as mentioned in the introduction. We fix an abelian surface $X$ and suppose $\NS(X)=\bbZ H$ in order to use Theorem \ref{thm:presentation}.

\subsection{An old conjecture of Mukai}
\label{subsec:mukai-conj}

In this subsection we give an answer to the following Conjecture~\ref{conj:mukai1'}, originally proposed by Mukai \cite[{Conjecture $1'$}]{Mukai:1980}.

\begin{conj}\label{conj:mukai1'}
Let $X$ be an abelian surface with $\NS(X)=\bbZ H$ and $v$ be a Mukai vector with $\ell\seteq \mpr{v^2}/2>0$. Suppose that $v$ has at least one solution of the numerical equation \eqref{eq:numerical_equation}.
\begin{enumerate}
\item
Among the numerical solutions for $v$, take the solution $(v_1,v_2,\ell_1,\ell_2)$ such that $\rk(v_i)$ is minimum, where the index $i\in\{1,2\}$ is determined by $\ell_i=\ell$. Then for a general member $\bigl[E^{-1}\xrightarrow{f} E^0\bigr]$ of $\frakM^{+}(v_1,v_2,\ell_1,\ell_2)$, $f$ is either surjective or injective. 
\item
In the situation of (1), the kernel or cokernel of $f$ is stable.
\item
A general member of $M_X^H(v)$ has a semi-homogeneous presentation corresponding to the numerical solution $(v_1,v_2,\ell_1,\ell_2)$ of $v$ such that $\rk(v_i)$ is minimum, where the index $i\in\{1,2\}$ is determined by $\ell_i=\ell$.
\end{enumerate}
\end{conj}

\begin{thm}\label{thm:mukai-conj1'}
Conjecture \ref{conj:mukai1'} is true.
\end{thm}

\begin{proof}
(1)
If $v=e^{d H}(1,0,-\ell)$ or $v=e^{d H}(\ell,0,-1)$ with $d\in\bbZ$,
then we have obvious solutions
$v=-\bigl(\ell e^{d H}(0,0,1)-e^{d H}(1,0,0)\bigr)$ 
or $v=\ell e^{d H}(1,0,0)-e^{d H}(0,0,1)$, and the claim holds.   

Hence we may assume that $v \ne e^{d H}(1,0,-\ell)$ nor 
$v \ne e^{d H}(\ell,0,-1)$.
Let $(v_1,v_2,\ell_1,\ell_2)$ be any numerical solution of $v$.
Let ${\calV}^{\bullet}:{\calV}^{-1} \to {\calV}^0$ 
be the family of simple complexes
on $X$ in Proposition~\ref{prop:univ-cpx}~(2).
Namely, we have  
${\calV}_s^{-1} \in M_X^H(\ell_1 v_1)^{ss}$ and 
${\calV}_s^0 \in M_X^H(\ell_2 v_2)^{ss}$ for 
$s \in {\frakS}\seteq\frakM^{+}(v_1,v_2,\ell_1,\ell_2)$.

In order to prove the statement, it is sufficient to 
prove the following:
if ${\calV}_s^{\bullet}[k]$, $k=-1,0$ is not a sheaf 
for all $s \in {\frakS}$, 
then there is another numerical solution
$(v_1',v_2',\ell_1',\ell_2')$ of $v$ 
such that $\rk(v_i') < \min\{\rk(v_1),\rk(v_2) \}$ 
for $i=1,2$.

Assume that ${\calV}_s^{\bullet}[k]$, $k=-1,0$ is not a sheaf 
for all $s \in \frakS$. 
By Fact 2.12, there is another numerical solution
$(v_1',v_2',\ell_1',\ell_2')$ such that
$H^{-1}({\calV}_s^{\bullet}) \in M_X^H(\ell_1' v_1')^{ss}$
and 
$H^0({\calV}_s^{\bullet}) \in M_X^H(\ell_2' v_2')^{ss}$
for a general $s \in \frakS$.
By Lemma \ref{lem:simple-complex2},
$\Phi_{X \to Y_1}^{{\bfE}_1^{\vee}[1]}({\calV}^{\bullet})$
is a family of simple complex with Mukai vector $(\ell_2,0,-\ell_1)$.
Since 
$\calV^{\bullet}_s \cong 
 \bigl[H^0(\calV^{\bullet}_s) \to H^{-1}(\calV^{\bullet}_s)[2]\bigr][-1]$,
we have
\begin{align*}
\Phi_{X \to Y_1}^{{\bfE}_1^{\vee}}({\calV}_s^{\bullet})
\cong 
& 
\bigl[
 \Phi_{X \to Y_1}^{{\bfE}_1^{\vee}}(H^0(\calV^{\bullet}_s)) 
 \to
 \Phi_{X \to Y_1}^{{\bfE}_1^{\vee}}(H^{-1}(\calV^{\bullet}_s))[2]
\bigr][-1]
\\
\cong &
\bigl[\Phi_{X \to Y_1}^{{\bfE}_1^{\vee}}(H^0(\calV^{\bullet}_s)) \to
      \Phi_{X \to Y_1}^{{\bfE}_1^{\vee}[2]}(H^{-1}(\calV^{\bullet}_s))
\bigr][-1].
\end{align*}
Since $\mu(v_1')<\mu(v_1)<\mu(v_2')$,
$\Phi_{X \to Y_1}^{{\bfE}_1^{\vee}[2]}\bigl(H^{-1}(\calV^{\bullet}_s)\bigr)$ 
and
$\Phi_{X \to Y_1}^{{\bfE}_1^{\vee}}\bigl(H^0(\calV^{\bullet}_s)\bigr)$
 are semi-homogeneous vector bundles with
Mukai vector $\ell_1' u_1$ and $\ell_2' u_2$ respectively. 
Here we defined 
$u_1 \seteq \Phi_{X \to Y_1}^{{\bfE}_1^{\vee}}(v_1')$ and
$u_2 \seteq \Phi_{X \to Y_1}^{{\bfE}_1^{\vee}}(v_2')$.

Now we show the inequality $\rk(v_1)-\rk(v_i')>0$.
We note that  
$0=\mu\bigl(\Phi_{X \to Y_1}^{{\bfE}_1^{\vee}}(v_2)\bigr)<
\mu(u_2) \leq \mu\bigl( {\bfE}_1|_{\{ x \} \times Y_1}^{\vee}\bigr) 
<\mu(u_1)$. 
If $u_2=v\bigl({\bfE}_1|_{\{ x \} \times Y_1}^{\vee}\bigr)$ and
$\rk(u_2)=\rk(v_1)=1$, then $\rk(v_1')=1$ and $v_2'=(0,0,1)$,
which means that $v=e^{d H}(1,0,-\ell)$, or $e^{d H}(\ell,0,-1)$, 
contradicting the assumption.
Therefore $u_2 \ne v({\bfE}_1|_{\{ x \} \times Y_1}^{\vee})$ or
$\rk(u_2) \ne 1$.
Since 
\begin{align*}
\rk(v_1)-\rk(v_i')
&=
\rk\bigl(\Phi_{Y_1 \to X}^{{\bfE}_1[2]}((0,0,1))\bigr)-
\rk\bigl(\Phi_{Y_1 \to X}^{{\bfE}_1[2]}(u_i)\bigr)
\\
&=
-\mpr{v\bigl({\bfE}_1|_{\{ x \} \times Y_1}^{\vee}\bigr),(0,0,1)}
+\mpr{v\bigl({\bfE}_1|_{\{ x \} \times Y_1}^{\vee}\bigr),u_i},
\end{align*}
the claim follows from Lemma \ref{lem:rank} below.

Next we deal with the inequality $\rk(v_2)-\rk(v_i')>0$.
By Lemma \ref{lem:simple-complex2},
$\Phi_{X \to Y_2}^{{\bfE}_2^{\vee}[1]}({\calV}^{\bullet})$
is a family of stable sheaves with Mukai vector $(\ell_1,0,-\ell_2)$.
Since 
$\calV^{\bullet}_s \cong 
 \bigl[H^0(\calV^{\bullet}_s) \to H^{-1}(\calV^{\bullet}_s)[2]\bigr][-1]$,
we have 
\begin{align*}
\Phi_{X \to Y_2}^{{\bfE}_2^{\vee}[1]}({\calV}_s^{\bullet})
\cong & 
\bigl[\Phi_{X \to Y_2}^{{\bfE}_2^{\vee}[1]}(H^0(\calV^{\bullet}_s))[-1] \to
\Phi_{X \to Y_2}^{{\bfE}_2^{\vee}[1]}(H^{-1}(\calV^{\bullet}_s))[1]\bigr]\\
\cong &
\bigl[\Phi_{X \to Y_2}^{{\bfE}_2^{\vee}}(H^0(\calV^{\bullet}_s)) \to
\Phi_{X \to Y_2}^{{\bfE}_2^{\vee}[2]}(H^{-1}(\calV^{\bullet}_s))\bigr].
\end{align*}
Since $\mu(v_1')<\mu(v_2)<\mu(v_2')$,
$\Phi_{X \to Y_2}^{{\bfE}_2^{\vee}[2]}(H^{-1}({\calV}^{\bullet}_s))$ and
$\Phi_{X \to Y_2}^{{\bfE}_2^{\vee}}(H^0({\calV}^{\bullet}_s))$
 are semi-homogeneous vector bundles with
Mukai vector $\ell_1' u'_1$ and $\ell_2' u'_2$
respectively, where
$u'_1 \seteq \Phi_{X \to Y_2}^{{\bfE}_2^{\vee}}(v_1')$ and
$u'_2 \seteq \Phi_{X \to Y_2}^{{\bfE}_2^{\vee}}(v_2')$.
Note that  
$\mu(u_2') \le \mu\bigl( {\bfE}_2|_{\{ x \} \times Y_2}^{\vee}\bigr) 
<\mu(u_1')<0$.
Then we can prove the inequality $\rk(v_2)-\rk(v_i')>0$ using 
Lemma \ref{lem:rank} below and a similar argument as above.

(2), (3)
These are the consequences of (1), Proposition~\ref{prop:cpx-stable}  
and the irreducibility of  $M_X^H(v)$ (Fact \ref{fact:moduli}). 
\end{proof}

\begin{lem}\label{lem:rank}
Let $Y$ be an abelian surface with $\NS(Y)=\bbZ H'$ and 
$w_i\seteq (r_i,d_i H',d_i^2 (H'^2)/(2r_i))$ ($i=1,2,3$) 
be primitive and isotropic Mukai vectors on $Y$ such that 
$r_i>0$,
$\langle w_1,w_3 \rangle =-1$,
$\mu(w_1) \leq \mu(w_2) \leq \mu(w_3)$
and $d_1 d_3>0$.
Then we have
$-\langle w_2,(0,0,1) \rangle+\langle w_2,w_i \rangle \geq 0$.
If the equality holds, then $\rk w_2=1$ and
$w_j=w_2$, where $\{i,j\}=\{1,3\}$.
\end{lem}

\begin{proof}
For $i,j$ with $\{ i, j \}=\{1,3\}$,
\begin{align*}
-\langle w_2,(0,0,1) \rangle+\langle w_2,w_i \rangle 
= & r_2-r_2 r_i 
\left(\dfrac{d_2}{r_2}-\dfrac{d_i}{r_i}\right)^2 \dfrac{(H'^2)}{2} \\
\ge & 
r_2-r_2 r_i 
\left(\dfrac{d_j}{r_j}-\dfrac{d_i}{r_i}\right)^2 \dfrac{(H'^2)}{2} 
=
r_2+\frac{r_2}{r_j}\langle w_i,w_j \rangle 
=
r_2\left(1-\dfrac{1}{r_j}\right) \geq 0.
\end{align*}
If the equality holds, then
$\mu(w_2)=\mu(w_j)$ and $r_j=1$.
Hence $w_j=w_2$ and $\rk(w_2)=1$. 
\end{proof}

\subsection{Sheaf-preserving Fourier-Mukai transforms}
\label{subsec:sheaf-presen}

Since FMTs are defined in the derived category, we cannot expect that a sheaf is transformed into a sheaf in general. In Lemma~\ref{lem:sheaf-criterion} of the previous section, we have already encountered such a phenomenon. Similar statements are dealt in the next proposition.

\begin{prop}\label{prop:sheaf-consv}
Suppose that a coherent sheaf $E$ has both kernel and cokernel presentations.
\begin{align*}
0 \to  E_1 \to  E_2 \to  E \to  0,\quad
0 \to  E   \to  F_1 \to  F_2 \to  0.
\end{align*}
Let $\bfE$ be a universal family of simple semi-homogeneous sheaves. Let $\mu\seteq\mu(\bfE_y)$, $y\in Y$. The image of $E$ by $\Phi\seteq\Phi_{X\to Y}^{\bfE^\vee}$ is classified as follows.

\begin{center}
\begingroup
\renewcommand{\arraystretch}{1.2}
\begin{tabular}{l | l l}
range of slope           &image of $E$ &non-vanishing
\\[-3pt]
                         &             &cohomology sheaf
\\
\hline
$\phantom{M M m m}\mu<\mu(E_1)$&genuine sheaf&$\Phi^0(E)$\\
$\mu(E_1)\le\mu<\mu(E_2)$   &2-term cpx.  &$\Phi^0(E)$, $\Phi^1(E)$\\
$\mu(E_2)\le\mu<\mu(F_1)$   &genuine sheaf&$\Phi^1(E)$\\
$\mu(F_1)\le\mu<\mu(F_2)$   &2-term cpx.  &$\Phi^1(E)$, $\Phi^2(E)$\\
$\mu(F_2)\le\mu$            &genuine sheaf&$\Phi^2(E)$
\end{tabular}
\end{center}
\endgroup
\end{prop}
\begin{proof}
This is the consequence of WIT for semi-homogeneous sheaves (Proposition~\ref{prop:semihom-fmt}). Since the proof is similar to that of Lemma~\ref{lem:FM-sheaf-criterion}, we omit the detail. 
\end{proof}

\begin{eg}
Suppose $\NS(X)=\bbZ H$ and put $n\seteq(H^2)/2\in\bbZ_{>0}$. 
We can give the condition 
when a sheaf is converted into a complex 
under the FMT $\Phi\seteq \Phi^\bfP_{X\to \widehat{X}}$, 
where $\bfP$ is the Poincar\'{e} bundle of $X$. 
Assume that the Mukai vector $v$ with $\mpr{v^2}>0$ has a numerical solution. 
Then for a general member $E\in M_X^H(v)$, 
the image of $E$ under $\Phi$ is a two-term complex 
only if $v(E)$ can be written by one of the following four forms.
\begin{align*}
\begin{array}{l l l l l l}
 \text{(a)}&(\ell d^2n-1,\ell d H,\ell) &\mbox{with}\ d>0,
&\text{(b)}&(d^2n-\ell,d H,1)           &\mbox{with}\ d>0,\\
 \text{(c)}&(\ell-d^2n,-d H,-1)         &\mbox{with}\ d>0,
&\text{(d)}&(1,0,-\ell).
\end{array}
\end{align*}
\end{eg}
\begin{proof}
By Theorem \ref{thm:presentation}, we can assume that $E$ has a semi-homogeneous presentation. We divide the argument into two cases according to the kind of the presentation.

(1)
Assume that $E$ has a cokernel presentation
\begin{align*}
0\to E_1\to E_2\to E\to 0
\end{align*}
and that the image of $E$ by $\Phi$ is not a sheaf under any shift.

(i)
First consider the case 
\begin{align*}
v(E_1)=v_1=(r_1,d_1H,a_1),\quad v(E_2)=\ell v_2=\ell(r_2,d_2H,a_2),
\end{align*}
where $v_1$ and $v_2$ are primitive. Then 
\begin{align*}
&\mpr{v_1,v_2}=-1\iff 2d_1d_2n=r_1a_2+r_2a_1-1,\\
&\mpr{v_1^2}=0\iff d_1^2n=r_1a_1,\quad
 \mpr{v_2^2}=0\iff d_2^2n=r_2a_2.
\end{align*}
By the result of Proposition \ref{prop:sheaf-consv} 
and the equality $\mu(\bfP_y)=0$, 
we need to determine the condition for $\mu(E_1) \le 0 <\mu(E_2)$. 
If $\mu(E_1)\neq0$, then we have
\begin{align*}
r_1a_2+r_2a_1-1<0,\quad  r_1a_1>0\quad \text{and}\quad r_2a_2>0.
\end{align*}
These inequities have no solution, so that we may assume $\mu(E_1)=0<\mu(E_2)$. This condition yields
\begin{align*}
r_1>0,\quad d_1=0,\quad  d_2\ge0\quad \text{and}\quad r_2\ge 0.
\end{align*}
These ineqalities have only one solution:
\begin{align*}
v_1=(1,0,0),\quad v_2=(d^2n,d H,1),\quad d\ge0.
\end{align*}
{}From the condition $\rk(E)>0$ we have $d>0$. 
Thus we obtained the first solution (a).

(ii)
Next consider the case 
\begin{align*}
v(E_1)=\ell v_1=\ell (r_1,d_1H,a_1),\quad v(E_2)=v_2=(r_2,d_2H,a_2).
\end{align*} 
By the result of Proposition \ref{prop:sheaf-consv}, we need to determine the condition for $\mu(E_1)\le0<\mu(E_2)$. As in the first case, $\mu(E_1)<0<\mu(E_2)<\infty$ has no solution. Then we have to consider the case $\mu(E_1)=0<\mu(E_2)$. This case has the next solution
\begin{align*}
v_1=(1,0,0),\quad v_2=(d^2n,d H,1)\quad d\ge0,
\end{align*}
and we can exclude the case $d=0$ by the condition $\rk(E)>0$. Hence the second solution (b) is obtained.

(2)
Assume that $E$ has a kernel presentation
\begin{align*}
0\to E\to E_2\to E_1\to 0
\end{align*}
and the image of $E$ is a two-term complex. We can put $v(E_1)=\ell_1v_1$ and $v(E_2)=\ell_2v_2$ with $(\ell_1,\ell_2)=(1,\ell)$ or $(\ell,1)$. The condition we have to consider is $\mu(E_1)\le0<\mu(E_2)$. Then a simple calculation as in (1) shows that the solutions are
\begin{align*}
&v_1=(d^2n,d H,1),\quad v_2=(1,0,0),\quad d\ge0.
\end{align*}
Considering the condition $\rk(E)>0$, we obtain the third solution (c) and the last one (d).
\end{proof}

\begin{rmk}
In the cases (a) and (b), 
one may easily find 
that $\calD_{\widehat{X}}\circ\Phi$ induces a birational map 
$M_X^H(v)\cdots\to X\times\Hilb{\ell}{\widehat{X}}$. 
Thus we recover the result 
of \cite[Theorem 3.14.~(1)]{Yoshioka:2009:stabilityII}.
\end{rmk}

As a consequence of Theorem \ref{thm:presentation} 
and Proposition \ref{prop:sheaf-consv}, we have the next proposition.

\begin{prop}
Suppose that $\NS(X)=\bbZ H$, $(H^2)=2n$ 
and that $v$ is a primitive Mukai vector. 
Then for any universal family $\bfE$ of 
simple semi-homogeneous sheaves over $Y=M_X^H(v')$ 
with $v'=(r',d' H,a')$ satisfying
\begin{align}
\label{ineq:cpx_trans}
\dfrac{d'}{r'}<\dfrac{d}{r}-\dfrac{\sqrt{r+1}}{r}\sqrt{\dfrac{\ell}{n}},
\end{align}
or
\begin{align*}
\dfrac{d'}{r'}>\dfrac{d}{r}+\dfrac{\sqrt{r+1}}{r}\sqrt{\dfrac{\ell}{n}},
\end{align*}
the image $\Phi^{\bfE}_{X\to Y}(E)$ of a general member $E\in M_X^H(v)$ 
is a stable sheaf up to shift.
\end{prop}

\begin{proof}
We only show the proof for the case \eqref{ineq:cpx_trans}, 
since the other case can be shown quite similarly.

By Fact~\ref{fact:yoshioka-thm2-1}, 
we may assume that $v$ has a numerical solution.
By Theorem \ref{thm:presentation}, 
we can assume that $E$ has a semi-homogeneous presentation
\begin{align*}
\begin{array}{l l l l l l l l l l l}
(1)&0&\to&E&\to &E_1&\to &E_2&\to&0,&\quad\mbox{or}\\
(2)&0&\to&E_2&\to &E_1&\to &E&\to&0,&\quad\mbox{or}\\
(3)&0&\to&E&\to &E_2&\to &E_1&\to&0,&\quad\mbox{or}\\
(4)&0&\to&E_1&\to &E_2&\to &E&\to&0&
\end{array}
\end{align*}
with $v(E_1)=\ell v_1$, $v(E_2)=v_2$.
By Lemma \ref{lem:pmvpr} proved later, we can set 
\begin{align*}
&v_1=(p_1^2r_1,p_1q_1H,q_1^2r_2),
\quad
v_2=(p_2^2r_2,p_2q_2H,q_2^2r_1),\\
&r_1r_2=n,\quad p_1>0,\quad  p_1q_2r_1-p_2q_1r_2=1.
\end{align*}

The cases (1) and (3) follows from Proposition \ref{prop:sheaf-consv}, since for each case we have 
\begin{align*}
\mu(v')<\mu(v)<\mu(v_1)<\mu(v_2),\quad
\mu(v')<\mu(v)<\mu(v_2)<\mu(v_1).
\end{align*}

For the case (2), we estimate the value of $\mu(v)-\mu(v_2)$.
{}From $v=\ell v_1-v_2$ 
we have $r=\ell p_1^2r_1-p_2^2r_2$ and $d=\ell p_1q_1-p_2q_2$. 
Then we have
\begin{align*}
\mu(v)-\mu(v_2)
&=\dfrac{\ell p_1q_1-p_2q_2}{\ell p_1^2r_1-p_2^2r_2}-\dfrac{q_2}{p_2r_2}
=\dfrac{-\ell}{r}\dfrac{p_1}{p_2r_2}.
\end{align*}
Notice that $p_2<0$ holds, since 
\begin{align*}
\mu(v_2)<\mu(v_1)
&\iff
\dfrac{q_2}{p_2r_2}<\dfrac{q_1}{p_1r_1}
\iff
0<\dfrac{q_1}{p_1r_1}-\dfrac{q_2}{p_2r_2}
=\dfrac{p_2q_1r_2-p_1q_2r_1}{p_1p_2r_1r_2}\\
&\iff
0<-\dfrac{1}{p_1p_2n}\iff p_2<0.
\end{align*}
We also find from $r=\ell p_1^2r_1-p_2^2r_2$ that $p_1^2/p_2^2r_2^2=(1+r/p_2^2r_2)/\ell n$, which leads to $p_1^2/p_2^2r_2^2\le (r+1)/\ell n$. The equality holds if and only if $p_2^2r_2=1$. Hence we have
\begin{align*}
\mu(v)-\mu(v_2)\le\dfrac{\sqrt{r+1}}{r}\sqrt{\dfrac{\ell}{n}}.
\end{align*}
Therefore the conclusion follows from Proposition \ref{prop:sheaf-consv}, since
the assumption \eqref{ineq:cpx_trans} means
\begin{align*}
\mu(v')<\mu(v)-\dfrac{\sqrt{r+1}}{r}\sqrt{\dfrac{\ell}{n}}
       \le\mu(v_2)
\end{align*}
for any choice of $v_1$ and $v_2$.

For the case (4), we have $r=p_2^2r_2-\ell p_1^2r_1$ and $d=p_2q_2-\ell p_1q_1$. Then 
\begin{align*}
\mu(v)-\mu(v_1)
&=\dfrac{p_2q_2-\ell p_1q_1}{p_2^2r_2-\ell p_1^2r_1}-\dfrac{q_1}{p_1r_1}
=\dfrac{1}{r}\dfrac{p_2}{p_1r_1}.
\end{align*}
In this case we have $p_2>0$ and $p_2^2/p_1^2r_1^2=(\ell+r/p_1^2r_1)/n$, which means that $(r+\ell)p_2^2/n p_1^2r_1^2>\ell/n$. Hence we have
\begin{align*}
\mu(v)-\mu(v_1)\le\dfrac{1}{r}\sqrt{\dfrac{\ell+r}{n}}.
\end{align*}
Since $\sqrt{\ell(r+1)/n}\ge\sqrt{(\ell+r)/n}$, the assumption \eqref{ineq:cpx_trans} means that
\begin{align*}
\mu(v')<\mu(v)-\dfrac{1}{r}\sqrt{\dfrac{\ell+r}{n}}
       \le\mu-(\mu-\mu_1)=\mu_1.
\end{align*}
Therefore the conclusion holds 
by Proposition \ref{prop:sheaf-consv} and Fact~\ref{fact:yoshioka-lem2-2}~(1).
\end{proof}

\section{Matrix description of cohomological FMT}
\label{section:matrix}

This section is devoted to the description of cohomological FMT via quadratic matrices. The original idea is due to Mukai \cite{Mukai:1979,Mukai:1980}. We fix an abelian surface $X$ with $\NS(X)=\bbZ H$ and set $(H^2)=2n$.

Our study can be summarized in the next diagram.
\begin{align}
\label{diag:groups}
\xymatrix{
\calE(X)     \ar[r]^{\theta}
&
\LieO(B)     \ar[r]_{\sim}^{\beta}
&
\LieO(2,1)   \ar@{<-^{)}}[d]
\\
\ar@{}[rd]|{\circlearrowright}
G            \ar@{>>}[d]_{2:1} \ar@{^{(}->}[r]
&
\SL(2,\bbR)  \ar@{>>}[d]_{2:1} \ar[dr]_{\Ad} \ar[u]^{\alpha}
\ar@{}[r]|{\circlearrowright}
&
\SO(2,1)     \ar@{<-^{)}}[d]
\\
G/\{\pm1\}   \ar@{^{(}->}[r]
&
\PSL(2,\bbR) \ar[r]_{\sim}
&
\SO_0(2,1)   
}
\end{align}
Here $\LieO(2,1),\SO(2,1),\SL(2,\bbR),\PSL(2,\bbR)$ are the familiar  real Lie groups, regarded as subgroups of the group of invertible $3 \times 3$ matrices $\GL(3,\bbR)$. $\SO(2,1)$ has two connected components and $\SO_0(2,1)$ denotes the component of $\SO(2,1)$ including the identity matrix.  We are interested in the set
\begin{align*}
\calE(X)\seteq
\bigcup_{Y\in\FM(X)} \Eq_0(\bfD(Y),\bfD(X)), 
\end{align*}
where $\Eq_0(\bfD(Y),\bfD(X))$ is a subset of $\Eq(\bfD(Y),\bfD(X))$ 
defined in Definition~\ref{dfn:Eq0&theta}~(1). 
Although this set cannot be studied directly, 
we can investigate it through its cohomological representation. 
The map $\theta:\calE(X)\to\LieO(B)$, 
defined in Definition~\ref{dfn:Eq0&theta}~(2), realizes this idea. 
Here $B$ is a lattice of rank three, 
which is essentially the same as the Mukai paring 
on the abelian surface of Picard number one, 
as defined in Definition~\ref{dfn:Sym2&B}. 
The group $\LieO(B)$ consists of the isometry of the lattice $B$. 

The main result of this section, Theorem \ref{thm:bij}, claims the existence of bijection
\begin{align}\label{eq:bij}
\Img \theta \cong G/\{\pm1\},
\end{align}
where $G$ is a subgroup of $\SL(2,\bbR)$ introduced in Definition~\ref{dfn:G}. 
The element of $G$ encodes 
the data of Mukai vectors $\Phi((0,0,1))$ and $\Phi((1,0,0))$. 
The bijection \eqref{eq:bij} indicates 
that the cohomological action of $\Phi$ is uniquely determined 
by the image of the skyscraper sheaf and the structure sheaf.

\subsection{Cohomological correspondence}\label{subsec:matrix:prelim}

In this preliminary subsection we give a description of the Mukai lattice on the Fourier-Mukai partner $Y$ of $X$.

In the paper \cite[\S\S\,1.3]{Yoshioka:2009:stabilityII}, 
one of the author deals with the ample divisor $\widehat{H}$. 
Let us recall the argument. For a while we do not assume $\NS(X)=\bbZ H$. 
Let $v_0=(r_0,\xi_0,a_0)$ be an isotropic Mukai vector 
and set $Y\seteq M_X^H(v_0)$. 
Assume that there exists a universal family $\bfE$ on $Y\times X$. 
Let $x\in X$ and $y\in Y$ be arbitrary closed points. 
Set 
$w_0\seteq v_0(\bfE|_{\{x\}\times Y})
=r_0+\widetilde{\xi}_0+\widetilde{a}_0\rho_Y$, 
$\widetilde{\xi}_0\in\NS(Y)$. 
Then $\Phi(\bfE^\vee|_{X\times\{y\}})=\bbC_y$ and 
$\Phi(\bbC_x)=\bfE|_{\{x\}\times Y}$, where $\Phi\seteq\Phi_{X\to Y}^{\bfE}$. 
Hence for the cohomological FMT, 
we get $\Phi^H(v_0^\vee)=\rho_Y$ and $\Phi^H(\rho_X)=w_0$.

Now we introduce a map $\nu\colon H^2(X,\bbQ) \to H^2(Y,\bbQ)$ by
\begin{align*}
-\nu(D)
\seteq
-\left[\Phi^H\left(D+\dfrac{(D,-\xi_0)}{r_0}\rho_X\right)\right]_1
=\left[p_{Y*}
  \left(\biggl(c_2(\bfE)-\dfrac{r_0-1}{2r_0}c_1(\bfE)^2\biggr)\cup p_X^*(D)
  \right)
 \right]_1,
\end{align*}
where $[\ ]_1$ means the natural projection from $H^{\rm ev}(Y,\bbQ)$ to $H^2(Y,\bbQ)$. Then by \cite[Lemma 1.4]{Yoshioka:2009:stabilityII}, $-r_0\nu(H)$ is represented by an ample divisor on $Y$. Thus $\widehat{H}$ is proportional to $-r_0\nu(H)$.

Using this divisor $\widehat{H}$ we can describe the cohomological behaviour of FMT in the following manner.

\begin{fact}[{\cite[Proposition 1.5]{Yoshioka:2009:stabilityII}}]\label{fact:fmt:cohomology}
\begin{enumerate}
\item
Every Mukai vector $v=(r,\xi,a)\in H^{\ev}(X,\bbZ)_\alg$ is uniquely written as
\begin{align*}
v=\ell v_0^\vee+\chi\rho_X+d H+D-\dfrac{1}{r_0}(d H+D,\xi_0)\rho_X
\end{align*}
with
\begin{align*}
\ell\seteq\dfrac{r}{r_0}\in\dfrac{1}{r_0}\bbZ,\quad
\chi
\seteq\dfrac{r a_0+(\xi_0,\xi)+r_0a}{r_0}
\in\dfrac{1}{r_0}\bbZ,\quad
&
d\seteq\dfrac{(r_0\xi+r\xi_0,H)}{r_0(H^2)}
\in\dfrac{1}{r_0(H^2)}\bbZ
\end{align*}
and $D\in\NS(X)\otimes_\bbZ\bbQ\cap H^\perp$.
\item
\begin{align}
\label{eq:cohom_fmt}
\begin{split}
&\Phi^H(\ell v_0^\vee+\chi\rho_X+d H+D-\dfrac{1}{r_0}(d H+D,\xi_0)\rho_X))
\\
&=\ell \rho_Y+\chi w_0+d\nu(H)+\nu(D)+\dfrac{1}{r_0}(d\nu(H)+\nu(D),
\widetilde{\xi}_0)\rho_Y,
\end{split}
\end{align}
where $D\in\NS(X)\otimes_\bbZ\bbQ\cap H^\perp$.
\end{enumerate}
\end{fact}

Now we return to the case $\NS(X)=\bbZ H$. Suppose that for an isotropic Mukai vector $v$, $Y\seteq M_X^H(v)\in\FM(X)$ and that there exists a universal family $\bfE$ on $Y\times X$. Set $\NS(Y)=\bbZ \widehat{H}$ and $\Phi\seteq\Phi_{X\to Y}^{\bfE}$. 

\begin{lem}
In this setting we have $-\nu(H)=\widehat{H}$.
\end{lem}
\begin{proof}
We can write $\widehat{H}=-c \nu(H)$ with $c\in\bbQ_{>0}$ as above. Lemma \ref{lem:picard} shows that $(\widehat{H}^2)=2n$. Thus we have $c=1$. (See also \cite[Lemma Appendix A.1]{Yoshioka:2009:stabilityII} for a general treatment.)
\end{proof}

Consider the action of the cohomological FMT $\Phi^H$ on the lattices $H^{\ev}(X,\bbZ)_\alg\to H^{\ev}(Y,\bbZ)_\alg$ of rank three. Take a basis $\langle 1, H, \rho_X\rangle$ of $H^{\ev}(X,\bbZ)_\alg$ and a basis $\langle 1, \widehat{H}, \rho_Y\rangle$ of $H^{\ev}(Y,\bbZ)_\alg$. Then $\Phi^H$ can be written by a $3\times3$ matrix.

\begin{lem}\label{lem:det:positive}
The determinant of $\Phi^H$ is $1$.
\end{lem}
\begin{proof}
By Lemma \ref{lem:picard}, we may put $(H^2)=(\widehat{H}^2)=2n$. 
Set $v_0=(r_0,d_0 H,a_0)$ and 
$w_0=(r_0,\widetilde{d}_0\widehat{H},\widetilde{a}_0)$. 
{}From \eqref{eq:cohom_fmt}, $\Phi^H$ has an effect
\begin{align*}
\ell v_0^\vee+d H+(\chi-\dfrac{2 n d_0}{r_0}d)\rho_X
\mapsto
\chi w_0-d\widehat{H}+(\ell-\dfrac{2n\widetilde{d}_0}{r_0}d)\rho_Y,
\end{align*}
so that if we present $\Phi^H$ by the basis $\langle v_0^\vee,H,\rho_X\rangle$ and $\langle w_0,\widehat{H},\rho_Y\rangle$, then its matrix is
\begin{align*}
M\seteq
\begin{bmatrix}
0 & 2n d_0/r_0 & 1\\ 0 & -1         & 0\\
1 & -2n\widetilde{d}_0/r_0 & 0
\end{bmatrix}.
\end{align*}
On the other hand, the base change $\langle 1,H,\rho_X\rangle \mapsto \langle v_0^\vee,H,\rho_X\rangle$ can be written as
\begin{align*}
P_1\seteq
\begin{bmatrix}
1/r_0 & 0 & 0\\ d_0/r_0 & 1 & 0\\ -a_0/r_0 & 0 & 1
\end{bmatrix}.
\end{align*}
Similarly the base change $\langle w_0,\widehat{H},\rho_Y\rangle \mapsto \langle 1,\widehat{H},\rho_Y\rangle$ can be written as
\begin{align*}
P_2\seteq
\begin{bmatrix}
1/r_0 & 0 & 0\\ -\widetilde{d}_0/r_0 & 1 & 0\\ -\widetilde{a}_0/r_0 & 0 & 1
\end{bmatrix}^{-1}.
\end{align*}
The determinant of the multiplication of these three matrices is $1$. Thus the conclusion holds.
\end{proof}

\begin{rmk}\label{rmk:g33}
For later use, we display the resulting matrix $P_2 M P_1$.
\begin{align}\label{eq:g33}
P_2 M P_1
=
\begin{bmatrix}
a_0 & 2n d_0 & r_0\\ 
(a_0 \widetilde{d}_0 -d_0)/r_0& -1+2n d_0\widetilde{d}_0/r_0  
 & \widetilde{d}_0\\
(-2n d_0 \widetilde{d}_0+r_0+a_0 \widetilde{a}_0 r_0)/r_0 ^2
 & 2 n (b_0 d_0 -\widetilde{d}_0)/r_0 & \widetilde{a}_0
\end{bmatrix}.
\end{align}
\end{rmk}

\subsection{Some calculations for matrix description}
\label{subsec:prelim:matrix}

In this subsection we give elementary remarks on the groups appeared in the diagram \eqref{diag:groups}. We begin with the description of the indefinite orthogonal group $\LieO(2,1)$. The following notation will be used.
\begin{align*}
&I_{2,1}\seteq
 {\rm diag}(1,1,-1)\,\in\GL(3,\bbR),\\
&\LieO(2,1)\seteq
 \{g\in\GL(3,\bbR)\mid {}^t g I_{2,1} g= I_{2,1}\},\\
&\SO(2,1)\seteq
 \SL(3,\bbR)\cap \LieO(2,1).
\end{align*}
The connected components of $\LieO(2,1)$ are described as follows.

\begin{lem}\label{lem:Opq}
There are four connected components of $\LieO(2,1)$, and these are given by

(0) \quad $\{g=(g_{i j})\in\LieO(2,1)\mid \det g=+1,\ g_{3 3}\ge 1\}$,

(1) \quad $\{g=(g_{i j})\in\LieO(2,1)\mid \det g=+1,\ g_{3 3}\le -1\}$,

(2) \quad $\{g=(g_{i j})\in\LieO(2,1)\mid \det g=-1,\ g_{3 3}\ge 1\}$,

(3) \quad $\{g=(g_{i j})\in\LieO(2,1)\mid \det g=-1,\ g_{3 3}\le -1\}$.
\end{lem}

\begin{proof}
This is a well-known fact and we only show the sketch. 
First note that $g\in\LieO(2,1)$ has $\det g=\pm1$. 
Introduce the bilinear form $(\cdot,\cdot)$ on $\bbR^3$ by
\begin{align*}
(x,y)\seteq x_1 y_1+x_2 y_2-x_3 y_3.
\end{align*}
Then for $g\in\GL(3,\bbR)$, 
the condition ${}^t g I_{2,1} g=I_{2,1}$ is equivalent to 
$[(g_i,g_j)]=I_{2,1}$, 
where $g_i$ denotes the $i$-th column of the matrix $g$. 
Hence the condition is rewritten into
\begin{align*}
\mbox{(i)  }    (g_i,g_j)=0\ (1\le i,j\le 3,\ i\neq j),
\quad
\mbox{(ii)  }   (g_i,g_i)=1\ (i=1,2),
\quad
\mbox{(iii)  }  (g_3,g_3)=-1.
\end{align*}
Then (iii) implies $g_{13}^2+g_{23}^2-g_{33}^2=-1$, which yields $g_{33}^2\ge 1$.

Considering the continuous functions $\det g,g_{33}$ on the real Lie group $\LieO(2,1)$, we find that the subsets (0)-(3) belong to mutually different components. Define $I_0\seteq{\rm diag}(1,1,1)$, $I_1\seteq{\rm diag}(-1,1,-1)$, $I_2\seteq{\rm diag}(1,-1,1)$, $I_3\seteq{\rm diag}(1,1,-1)$. Then each $I_j$ belongs to the subset ($j$). One can explicitly construct a path from any element of the subset ($j$) to the element $I_j$. Thus the conclusion follows.
\end{proof}

The component (0) in the above classification will be denoted by $\SO_0(2,1)$, which contains the identity. As is well-known, the adjoint map 
\begin{align}\label{eq:ad}
\begin{array}{c c c l}
\overline{\Ad} \colon 
      & \SL(2,\bbR) & \to     & \SO_0(2,1)\\
      &     g       & \mapsto & (x \in\Liesl(2,\bbR) \mapsto g^{-1} x g)
\end{array}
\end{align}
induces an \emph{anti}-isomorphism $\PSL(2,\bbR)\seteq\SL(2,\bbR)/\{\pm 1\}\simto\SO_0(2,1)$. Here the image $\overline{\Ad}(g)=(x\mapsto  g^{-1} x g)$ is regarded as the matrix presented in the basis 
$\langle 
\begin{bmatrix}1&0\\ 0&-1\end{bmatrix},
\begin{bmatrix}0&1\\ 1& 0\end{bmatrix},
\begin{bmatrix}0&1\\-1& 0\end{bmatrix}
\rangle$ of $\Liesl(2,\bbR)$. 
A straightforward computation yields that for $g=\begin{bmatrix}p& q\\ r&s\end{bmatrix}$ the matrix becomes
\begin{align}\label{eq:ad_g}
\overline{\Ad}(g)=
\begin{bmatrix}
p s+q r& -p q +r s& p q+r s\\
q s-p r&( p^2-q^2-r^2+s^2)/2 & (-p^2+q^2-r^2+s^2)/2\\
q s+p r&(-p^2-q^2+r^2+s^2)/2 & ( p^2+q^2+r^2+s^2)/2
\end{bmatrix}.
\end{align} 
The overline in the notation indicates that we do not apply the usual definition $\Ad(g)\seteq (x \mapsto g  x g^{-1})$, which gives an isomorphism $\PSL(2,\bbR)\simto\SO_0(2,1)$.
The quadratic form used in this correspondence is $x_1^2+x_2^2-x_3^2$, which is different from the Mukai paring. This discrepancy is solved as follows.

\begin{dfn}\label{dfn:Sym2&B}
We define
\begin{align*}
\Sym_2(\bbR)\seteq
\left\{\begin{bmatrix}x&y\\y&z\end{bmatrix}\, \Bigg|\, x,y,z\in \bbR\right\}.
\end{align*}
We also introduce a bilinear form $B$ on $\Sym_2(\bbR)$ as
\begin{align*}
B(X,X')\seteq-\tr(\widetilde{X}X')=2y y'-(x z'+z x'),\quad
X =\begin{bmatrix}x &y \\y &z \end{bmatrix},\ 
X'=\begin{bmatrix}x'&y'\\y'&z'\end{bmatrix}
\in\Sym_2(\bbR).
\end{align*}
Here we denoted $\widetilde{X}\seteq\begin{bmatrix}z&-y\\-y&x\end{bmatrix}$ for $X=\begin{bmatrix}x&y\\y&z\end{bmatrix}$.
\end{dfn}

The signature of $(\Sym(2,\bbR),B)$ is $(2,1)$. We define an action of $\GL(2,\bbR)$ on $\Sym_2(\bbR)$ by
\begin{align}
\label{eq:G_action}
X\mapsto X\cdot g \seteq 
 {}^t g X g,\quad X\in\Sym_2(\bbR),\ g \in\GL(2,\bbR).
\end{align}
Since $\widetilde{X}=(\det X)X^{-1}$, we see that
$\widetilde{{}^t g X g}
=(\det g)^2 g^{-1}\cdot\widetilde{X}\cdot{}^{t} g^{-1}$. Hence
\begin{align*}
B({}^t g X g,{}^t g X' g)=(\det g)^2B(X,X').
\end{align*}
In particular, $\SL(2,\bbR)$ preserves the bilinear form $B$. Thus we get a map $\alpha:\SL(2,\bbR)\to \LieO(B)$, where $\LieO(B)$ is the set of linear transformations preserving the paring $B$.

To present the homomorphism $\alpha$ explicitly, we choose a basis 
$\langle 
\begin{bmatrix}1&0\\ 0& 0\end{bmatrix},
\begin{bmatrix}0&1\\ 1& 0\end{bmatrix},
\begin{bmatrix}0&0\\ 0& 1\end{bmatrix}
\rangle$ of $\Sym_2(\bbR)$. Then the paring $B$ corresponds to the symmetric matrix 
\begin{align*}
J\seteq
\begin{bmatrix}
 0& 0&-1\\
 0& 2& 0\\
-1& 0& 0
\end{bmatrix}.
\end{align*}
Using the same basis, we can regard the element $M\in \LieO(B)$ as a $3\times 3$ matrix. There is an isomorphism $\beta:\LieO(B)\to\LieO(2,1)$ defined by
\begin{align}\label{eq:B_21}
\beta \colon \LieO(B)\ni M \mapsto P^{-1} M P\in\LieO(2,1),
\end{align}
where we set
\begin{align}\label{eq:P}
P\seteq
\dfrac{1}{\sqrt{2}}
\begin{bmatrix} -1& 0&1\\  0& 1&0\\  1& 0&1 \end{bmatrix}.
\end{align}
Note that the map $\beta$ is well-defined since ${}^t P J P=I_{2,1}$ holds.

\begin{lem}
The following diagram is commutative.
\begin{align*}
\xymatrix{\ar@{}[rd]|{\circlearrowright}
 \SL(2,\bbR) \ar[r]^\alpha \ar@{>>}[d]_{\overline{\Ad}} 
&\LieO(B)    \ar[d]^\beta_{\wr}  \\
 \SO_0(2,1)  \ar@{^{(}->}[r] 
&\LieO(2,1) }
\end{align*}
In particular the map $\beta\circ\alpha\colon\SL(2,\bbZ)\to\LieO(2,1)$ is a surjective anti-homomorphism $\SL(2,\bbZ)\to\SO_0(2,1)$.
\end{lem}

\begin{proof}
For $g=\begin{bmatrix}p&q\\r&s\end{bmatrix}\in\SL(2,\bbR)$, the image $\alpha(g)=(X\mapsto {}^t g X g)$ in \eqref{eq:G_action} is written in a matrix form via the basis $\langle 
\begin{bmatrix}1&0\\ 0& 0\end{bmatrix},
\begin{bmatrix}0&1\\ 1& 0\end{bmatrix},
\begin{bmatrix}0&0\\ 0& 1\end{bmatrix}
\rangle$ of $\Sym_2(\bbR)$ as
\begin{align}\label{eq:ad_m}
\alpha(g)\seteq
\begin{bmatrix}
p^2 &  2 p r& r^2\\
p q &p s+q r& r s\\
q^2 &  2 q s& s^2
\end{bmatrix}.
\end{align}
It satisfies ${}^t \alpha(g) J \alpha(g)=J$. Then a direct calculation using \eqref{eq:ad_g} and \eqref{eq:ad_m} yields the desired equality $\overline{\Ad}(g)=P^{-1} \alpha(g) P$.
\end{proof}


We continue the study on $\Sym_2(\bbR)$. Set
\begin{align*}
\Sym_2(\bbZ,n)\seteq
\left\{
 \begin{bmatrix}x&y\sqrt{n}\\y\sqrt{n}&z\end{bmatrix}
 \, \Bigg|\, x,y,z\in \bbZ
\right\}
\subset\Sym_2(\bbR).
\end{align*}
Then $B$ is an integral bilinear form on $\Sym_2(\bbZ,n)$. It is nothing but an analogue of the Mukai pairing for abelian surfaces of Picard number one.

\begin{lem}\label{lem:group}
Let $x,y,z,w\in\bbR$ be algebraic numbers such that
\\
(i)   $x^2,x y/\sqrt{n},y^2\in\bbZ$.\\
(ii)  $z^2,z w/\sqrt{n},w^2\in\bbZ$.\\
(iii) $x w-y z=\pm1$.\\
Then there are unique integers $a,b,c,d$ and $r,s>0$ such that
\begin{align}\label{eq:matrix}
 \begin{bmatrix}x&y\\z&w\end{bmatrix}
=\begin{bmatrix}a\sqrt{r}&b\sqrt{s}\\c\sqrt{s}&d\sqrt{r}\end{bmatrix}.
\end{align}
with $r s=n$ and $a d r-b c s=\pm1$.
\end{lem}
\begin{proof}
The case $x y=0$ is trivial, since then we can write
\begin{align*}
\begin{bmatrix}x&y\\z&w\end{bmatrix}
=\pm\begin{bmatrix}0&1\\-1&d\sqrt{n}\end{bmatrix}
\ \mbox{or}\ 
\pm\begin{bmatrix}1&0\\c\sqrt{n}&1\end{bmatrix},\quad
c,d\in\bbZ.
\end{align*}
The case $z w=0$ can be treated similarly. 
Hence we may assume that $x y z w\neq0$.

We set $a\seteq\sgn(x)\gcd(x^2,x y/\sqrt{n})$ and 
$b\seteq\sgn(y)\gcd(y^2,x y/\sqrt{n})$. 
The condition (iii) yields
\begin{align}\label{eq:square}
x^2w^2+y^2z^2-2\dfrac{x y}{\sqrt{n}}\dfrac{z w}{\sqrt{n}}n=1.
\end{align}
Hence $\gcd(a,y^2)=1$, $\gcd(b,x^2)=1$, and $\gcd(a,b)=1$. 
We can write $x^2=a^2r$, $y^2=b^2s$ and $x y/\sqrt{n}=ab\lambda$, 
where $\lambda,r,s\in\bbZ$. 
Then $x^2 y^2 =(x y/\sqrt{n})^2 n \iff a^2 r b^2 s=a^2 b^2\lambda^2 n$ 
implies that $r s=\lambda^2 n$. 
Since $|a|=\gcd(x^2,x y/\sqrt{n})=\gcd(a^2r,ab\lambda)$, 
we have $\gcd(a r,b\lambda)=1$. 
In the same way, we get $\gcd(b s,a \lambda)=1$.
Therefore we should have $\lambda=\pm 1$  and $r s=n$. 
A similar argument using (ii) and (iii) works on $z,w$. 

In the consequence we have 
\begin{align*}
\begin{bmatrix}x&y\\z&w\end{bmatrix}
=\begin{bmatrix}a\sqrt{r}&b\sqrt{s}\\c\sqrt{t}&d\sqrt{u}\end{bmatrix},
\quad
a,b,c,d\in\bbZ,\; r,s,t,u\in\bbZ_{>0},\; r s=t u=n. 
\end{align*}
Then \eqref{eq:square} becomes $a^2d^2r u+b^2c^2s t-2a b c d n=1$. 
By $r s=t u=n$, we have $\gcd(s,u)=\gcd(t,r)=1$. 
Hence we have $s=t$ and $r=u$. 
Then (iii) is equivalent to the requirement $a d r-b c s=\pm1$.

Finally we prove the uniqueness. 
By $a d r-b c s=\pm1$ 
we have $\gcd(x^2,x y/\sqrt{n})=|a| \gcd(a r,b)=|a|$. 
Since the sign of $a$ should equal to that of $x$, 
$a$ is uniquely determined. 
In the same way, $b,c,d$ are uniquely determined. 
Then $r$ and $s$ are characterized as $s=x^2/a^2$ and $r=y^2/b^2$. 
Therefore the uniqueness is shown.
\end{proof}

The condition appeared in Lemma \ref{lem:group} defines a subgroup of $\SL(2,\bbR)$.

\begin{dfn}\label{dfn:G}
We set
\begin{align*}
&\widehat{G}\seteq
\left\{
\begin{bmatrix}x&y\\z&w\end{bmatrix}\in\GL(2,\bbR)\, \Bigg|\, 
x^2,y^2,z^2,w^2,\dfrac{x y}{\sqrt{n}},\dfrac{z w}{\sqrt{n}}\in \bbZ\right\},
\\
&G\seteq \widehat{G}\cap \SL(2,\bbR).
\end{align*}
\end{dfn}

\begin{lem}
(1)
$G$ (resp. $\widehat{G}$) is a subgroup of $\SL(2,\bbR)$ 
(resp. $\GL(2,\bbR)$).

(2)
\begin{align*}
&\widehat{G}= G \rtimes 
\left \langle
\begin{bmatrix} 1 & 0 \\ 0 & -1 \end{bmatrix}
\right \rangle,
\\
&\begin{bmatrix} 1 & 0 \\ 0 & -1 \end{bmatrix}
\begin{bmatrix} x & y \\ z & w  \end{bmatrix}
\begin{bmatrix} 1 & 0 \\ 0 & -1 \end{bmatrix}
=
\begin{bmatrix} x & -y\\ -z & w \end{bmatrix},\quad
\begin{bmatrix} x & y \\  z & w \end{bmatrix} 
\in \GL(2,\bbR).
\end{align*}

(3)
$\widehat{G}$ preserves $\Sym_2(\bbZ,n)$ 
under the action \eqref{eq:G_action} of $\GL(2,\bbR)$ on $\Sym_2(\bbR)$.
\end{lem}

\begin{proof}
We note that
\begin{align}
\label{eq:multiplication}
\begin{bmatrix}a\sqrt{r}  &b\sqrt{s}  \\c\sqrt{s}  &d\sqrt{r}\end{bmatrix}
\begin{bmatrix}a'\sqrt{r'}&b'\sqrt{s'}\\c'\sqrt{s'}&d'\sqrt{r'}\end{bmatrix}
=
\begin{bmatrix}
 a a'\sqrt{r r'}+b c'\sqrt{s s'}
&a b'\sqrt{r s'}+b d'\sqrt{s r'}\\
 c a'\sqrt{s r'}+d c'\sqrt{r s'}
&c b'\sqrt{s s'}+d d'\sqrt{r r'}
\end{bmatrix}
\end{align}
for two elements of $\widehat{G}$. 
Since $\sqrt{s r}=\sqrt{s' r'}=\sqrt{n}$, 
$\widehat{G}$ is closed under the multiplication. 
Obviously the inverse is also well-defined on $\widehat{G}$. 
Hence $\widehat{G}$ is a subgroup of $\SL(2,\bbR)$. 

(2) is obvious and (3) is a consequence of the following calculation.
\begin{align*}
&\begin{bmatrix}a\sqrt{r}&c\sqrt{s}\\b\sqrt{s}&d\sqrt{r}\end{bmatrix}
\begin{bmatrix}x&y\sqrt{n}\\y\sqrt{n}&z\end{bmatrix}
\begin{bmatrix}a\sqrt{r}&b\sqrt{s}\\c\sqrt{s}&d\sqrt{r}\end{bmatrix}\\
=&\begin{bmatrix}a^2r x+2a c n y+c^2s z&(a b x+(a d r+b c s)y+c d z)\sqrt{n}\\
(a b x+(a d r+b c s)y+c d z)\sqrt{n}&b^2s x+2b d n y+d^2r z\end{bmatrix}.
\end{align*}
\end{proof}

\begin{rmk}
If $(H^2)=2$, that is, if the surface is principally polarized, 
then $G=\SL(2,\bbZ)$. 
This group is mentioned in \cite[Example 1.15 (1)]{Mukai:1998}. 
\end{rmk}

For stating Lemma \ref{lem:SO_FMT}, we introduce a terminology in $\Sym_2(\bbR)$, corresponding to the notion of the positive Mukai vector.

\begin{dfn}
An element $\begin{bmatrix}x&y\\y&z\end{bmatrix}\in\Sym_2(\bbR)$ 
is called \emph{positive} if $x>0$, or if $x=0$ and $y>0$, 
or if $x=y=0$ and $z>0$.
\end{dfn}

\begin{lem}
\label{lem:SO_FMT}
Consider the action \eqref{eq:G_action} of $\PSL(2,\bbR)$ on $\Sym_2(\bbZ,n)$.
\begin{enumerate}
\item
If $\gamma_1,\gamma_2\in\PSL(2,\bbR)$ have the same action on $\begin{bmatrix}1&0 \\ 0&0\end{bmatrix}$ and $\begin{bmatrix}0&0 \\ 0&1\end{bmatrix}$, then $\gamma_1=\gamma_2$.
\item
Assume that we are given two positive elements $X_1,X_2\in\Sym_2(\bbZ,n)$ with $B(X_1,X_1)=B(X_2,X_2)=0$ and $B(X_1,X_2)=-1$. Then there exists a unique element $\gamma\in G/\{\pm1\}\subset\PSL(2,\bbR)$ such that $\begin{bmatrix}1&0 \\ 0&0\end{bmatrix}\cdot \gamma=X_1$ and $\begin{bmatrix}0&0 \\ 0&1\end{bmatrix}\cdot \gamma=X_2$.
\end{enumerate}
\end{lem}

\begin{proof}
(1)
Take elements $g_1\seteq\begin{bmatrix}x&y\\z&w\end{bmatrix},g_2\seteq\begin{bmatrix}p&q\\r&s\end{bmatrix}\in\SL(2,\bbR)$ such that their equivalent classes in $\PSL(2,\bbR)$ equal to $\gamma_1,\gamma_2$. Then the condition is rewritten into 
\begin{align*}
 {}^t g_1\begin{bmatrix}1&0 \\ 0&0\end{bmatrix}g_1
={}^t g_2\begin{bmatrix}1&0 \\ 0&0\end{bmatrix}g_2,\quad
 {}^t g_1\begin{bmatrix}0&0 \\ 0&1\end{bmatrix}g_1
={}^t g_2\begin{bmatrix}0&0 \\ 0&1\end{bmatrix}g_2.
\end{align*}
 These equations yield $x^2=p^2$, $x y=p q$, $y^2=q^2$, $z^2=r^2$, $z w=r s$, $w^2=s^2$. Hence $(x,y)=\ep_1(p,q)$ and $(z,w)=\ep_2(r,s)$ where $\ep_1=1$ or $-1$ and $\ep_2=1$ or $-1$. Since $g_1,g_2\in\SL(2,\bbR)$, we should have $\ep_1=\ep_2$. Therefore $g_1=\pm g_2$, which means $\gamma_1=\gamma_2$.

(2)
The uniqueness follows from (1). We only have to show the existence.

Set $X_i=\begin{bmatrix}r_i&d_i\sqrt{n}\\d_i\sqrt{n}&a_i\end{bmatrix}$ ($i=1,2$). Then the assumption is equivalent to 
\begin{align*}
d_1^2n-r_1a_1=0,\quad d_2^2n-r_2a_2=0,\quad 2d_1d_2n-r_1a_2-r_2a_1=-1.
\end{align*}
By the first two equations we can take $x,y,z,w\in\bbR$ such that $x^2=r_1$, $x y/\sqrt{n}=d_1$, $y^2=a_1$, $z^2=r_2$, $z w/\sqrt{n}=d_2,w^2=a_2$. By the third equation we also have $x w-y z=\pm1$. Hence we may assume $x w-y z=+1$ by changing the signs of $z$ and $w$. Now set $g\seteq\begin{bmatrix}x&y\\z&w\end{bmatrix}$. Then $g \in G\subset\SL(2,\bbR)$. Then the equivalent class $\gamma\seteq\pm g\in\PSL(2,\bbZ)$ can be defined. This $\gamma$ satisfies the requirement by construction.
\end{proof}

\subsection{Matrix description}\label{subsec:matrix}

We introduce a description of Mukai vectors via $2\times 2$ matrix form, 
which is a generalization of the one defined in \cite{Mukai:1980}. 
We suppose $\NS(X)=\bbZ H$ and set $n\seteq(H^2)/2\in\bbZ_{>0}$.

First we define an isomorphism of modules
\begin{align*}
\begin{array}{l l l l}
\iota_X \colon
  &H^{\ev}(X,\bbZ)_\alg &\simto     &\Sym_2(\bbZ,n)\\
  &(r,d H,a)            &\mapsto        
  &\begin{bmatrix}r&d\sqrt{n}\\d\sqrt{n}&a\end{bmatrix}.
\end{array}
\end{align*}
Then it is clear that 
\begin{align*}
\mpr{v,v'}=B(\iota_X(v),\iota_X(v')).
\end{align*}
Thus we have an isometry of lattices:
\begin{align*}
\iota_X\colon (H^{\ev}(X,\bbZ)_\alg,\mpr{\cdot,\cdot})
\simto(\Sym_2(\bbZ,n),B).
\end{align*}

We also have the next restatement of Lemma \ref{lem:group} and Lemma \ref{lem:SO_FMT}.

\begin{lem}\label{lem:pmvpr}
If $X$ is an abelian surface with $\NS(X)=\bbZ H$ and $(H^2)=2n$, then any pair of positive primitive isotropic Mukai vectors $v_1,v_2\in H^{\ev}(X,\bbZ)_\alg$ satisfying $\mpr{v_1,v_2}=-1$ can be uniquely written as
\begin{align*}
&v_1=(a^2r,a b H,b^2s),\quad v_2=(c^2s,c d H,d^2r),\\
&rs=n, \quad a d r-b c s=1,\quad a>0.
\end{align*}
\end{lem}

\begin{proof}
This can be shown in the same manner as  Lemma \ref{lem:SO_FMT}(2), so we omit the detail.
\end{proof}

\begin{dfn}\label{dfn:Eq0&theta}
\begin{enumerate}
\item
For $Y,Z\in\FM(X)$, we set 
\begin{align*}
\Eq_0({\bf D}(Y),{\bf D}(Z))
&\seteq
 \{\Phi_{Y \to Z}^{{\bfE}[2k]} 
   \in \Eq({\bf D}(Y),{\bf D}(Z))
   \mid {\bfE} \in \Coh(Y \times Z), k \in {\bbZ} 
 \},
\\
\calE(Z)
&\seteq
 \bigcup_{Y\in\FM(X)}\Eq_0({\bf D}(Y),{\bf D}(Z)),
\\
\calE
&\seteq
 \bigcup_{Y,Z\in\FM(X)}\Eq_0({\bf D}(Y),{\bf D}(Z)).
\end{align*}
\item
Let $Y,Z\in\FM(X)$.
For a FMT $\Phi\seteq\Phi^{\bfE^\bullet}_{Y\to Z}\colon\bfD(Y)\to\bfD(Z)$ 
with $\bfE^\bullet\in\bfD(Y\times Z)$, set 
\begin{align*}
\theta(\Phi)\seteq \iota_Z\circ \Phi^H \circ\iota_Y^{-1},
\end{align*}
where $\Phi^H\colon H^\ev(Y,\bbZ)\to H^\ev(Z,\bbZ)$ is the cohomological FMT 
induced by $\Phi$ (see \S\S\,\ref{subsec:ms}). 
\end{enumerate}
\end{dfn}

Note that $\calE=\bigcup_{Z\in\FM(X)}\calE(Z)$ and 
that $\calE$ is a groupoid with respect to the composition of the equivalences.
Since $\iota_Y,\iota_Z$ and $\Phi^H$ are isometries, 
their composition $\theta(\Phi)$ is again an isometry 
on the lattice $(\Sym_2(\bbZ,n),B)$. 
Thus we have a map $\theta\colon\calE\to\LieO(B)$.

\begin{thm}\label{thm:bij}
For each $Z\in\FM(X)$, we have a bijection
\begin{align*}
\theta(\calE(Z)) \cong G/\{\pm 1\}.
\end{align*}
\end{thm}

\begin{proof}
It is enough to show for the case $Z=X$.
First we use Lemma~\ref{lem:SO_FMT} to make a map
\begin{align*}
\kappa\colon\theta(\calE(X))\longrightarrow G/\{\pm 1\}
\end{align*}
Take an element $\Phi=\Phi_{Y\to X}^{\bfE[2k]}\in\Eq(\bfD(Y),\bfD(X))$. 
By Fact~\ref{fact:FMT-abelian-surface}~(3), 
there are primitive isotropic Mukai vectors 
$v_i=(r_i,d_i H,a_i)\in H^{\ev}(X,\bbZ)_\alg$ for $i=1,2$ 
such that $v(\bfE_y)=v_2$, $Y\cong M_X^H(v_2)$, $\mpr{v_1,v_2}=-1$ and 
$\bfE$ is a universal family on $Y\times X$. 
Set $X_i\seteq \iota_X(v_i)\in\Sym_2(\bbZ,n)$. 
Then the assumption of Lemma~\ref{lem:SO_FMT}~(2) holds. 
Thus a map $\kappa$ is defined. 
The same Lemma~\ref{lem:SO_FMT} implies that $\kappa$ is injective.

For the converse relation, 
let us take $\gamma\in G/\{\pm 1\}$ and 
choose a representative $g\seteq\begin{bmatrix}x&y\\z&w\end{bmatrix}\in G$ 
of $\gamma$. 
Then  by Lemma~\ref{lem:group} there exists integers $a,b,c,d,r,s$ 
uniquely determined by $g$. 
If we set $v_1\seteq(a^2 r,a b H, b^2 s)$ and $v_2\seteq(c^2s,c d H,d^2r)$, 
then $\mpr{v_1^2}=\mpr{v_2^2}=0$ and $\mpr{v_1,v_2}=-(a d r-b c s)^2=-1$. 
Thus we can apply Proposition~\ref{prop:semihom-fmt}~(4) 
to the pair $v_1,v_2$. 
Hence there exists a FMT $\Psi\seteq\Phi_{X\to Y}^{\bfE^\vee[2]}$ 
such that $\Psi(E)=\bbC_y$ and $\Psi(F)=\calO_Y[-i]$ ($i=0$ or $2$) 
for $E\in M_X^H(v_2)$ and $F\in M_X^H(v_1)$. 
Then the inverse transform $\Phi\seteq\Phi_{Y\to X}^{\bfE}$ 
satisfies $\Phi^H((1,0,0))=v_1$ and $\Phi^H((0,0,1))=v_2$. 
Set $X_i\seteq \iota_X(v_i)$. Now using $\ref{eq:matrix}$, 
we find that $\theta(\Phi)$ and 
$\gamma$ have the same action on $X_1$ and $X_2$. 
Since $X_1$ and $X_2$ satisfy the assumption on Lemma~\ref{lem:SO_FMT}~(1), 
we find that $\kappa\circ\theta(\Phi_{Y\to X}^{\bfE})$ coincides with $\gamma$.
\end{proof}

\begin{rmk}
Lemma \ref{lem:det:positive} means 
that the image $\beta\circ\theta(\Phi)\in\LieO(2,1)$ 
is actually in $\SO(2,1)$. 
Theorem~\ref{thm:bij} further claims 
that $\beta\circ\theta(\Phi)\in\SO_0(2,1)$. 
We can prove it by using Lemma~\ref{lem:Opq} as follows. 
The matrix form of $\theta(\Phi)$ in the basis 
$\langle 
\begin{bmatrix}1&0\\ 0& 0\end{bmatrix},
\begin{bmatrix}0&1\\ 1& 0\end{bmatrix},
\begin{bmatrix}0&0\\ 0& 1\end{bmatrix}
\rangle$ 
of $\Sym_2(\bbR)$ is the $P_2 M P_1$ in \eqref{eq:g33} 
(see Remark \ref{rmk:g33}). 
The element $\beta\circ\theta(\Phi)$ is a matrix presented in the basis
$\langle 
\begin{bmatrix}1&0\\ 0&-1\end{bmatrix},
\begin{bmatrix}0&1\\ 1& 0\end{bmatrix},
\begin{bmatrix}0&1\\-1& 0\end{bmatrix}
\rangle$ of $\Liesl(2,\bbR)$. 
The transition matrix between these two bases is the $P$ in \eqref{eq:P}. 
Thus we should calculate $\beta\circ\theta(\Phi)=Q\seteq P^{-1} P_2 M P_1 P$. 
A straightforward computation yields 
\begin{align*}
Q_{33}
&=
\dfrac{1}{2 r_0^2}
(r_0+r_0^3
 +a_0 r_0^2+\widetilde{a}_0 r_0^2+a_0 \widetilde{a}_0 r_0
 -2 n d_0 \widetilde{d}_0)\\
&=\dfrac{1}{2 r_0^2}
(n(d_0-\widetilde{d}_0)^2+n(r_0-1)d_0^2+n(r_0-1)\widetilde{d}_0^2
+(a_0 \widetilde{a}_0+1)r_0)
\end{align*}
Since all the $a_0,\widetilde{a}_0,r_0$ are positive, we find that $Q_{33}>0$.
 Since we have known $\beta\circ\theta(\Phi)\in\SO(2,1)$, 
Lemma~\ref{lem:Opq} shows that $\beta\circ\theta(\Phi)\in\SO_0(2,1)$
\end{rmk}

Next we treat the composition of FMT and the dualizing functor.
For a FMT $\Phi_{X \to Y}^{\bfE^\bullet} \in \Eq(\bfD(X),\bfD(Y))$,
we set
\begin{align*}
\theta(\Phi_{X \to Y}^{\bfE}\calD_X)
\seteq
&
\iota_Y \circ (\Phi_{X \to Y}^{\bfE})^H \calD_X \circ \iota_X^{-1}
\colon
\Sym_2(\bbR)\to\Sym_2(\bbR).
\end{align*}
Here the $\calD_X$ in the last line means the cohomological action 
$H^{\ev}(X,\bbZ)_{\alg}\to H^{\ev}(X,\bbZ)_{\alg}$ of the dualizing functor.

\begin{lem}
\label{lem:matrix}
\begin{enumerate}
\item
\begin{align*}
\theta(\calD_X)=
\pm\begin{bmatrix} 1 & 0 \\ 0 & -1 \end{bmatrix}
\in \PSL(2,\bbR).
\end{align*}
We also have $\theta(\Phi_{X \to Y}^{\bfE}\calD_X)\in\PSL(2,\bbR)$ 
for a FMT $\Phi_{X \to Y}^{\bfE^\bullet} \in \Eq(\bfD(X),\bfD(Y))$.

\item
If $\theta(\Phi_{X \to Y}^{\bfE})=
\begin{bmatrix} a & b\\ c & d \end{bmatrix}$,
then 
\begin{align*}
\theta(\Phi_{Y \to X}^{\bfE})=
\pm\begin{bmatrix}d & b\\c & a \end{bmatrix},\quad
\theta(\Phi_{Y \to X}^{\bfE^{\vee}[2]})=
\pm\begin{bmatrix}d & -b\\-c & a \end{bmatrix},\quad
\theta(\Phi_{X \to Y}^{\bfE^{\vee}[2]})=
\pm\begin{bmatrix}a & -b\\-c & d \end{bmatrix}.
\end{align*}
\end{enumerate}
\end{lem}

\begin{proof}
(1)
This is trivial.

(2) The claims follows from the following relations.
\begin{align*}
\Phi_{Y \to X}^{\bfE^{\vee}[2]}
=&(\Phi_{X \to Y}^{\bfE})^{-1},
\\
\Phi_{Y \to X}^{\bfE}
=& \calD_X  \Phi_{Y \to X}^{\bfE^{\vee}[2]}\calD_Y
=  \calD_X (\Phi_{X \to Y}^{\bfE})^{-1}    \calD_Y,
\\
\Phi_{X \to Y}^{{\bf E}^{\vee}[2]}
=& (\Phi_{Y \to X}^{\bfE})^{-1}
=\calD_Y \Phi_{X \to Y}^{\bfE} \calD_X.
\end{align*}
The first equality is obvious and 
the second one is a consequence of the Serre duality (Lemma~\ref{lem:Serre}).
For the third equality, we take the inverse of the second equality.
\end{proof}

Now we summarize the argument of this section in the following manner.

\begin{prop}\label{prop:matrix}
Suppose $\NS(X)=\bbZ H$ and set $n\seteq(H^2)/2\in\bbZ_{>0}$. 
Consider a pair of primitive isotropic Mukai vectors $v_1$ and $v_2$ on $X$ 
such that $\mpr{v_1,v_2}=-1$. 
By Lemma \ref{lem:pmvpr}, we may put
\begin{align*}
&v_1=(p^2_1r_1,p_1q_1 H,q^2_1r_2),\quad 
 v_2=(p^2_2r_2,p_2q_2 H,q^2_2r_1),\\
&r_1r_2=n, \quad p_1q_2r_1-p_2q_1r_2=1,\quad p_1>0.
\end{align*}
Let us denote by $E_i$ a simple semi-homogeneous sheaf 
with Mukai vector $v_i$, where $i=1,2$. 
Proposition \ref{prop:semihom-fmt} shows 
that there exists a FMT $\Phi\seteq \Phi_{Y\to X}^{\bfE}$ 
which converts the structure sheaf $\calO_Y$ 
into $E_1$ and a skyscraper sheaf into $E_2$ up to even shift. 

Then the action of the cohomological FMT 
$\Phi^H : H^{\ev}(Y,\bbZ)_{\alg}\to H^{\ev}(X,\bbZ)_{\alg}$ 
can be written as follows. 
Let $\widehat{H}$ be the ample generator of $\NS(Y)$. 
For a Mukai vector $v=(r,d\widehat{H},a)\in H^{\ev}(Y,\bbZ)_\alg$, 
the image $\Phi^H(v)$ is given by
\begin{align*}
\Phi^H(v)=v\cdot g,\quad 
g\seteq\begin{bmatrix}
p_1\sqrt{r_1}&q_1\sqrt{r_2}\\p_2\sqrt{r_2}&q_2\sqrt{r_1}
\end{bmatrix}\in G.
\end{align*}
Here the action $\cdot$ is given by
\begin{align}
\label{eq:action_cdot}
 (r,d H,a)\cdot g\seteq(r',d' H,a'),\quad
 {}^t g
 \begin{bmatrix}r &d  \sqrt{n}\\d  \sqrt{n}&a \end{bmatrix} 
 g
=\begin{bmatrix}r'&d' \sqrt{n}\\d' \sqrt{n}&a'\end{bmatrix}.
\end{align}

Conversely, given $\gamma \in G/\{\pm 1\}$, 
there exists a FMT $\Phi=\Phi_{Y\to X}^\bfE$ 
such that $\theta(\Phi)=\gamma$, 
and this FMT is unique up to even shifts, 
automorphism of $X$ and tensor with a degree zero line bundle.
\end{prop}

\section{Tame complex}
\label{section:complex}

This section gives some examples of calculations involving the semi-homogeneous presentation and the matrix description.  We consider a certain sequence of FMTs associated to the power of a quadratic matrix. The resulting images of the original sheaf defines a sequence of complexes, which we call \emph{tame complexes}. These complexes are easy to handle in the sense that we can state explicit numerical criteria for a FMT to transform the complex into a stable sheaf. The criteria, which are the main results of this section, are given in Theorems \ref{thm:csv:1} and \ref{thm:csv:2}.

We fix an abelian surface $X$ with $\NS(X)=\bbZ H$. 
We define $n\seteq(H^2)/2\in\bbZ_{>0}$.

\subsection{The arithmetic group $G$ and numerical solutions for the ideal sheaf}\label{subsect:solution}

Before starting the calculation, 
we notice the next result on the number of numerical solutions 
for \eqref{eq:numerical_equation}.

\begin{prop}\label{prop:num_sol}
Let $E$ be a stable sheaf on $X$ with $\ell\seteq\mpr{v(E)^2}/2>0$.
Then the numerical equation \eqref{eq:numerical_equation} 
for the Mukai vector $v\seteq v(E)$ gives the next criterion.
\\
\qquad
If $n\ell$ is a square number, 
then the number of numerical solution is 0 or 1.
\\
\qquad
Otherwise, the number of numerical solutions is 0 or infinity.
\end{prop}
\begin{proof}
We may assume that there is at least one solution for $v$.
Then there exists a FMT $\Phi$ whose cohomological counterpart $\Phi^H$ 
maps $v$ to $(1,0,-\ell)$, 
and it is sufficient to prove the claim for $v=(1,0,-\ell)$.

We rewrite the numerical equation as $v=\pm (\ell u_1-u_2)$, 
where $u_1,u_2$ are positive isotropic Mukai vectors with $\mpr{u_1,u_2}=-1$.
Using Lemma~\ref{lem:pmvpr}, 
we can further rewrite the numerical equation as
\begin{align*}
&(1,0,-\ell)=\pm(\ell u_1-u_2),\quad 
u_1=(p_1^2r_1,p_1q_1H,q_1^2a_1),\quad\mpr{u_2^2}=0,
\end{align*}
where $p_1,q_1\in\bbZ$, $r_1,a_1\in\bbZ_{>0}$ and $r_1a_1=n$.
Then we have
\begin{align*}
0&=\mpr{u_2^2}=\mpr{(\ell p_1^2r_1\mp1,\ell p_1q_1H,\ell(q_1^2a_1\pm1)}^2\\
&=2n(\ell p_1q_1)^2-2(\ell p_1^2r_1\mp1)\ell(q_1^2a_1\pm1)\\
&=2n\ell^2p_1^2q_1^2
-2\ell(n\ell p_1^2q_1^2\mp q_1^2a_1\pm \ell p_1^2r_1-1)\\
&=2\ell(\pm q_1^2a_1\mp \ell p_1^2r_1+1).
\end{align*}
Therefore 
the above equation becomes
\begin{align*}
\ell p_1^2r_1-q_1^2a_1=\pm1.
\end{align*}

Let us consider the solution of this equation. 
Suppose that $n\ell$ is a square number. 
We may suppose $\gcd(\ell r_1,a_1)=1$. 
Then since $\ell r_1a_1=\ell n$ is a square number, 
both $\ell r_1$ and $a_1$ are square numbers. 
Hence the equation above becomes 
$(\sqrt{\ell r_1}p_1-\sqrt{a_1}q_1)(\sqrt{\ell r_1}p_1+\sqrt{a_1}q_1)=\pm1$. 
Then the solution is as follows.
\begin{align*}
\begin{array}{ll}
\ell u_1-u_2=(-1,0,\ell) \colon &u_1=(0,0,1),u_2=(1,0,0),\\
\ell u_1-u_2=(1,0,-\ell) \colon &\text{a solution exists only when }\ell=1\ 
\mbox{and}\\
&u_1=(1,0,0),u_2=(0,0,1).
\end{array}
\end{align*}

Next suppose that $n\ell$ is not a square number. 
Put $r_1=n$ and $a_1=1$. 
Then by the above calculation, 
the equation $\ell u_1-u_2=(-1,0,\ell)$ leads to $q_1^2-\ell n p_1^2=1$, 
which has infinite number of solutions.
\end{proof}

\begin{eg}
If $\rk(v)=0$, then we can write $v=(0,d H,a)$ with $\ell n=(n d)^2$.
In this case, there is no kernel presentation.
\end{eg}

Hereafter we assume that $\sqrt{\ell n}\notin\bbZ$. 
By Theorem \ref{thm:presentation}, 
a general member of $M_X^H(v)$ has both kernel and cokernel presentations 
if a numerical solution exists.

Our next task is to describe the numerical solution of the ideal sheaf of 0-dimensional subscheme. First we introduce an arithmetic group $S_{n,\ell}$.

\begin{dfn}
For $(x,y)\in\bbR^2$, set
\begin{align*}
P(x,y)\seteq\begin{bmatrix}y&\ell x\\ x&y\end{bmatrix}.
\end{align*}
We also set
\begin{align*}
S_{n,\ell}\seteq
\left\{\begin{bmatrix}y&\ell x\\ x&y\end{bmatrix}
\,\bigg|\,
x^2,y^2,\dfrac{x y}{\sqrt{n}}\in\bbZ,\ 
y^2-\ell x^2=\pm 1 \right\}.
\end{align*}
\end{dfn}

\begin{lem}
\begin{enumerate}
\item
$S_{n,\ell}$ is a commutative subgroup of $\GL(2,\bbR)$.
\item
We have a homomorphism
\begin{align*}
\begin{array}{c c c c}
\phi \colon & S_{n,\ell} & \longrightarrow & \bbR^\times \\
            & P(x,y)     & \mapsto         & y+x\sqrt{\ell}.
\end{array}
\end{align*}
\item
For $\ell>1$ $\phi$ is injective,  and for $\ell=1$ we have
\begin{align*}
\Ker\phi=\left\langle\begin{bmatrix} 0&1 \\ 1&0 \end{bmatrix}\right\rangle
\end{align*}
\item
We set a subgroup $G_{n,\ell}$ of $\widehat{G}$ (Definition~\ref{dfn:G}) to be 
\begin{align*}
G_{n,\ell}\seteq 
\left\{
 g \in \widehat{G} 
 \left|\; 
 {}^{t} g \begin{bmatrix} 1 & 0 \\ 0 & -\ell \end{bmatrix}  g
 =\pm\begin{bmatrix} 1 & 0 \\ 0 & -\ell \end{bmatrix}
 \right.
\right\}.
\end{align*}
Then
\begin{align*}
G_{n,\ell}
=S_{n,\ell} \rtimes 
\left \langle
 \begin{bmatrix} 1 & 0 \\ 0 & -1 \end{bmatrix}
\right\rangle.
\end{align*}
\end{enumerate}
\end{lem}

\begin{proof}
The proofs of (1) and (2) are straightforward.

For (3), assume that 
$x,y\in\bbR$ with $x^2,y^2,x y/\sqrt{n}\in\bbQ$ satisfy $y+x\sqrt{\ell}=1$. 
Then $(y^2+\ell x^2)+2(x y/\sqrt{n})\sqrt{\ell n}=(y+x\sqrt{\ell})^2=1$. 
Our assumptions yields $y^2+\ell x^2=1$ and $x y=0$. If $x=0$, then $y=\pm 1$. 
If $y=0$, then $\ell=1$ and $x=1$. 
Hence the conclusion holds.

(4) follows from direct computations.
\end{proof}

Then the Dirichlet unit theorem yields the following corollary.

\begin{cor}
If $\ell>1$, then $S_{n,\ell}\cong\bbZ\oplus\bbZ/2\bbZ$.
\end{cor}
\begin{proof}
Let $p_1,p_2,\ldots,p_m$ be the prime divisors of $\ell n$ and $\frako$ the ring of algebraic integers in $\bbQ(\sqrt{p_1},\sqrt{p_2},\ldots,\sqrt{p_m})$. By Dirichlet unit theorem $\frako^\times$ is a finitely generated abelian group whose torsion subgroup is $\{\pm 1\}$. Hence $\phi(S_{n,\ell})$ is a finitely generated abelian group whose torsion subgroup is $\{\pm 1\}$. For $A\in S_{n,\ell}$, we have $\phi(A^2)\in\bbZ[\sqrt{\ell n}]$. Since $\bbZ[\sqrt{n\ell}]^\times\cong\bbZ\oplus\bbZ/2\bbZ$, we get $S_{n,\ell}\cong\bbZ\oplus\bbZ/2\bbZ$.
\end{proof}

\begin{rmk}
If $n=1$ and $\ell>1$, 
then $S_{1,\ell}$ is the group of units of $\bbZ[\sqrt{\ell}]$. 
Moreover if $\ell$ is square free and $\ell\equiv 2,3 \pmod{4}$, 
then since $\bbZ[\sqrt{\ell}]$ is 
the ring of the integers of $\bbQ[\sqrt{\ell}]$, 
a generator of $S_{1,\ell}$ becomes a fundamental unit.
\end{rmk}

Now we will construct families of simple complexes for all numerical 
solutions of $(1,0,-\ell)$.

\begin{lem}\label{lem:numerical-P(p,q)}
For two positive isotropic Mukai vectors $w_0,w_1$ 
on the fixed abelian surface $X$, 
the condition
\begin{align*}
(1,0,-\ell)=\pm (\ell w_0-w_1),\quad \mpr{w_0,w_1}=-1
\end{align*}
is equivalent to
\begin{align*}
w_0=(p^2,-\tfrac{p q}{\sqrt{n}}H,q^2),\ 
w_1=(q^2,-\tfrac{\ell p q}{\sqrt{n}}H,\ell^2 p^2),\quad
P(p,q)\in S_{n,\ell}.
\end{align*}
\end{lem}

\begin{proof}
If there are isotropic Mukai vectors with the first condition, 
then we can write them as $w_0=(r,d H,(r\ell\pm 1))$ and 
$w_1=(r \ell\mp 1,d\ell H, r \ell^2)$, 
where $d^2(H^2)=2 r (r \ell\mp 1)$. 
We set $p\seteq \sqrt{r}$ and $q\seteq \sqrt{r \ell \mp 1}$.
Then $w_0=(p^2,-\tfrac{p q}{\sqrt{n}}H,q^2)$,
$w_1=(q^2,-\tfrac{\ell p q}{\sqrt{n}}H,\ell^2 p^2)$
and $\ell p^2-q^2=\pm 1$.
Thus $P(p,q)\in S_{n,\ell}$. 
The converse is obvious.
\end{proof}

\begin{cor}
Recall the action $\cdot$ of $\GL(2,\bbR)$ given in \eqref{eq:action_cdot}.
By the correspondence 
\begin{align*}
g \in G_{n,\ell}\mapsto (w_0,w_1),\;
w_0=(0,0,1)\cdot g, \, w_1=(1,0,0)\cdot g, 
\end{align*}
we have a bijective correspondence:
\begin{align*}
\begin{array}{c c c}
G_{n,\ell}
\left/
 \left\langle \pm \begin{bmatrix} 1 & 0\\ 0 & -1 \end{bmatrix} \right \rangle 
\right.
 \cong 
S_{n,\ell}/\{\pm 1\} 
& 
\longleftrightarrow 
&
\left\{ (w_0,w_1) 
 \left|
 \begin{aligned}
 \mpr{w_0,w_1}=-1,\, \mpr{ w_0^2 }=\mpr{w_1^2}=0,\\
 w_0,w_1>0,\, (1,0,-\ell)=\pm(\ell w_0-w_1)
 \end{aligned}
 \right. 
\right\}
\end{array}
\end{align*}
\end{cor}

\begin{dfn}
\label{dfn:ambm}
Assume that $\ell>1$.
Let $z\seteq q+\sqrt{\ell}p$ ($p,q>0$) 
be the generator of $S_{n,\ell}/\{\pm 1\}$.
We set $\epsilon\seteq q^2-\ell p^2\in\{\pm 1\}$.
For $m \in \bbZ$, we set
\begin{align*}
\begin{bmatrix}   q & \ell p  \\ p   & q   \end{bmatrix}^m= 
\begin{bmatrix} b_m & \ell a_m\\ a_m & b_m \end{bmatrix}.
\end{align*}
\end{dfn}

By the definition we have 
\begin{align*}
(a_0,b_0)=(0,1),\quad
(a_{-m},b_{-m})=\ep^m(-a_m,b_m),\, m\in\bbZ_{>0}
\end{align*}
and
\begin{align}
\label{eq:S_{n,l}}
S_{n,\ell}= 
\left\{ \left. \pm
 \begin{bmatrix} b_m & \ell a_m\\ a_m & b_m \end{bmatrix}
 \right| m \in \bbZ \right\}.
\end{align}

Next we consider the right action of $\GL(2,\bbR)$ on $\bbR^2$
\begin{align}
\label{eq:quad_action}
(x,y) \mapsto (x,y) X,
\quad X \in \GL(2,\bbR).
\end{align}
Then the quadratic map 
\begin{align*}
\begin{array}{c c c}
 \bbR^2 & \to     & \Sym_2(\bbR)\\[3pt]
 (x,y)  & \mapsto &
\begin{bmatrix} x \\ y \end{bmatrix}
\begin{bmatrix} x & y  \end{bmatrix}
\end{array}
\end{align*}
is $\GL(2,\bbR)$-equivariant.   

Using this action, we have the next descriptions of
the topological invariants of fine moduli spaces
$M_X^H(v)$ of dimension 2.
\begin{align}
\label{eq:top-inv}
\begin{split}
& 
\left\{ v 
 \left|
  \begin{aligned}
   v \in H^{\ev}(X,\bbZ)_{\alg},\, \mpr{v^2}=0,\, v>0,
   \\
   \mpr{w,v}=-1,\,\exists w \in H^{\ev}(X,\bbZ)_{\alg}
  \end{aligned}
 \right. 
\right\} 
\\ 
\overset{\varphi_1}{\longleftrightarrow} 
&
\left\{ 
 (a\sqrt{r}, b\sqrt{s}) \in \bbR^2
 \left|
  \begin{aligned}
  a,b\in\bbZ,\, r,s\in\bbZ_{>0},
  \\ 
  r s=n,\, \gcd(a r,b s)=1 
  \end{aligned}
 \right.
\right\}
/\{\pm 1\}
=
\left\{  (0,1) X \mid  X \in G \right\}/\{\pm 1\}
\\
\overset{\varphi_2}{\longleftrightarrow} 
&
\left\{
 \left. \dfrac{b \sqrt{n}}{a r} \in \bbP^1(\bbR) =\bbR \cup \{ \infty \} 
 \right| 
 r s=n,\; \gcd(a r,b s)=1
\right\},
\end{split}
\end{align}
where we used the correspondences 
\begin{align*}
v=(a^2 r, a b H,b^2 s) 
\overset{\varphi_1}{\longleftrightarrow}
 \pm(a\sqrt{r},b\sqrt{s})
\overset{\varphi_2}{\longleftrightarrow}
 \dfrac{\mu(v)}{\sqrt{n}}=\dfrac{b \sqrt{n}}{a r}.
\end{align*}
These correspondences are $\widehat{G}$-equivariant 
under the action \eqref{eq:quad_action}.

Lemma \ref{lem:numerical-P(p,q)}, \eqref{eq:S_{n,l}} and \eqref{eq:top-inv}
imply the following one to one correspondence:
\begin{align*}
\begin{array}{c c c}
\left\{ 
 \{ v_1,v_2 \} \left| 
 \begin{aligned}
  &\text{ There is a numerical solution}
  \\
  &\text{ $(v_1,v_2,\ell_1,\ell_2)$ of $(1,0,-\ell)$ }
 \end{aligned}
 \right. 
\right\} 
& \longleftrightarrow & 
\left\{ 
 \left. 
  \left\{\frac{b_m}{a_m}, \frac{\ell a_m}{b_m} \right\} \subset \bbP^1(\bbR) 
 \right| 
 m \in \bbZ 
\right\}
\\
\{ v_1,v_2 \} 
& \longleftrightarrow & 
\left\{\frac{\mu(v_1)}{\sqrt{n}},\frac{\mu(v_2)}{\sqrt{n}} \right\},
\end{array}
\end{align*}
where $(\ell_i,\ell_j)=(\ell,1)$ if and only if 
$(\frac{\mu(v_i)}{\sqrt{n}},\frac{\mu(v_j)}{\sqrt{n}})=
(\frac{b_m}{a_m},\frac{\ell a_m}{b_m})$.

\begin{dfn}
For the numerical solution $(v_1,v_2,\ell_1,\ell_2)$ 
corresponding to $\{\frac{b_{-m}}{a_{-m}},\frac{\ell a_{-m}}{b_{-m}} \}
=\{-\frac{b_m}{a_m},-\frac{\ell a_m}{b_m} \}$,
$\frakM^{\pm}_m$ denotes the moduli space
$\frakM^{\pm}(v_1,v_2,\ell_1,\ell_2)$. 
$\frakM^{\pm}_m[-1]$ denotes the moduli space
of $V^{\bullet}[-1]$, $V^{\bullet} \in \frakM^{\pm}_m$. 
\end{dfn}

In this section, we will explicitly construct families 
$\{F_m^\bullet\}_{m\in\bbZ}$
of complexes associated to all 
$\{-\frac{b_m}{a_m},-\frac{\ell a_m}{b_m} \}$, $m \in \bbZ$.

\begin{prop}\label{prop:explicit}
We have families of complexes $\{F_m^\bullet\}_{m\in\bbZ}$ such that
\begin{align*}
\begin{cases} 	
F_m^\bullet   \in \frakM^+_{m+1} \cap \frakM^-_{m} \; &  m>0,\\
F_0^\bullet\, \in \frakM^+_1 \cap \frakM^+_0[-1],\\
F_m^\bullet   \in \frakM^+_{m}[-1] \cap \frakM^-_{m+1}[-1] \; &  m<0,
\end{cases}
\end{align*}
and a sequence of isomorphisms:
\begin{align*}
\begin{array}{c c c c c c c}
  \cdots \overset{\Psi}{\to}
& \frakM^-_0[-1] \cap \frakM^+_{-1}[-1] 
& \overset{\Psi}{\to} 
& \frakM^+_1 \cap \frakM^+_0[-1]
& \overset{\Psi}{\to} \cdots 
\\
 \cdots \mapsto 
& F_{-1}^\bullet
& \mapsto 
& F_{0}^\bullet
& \mapsto\cdots
\\
\\
&
& \cdots \overset{\Psi}{\to} 
& \frakM^{+}_m \cap \frakM^{-}_{m-1} 
& \overset{\Psi}{\to}
& \frakM^{+}_{m+1} \cap  \frakM^{-}_{m} 
& \overset{\Psi}{\to}
\cdots
\\
&
& \cdots\mapsto 
& F_{m-1}^\bullet
& \mapsto 
& F_m^\bullet
& \mapsto\cdots
\end{array}
\end{align*}
where $\Psi$ is the isomorphism in Proposition~\ref{prop:dual}.
\end{prop}

We will call the complex $F_m^\bullet$ the \emph{tame complex}. 
Its construction and the proof of Proposition~\ref{prop:explicit}
will be given in the following subsections.

\subsection{Construction of the tame complex: case I}

Under the notation of the last subsection, we assume $\ep=q^2-\ell p^2=-1$.
In this case, the algebraic integers $a_m,b_m$ 
in Definition~\ref{dfn:ambm} satisfies the following relations.
\begin{align*}
 \dfrac{b_{2 k - 1}}{a_{2 k - 1}}<\dfrac{\ell a_{2 k}}{b_{2 k}}
<\dfrac{b_{2 k - 1}}{a_{2 k - 1}}<\sqrt{\ell}
<\dfrac{\ell a_{2 k + 1}}{b_{2 k + 1}}<\dfrac{b_{2 k}}{a_{2 k}}
<\dfrac{\ell a_{2 k - 1}}{b_{2 k - 1}}
\;
(k\in\bbZ_{>0}),
\quad
\lim_{k\to\infty}\dfrac{b_k}{a_k}=\sqrt{\ell}.
\end{align*}
We regard $(a_m:b_m)$ and $(b_m:\ell a_m)$, $m \in {\bbZ}$ 
as elements of ${\bbP}^1({\bbR})$. 
Then the inhomogeneous coordinates of these points give a sequence
\begin{multline}
\label{eq:interval:case1:050}
 -\infty=\dfrac{b_0}{a_0}
<-\dfrac{\ell p}{q}=\dfrac{\ell a_{-1}}{b_{-1}}
<\dfrac{b_{-2}}{a_{-2}}<\cdots<-\sqrt{\ell}<\cdots
<\dfrac{\ell a_{-2}}{b_{-2}}<\dfrac{b_{-1}}{a_{-1}}=-\dfrac{q}{p}
\\
<\dfrac{\ell a_0}{b_0}=0<\dfrac{b_1}{a_1}=\dfrac{q}{p}
<\dfrac{\ell a_2}{b_2}<\cdots<\sqrt{\ell}<\cdots<\dfrac{b_2}{a_2}
<\dfrac{\ell a_1}{b_1}=\dfrac{\ell p}{q}<\dfrac{b_0}{a_0}=\infty,
\end{multline}
where we write the inhomogeneous coordinate of $(0:1)$ as $\infty$ or $-\infty$.

In our case $\ep=-1$, we have the cokernel presentation
\begin{align}\label{eq:cokernel}
0 \to E_1 \to E_2 \to F \to 0
\end{align}
for a general $F\in M_X^H(1,0,-\ell)$. 
By Lemma~\ref{lem:numerical-P(p,q)}, 
we can write $v(E_1)=(q^2,-\tfrac{\ell p q}{\sqrt{n}}H,\ell^2 p^2)$ 
and $v(E_2)=\ell (p^2,-\tfrac{p q}{\sqrt{n}} H,q^2)$. 
Now we set
\begin{align*}
X_1\seteq X,\quad X_0\seteq M_X^H(p^2,-\tfrac{p q}{\sqrt{n}}H,q^2).
\end{align*} 
Since
\begin{align*}
\begin{bmatrix}q & -\ell p \\ p & -q \end{bmatrix} \in G,
\end{align*}
there is a universal family $\bfE$ on $X_0\times X_1$ such that the FMT $\Phi_{X_0\to X_1}^\bfE \colon \bfD(X_0)\to \bfD(X_1)$ satisfies
\begin{align*}
\theta(\Phi_{X_0\to X_1}^\bfE)=
\begin{bmatrix}q & -\ell p \\ p & -q \end{bmatrix}.
\end{align*}
Then the cokernel presentation \eqref{eq:cokernel} is rewritten as 
\begin{align}\label{eq:039}
0 \to \Phi_{X_0\to X_1}^\bfE(B_0) \to \Phi_{X_0\to X_1}^\bfE(A_0) 
  \to F                            \to 0,
\end{align}
where $B_0\in M_{X_0}^{\widehat{H}}(1,0,0)$ and $A_0\in M_{X_0}^{\widehat{H}}(0,0,\ell)$. Here $\widehat{H}$ is the ample generator of $\NS(X_0)$. On the other hand, every $F\in M_X^H(1,0,-\ell)$ has the kernel presentation
\begin{align}\label{eq:040}
0 \to F \to B_1 \to A_1 \to 0,
\end{align}
with $B_1\in \widehat{X}=M_X^H(1,0,0)$ and $A_1\in M_X^H(0,0,\ell)$.

For $i\in\bbZ$, we set 
\begin{align*}
X_i\seteq X_{i\,{\rm mod}\,2},\quad
A_i\seteq A_{i\,{\rm mod}\,2},\quad
B_i\seteq B_{i\,{\rm mod}\,2}.
\end{align*}

\begin{dfn}\label{defn:E_{i,m}}
We define ${\bfE}_{i,m} \in {\bf D}(X_i \times X_{i+m})$, $m>0$
inductively by
\begin{align*}
{\bfE}_{i,0}[1-(-1)^i]&={\calO}_{\Delta},
\quad
\Phi_{X_{i+m-1} \to X_{i+m}}^{{\bfE}[1-(-1)^{i+m-1}]}\,
\Phi_{X_i \to X_{i+m-1}}^{{\bfE}_{i,m-1}[1-(-1)^i]}
=
\Phi_{X_i \to X_{i+m}}^{{\bfE}_{i,m}[1-(-1)^i]}\quad (m>0).
\end{align*}
Here the symbol $\Delta$ denotes the diagonal of $X_i\times X_i$.
\end{dfn}

For ${\bfE}_{i,m} \in \Coh(X_i \times X_{i+m})$
and the isomorphism 
\begin{align*}
\xi\colon X_i \times X_{i+m} \to X_{i+m} \times X_{i},\quad
         (x_i, x_{i+m}) \mapsto (x_{i+m},x_i),
\end{align*}
we also denote $\xi_*({\bfE}_{i,m}) \in \Coh( X_{i+m} \times X_{i})$
by ${\bfE}_{i,m}$.

\begin{lem}\label{lem:inverse}
\begin{align}
&\label{eq:E_{i,m}}
\Phi_{X_i \to X_{i+m}}^{{\bfE}_{i,m}[1-(-1)^i]}
=\Phi_{X_{i+m-1} \to X_{i+m}}^{{\bfE}[1-(-1)^{i+m-1}]}
  \dotsm 
  \Phi_{X_{i+1} \to X_{i+2}}^{{\bfE}[1-(-1)^{i+1}]}
  \Phi_{X_{i} \to X_{i+1}}^{{\bfE}[1-(-1)^{i}]},
  \quad (m>0),
\\[\baselineskip]
&\label{eq:E_{i+m,m}^{vee}}
\Phi_{X_i \to X_{i+m}}^{{\bfE}_{i+m,m}^{\vee}[1+(-1)^{i+m}]}
=\bigl(
   \Phi_{X_{i+m} \to X_{i+2m}}^{{\bfE}_{i+m,m}[1-(-1)^{i+m}]}
 \bigr)^{-1}
\\
&\hphantom{\Phi_{X_i \to X_{i+m}}^{{\bfE}_{i+m,m}^{\vee}[1+(-1)^{i+m}]}
}
=\Phi_{X_{i+m-1} \to X_{i+m}}^{{\bfE}^{\vee}[1+(-1)^{i+m}]}
  \dotsm 
  \Phi_{X_{i+1} \to X_{i+2}}^{{\bfE}^{\vee}[1+(-1)^{i+2}]}
  \Phi_{X_{i} \to X_{i+1}}^{{\bfE}^{\vee}[1+(-1)^{i+1}]},\quad (m>0).
\nonumber
\end{align}
Under the identification $\xi\colon X_i \times X_{i+m} \to X_{i+m} \times X_{i}$, we have
\begin{align}
\label{eq:E_{i+m,m}=E_{i,m}}
{\bfE}_{i+m,m}={\bfE}_{i,m}.
\end{align}
\end{lem}

\begin{proof}
\eqref{eq:E_{i,m}} is derived from Definition \ref{defn:E_{i,m}}.
The claim \eqref{eq:E_{i+m,m}^{vee}} follows from the following computation.
\begin{align*}
\Phi_{X_i \to X_{i+m}}^{{\bfE}_{i+m,m}^{\vee}[1+(-1)^{i+m}]}
=&\bigl(
   \Phi_{X_{i+m} \to X_{i+2m}}^{{\bfE}_{i+m,m}[1-(-1)^{i+m}]}
  \bigr)^{-1}
\\[2pt]
=&\bigl(\Phi_{X_{i+2m-1} \to X_{i+2m}}^{{\bfE}[1-(-1)^{i+2m-1}]}\,
   \Phi_{X_{i+2m-2} \to X_{i+2m-1}}^{{\bfE}[1-(-1)^{i+2m-2}]}\,
   \dotsm\,
   \Phi_{X_{i+m} \to X_{i+m+1}}^{{\bfE}[1-(-1)^{i+m}]}
  \bigr)^{-1}
\\[2pt]
=& \Phi_{X_{i+m+1} \to X_{i+m}}^{{\bfE}^{\vee}[1+(-1)^{i+m}]}\,
   \dotsm \,
   \Phi_{X_{i+2m-1} \to X_{i+2m-2}}^{{\bfE}^{\vee}[1+(-1)^{i+2m-2}]}\,
   \Phi_{X_{i+2m} \to X_{i+2m-1}}^{{\bfE}^{\vee}[1+(-1)^{i+2m-1}]}
\\[2pt]
=& \Phi_{X_{i+m-1} \to X_{i+m}}^{{\bfE}^{\vee}[1+(-1)^{i+m}]}\,
   \dotsm
   \Phi_{X_{i+1} \to X_{i+2}}^{{\bfE}^{\vee}[1+(-1)^{i+2}]}\,
   \Phi_{X_{i} \to X_{i+1}}^{{\bfE}^{\vee}[1+(-1)^{i+1}]}.
\end{align*}
We prove the last claim. By using Lemma \ref{lem:Serre}, we see that
\begin{align*}
\Phi_{X_i \to X_{i+m}}^{{\bfE}_{i,m}^{\vee}[1+(-1)^{i}]}
=&\calD_{X_{i+m}}\,
  \Phi_{X_i \to X_{i+m}}^{{\bfE}_{i,m}[1-(-1)^{i}]}\,
  \calD_{X_i}
\\[2pt]
=&\calD_{X_{i+m}}\,\Phi_{X_{i+m-1} \to X_{i+m}}^{{\bfE}[1-(-1)^{i+m-1}]}
  \dotsm
  \Phi_{X_{i+1} \to X_{i+2}}^{{\bfE}[1-(-1)^{i+1}]}
  \Phi_{X_{i}   \to X_{i+1}}^{{\bfE}[1-(-1)^{i}]}
  \calD_{X_i}
\\[2pt]
=&\Phi_{X_{i+m-1} \to X_{i+m}}^{{\bfE}^{\vee}[1+(-1)^{i+m-1}]}
  \dotsm 
  \Phi_{X_{i+1} \to X_{i+2}}^{{\bfE}^{\vee}[1+(-1)^{i+1}]}
  \Phi_{X_{i}   \to X_{i+1}}^{{\bfE}^{\vee}[1+(-1)^{i}]}.
\end{align*}
Hence 
$\Phi_{X_i \to X_{i+m}}^{{\bfE}_{i+m,m}^{\vee}}
=\Phi_{X_i \to X_{i+m}}^{{\bfE}_{i,m}^{\vee}}$,
which implies \eqref{eq:E_{i+m,m}=E_{i,m}}.
\end{proof}

If $m$ is even, then
\begin{align*}
\Phi_{X_i \to X_i}^{{\bfE}_{i,m}^{\vee}[1+(-1)^i]}
=\bigl(\Phi_{X_i \to X_i}^{{\bfE}_{i,m}[1-(-1)^i]}\bigr)^{-1}
=\bigl(\Phi_{X_i \to X_i}^{{\bfE}_{i,2}[1-(-1)^i]}\bigr)^{-m/2}.
\end{align*}
If $m$ is odd, then
\begin{align*}
\Phi_{X_i \to X_{i+1}}^{{\bfE}_{i,m}^{\vee}[1+(-1)^{i+1}]}
=&\Phi_{X_i \to X_{i+1}}^{{\bfE}_{i+1,m}^{\vee}[1+(-1)^{i+1}]}
= \bigl(\Phi_{X_{i+1} \to X_{i}}^{{\bfE}_{i+1,m}[1-(-1)^{i+1}]}\bigr)^{-1}
\\[2pt]
=&\bigl(
   \Phi_{X_{i+1} \to X_{i+1}}^{{\bfE}_{i+1,2}[1-(-1)^{i+1}]}
  \bigr)^{-(m-1)/2}
  \bigl(\Phi_{X_{i+1} \to X_{i}}^{{\bfE}[1-(-1)^{i+1}]}\bigr)^{-1}.
\end{align*}

The next lemma compute the cohomological action of $\bfE_{i,m}$.

\begin{lem}\label{lem:cohom_e_im}
(1)
For $m\in\bbZ_{\ge0}$ we have
\begin{align*}
\theta(\Phi^{\bfE_{0,m}}_{X_0\to X_m})
&=
 \left\{
 \begin{array}{l l}
 \pm
 \begin{bmatrix}q   & \ell p   \\ p   & q   \end{bmatrix}^m
=\pm
 \begin{bmatrix}b_m & \ell a_m \\ a_m & b_m \end{bmatrix}
& m \in 2\bbZ_{\ge0},
\\
\\
 \pm
 \begin{bmatrix}q   & \ell p   \\ p   & q   \end{bmatrix}^m
 \begin{bmatrix}1   & 0        \\ 0   & -1  \end{bmatrix}
=\pm
 \begin{bmatrix}b_m &-\ell a_m \\ a_m &-b_m \end{bmatrix}
& m \in 2\bbZ_{\ge0}+1,
\end{array}
\right.
\\
\theta(\Phi^{\bfE_{1,m}[2]}_{X_1\to X_{1+m}})
&=
\left\{
 \begin{array}{l l}
 \pm
 \begin{bmatrix}-q  & \ell p   \\ p   & -q  \end{bmatrix}^m
=\pm
 \begin{bmatrix}-b_m&-\ell a_m \\ a_m &-b_m \end{bmatrix}
& m \in 2\bbZ_{\ge0},
\\
\\
 \pm
 \begin{bmatrix}-q  & \ell p   \\ p   & -q  \end{bmatrix}^m
 \begin{bmatrix}1   & 0        \\ 0   & -1  \end{bmatrix}
=\pm
 \begin{bmatrix}-b_m&-\ell a_m \\ a_m & b_m \end{bmatrix}
& m \in 2\bbZ_{\ge0}+1.
\end{array}
\right.
\end{align*}
In these equations, the signs in front of the matrices are irrelevant,
since the images of $\theta$ in $\PSL(2,\bbR)$.

(2)
$\bfE_{i,m}$ is a locally free sheaf for $m\in\bbZ_{>0}$.
\end{lem}

\begin{proof}
(1)
By Lemma~\ref{lem:matrix} we have
\begin{align*}
\theta(\Phi_{X_1\to X_0}^{\bfE})
=\pm\begin{bmatrix}-q &-\ell p \\ p & q \end{bmatrix}.
\end{align*}
Hence we find that
\begin{align*}
&\theta(\Phi_{X_1\to X_0}^\bfE \Phi_{X_0\to X_1}^\bfE)
=\theta(\Phi_{X_0\to X_1}^\bfE)\theta(\Phi_{X_1\to X_0}^\bfE)
=\pm\begin{bmatrix}q&\ell p \\ p & q\end{bmatrix}^2,
\\
&
\theta(\Phi_{X_0\to X_1}^\bfE \Phi_{X_1\to X_0}^\bfE)
=\pm\begin{bmatrix}q&-\ell p \\ -p & q\end{bmatrix}^2.
\end{align*}
{}From these relations we obtain the conclusion.

(2)
Assume that $\bfE_{i,m-1}$ is a locally free sheaf. By Lemma \ref{lem:cohom_e_im}~(1), we have
\begin{align*}
\begin{array}{l l}
\mu(\bfE_{i,m}|_{X_i \times \{x_{i+m}\}})   
 &=(-1)^i\dfrac{b_m}{a_m}\sqrt{n},\\[6pt]
\mu(\bfE_{i,m-1}|_{\{x_i\} \times X_{i+m-1}})
 &=(-1)^{i+m-1}\dfrac{b_{m-1}}{a_{m-1}}\sqrt{n},\\[6pt]
\mu(\bfE|_{X_{i+m-1} \times \{x_{i+m}\}})
 &=(-1)^{i+m-1}\dfrac{q}{p}\sqrt{n},
\end{array}
\end{align*}
Here $x_i$ is a closed point of $X_i$. 
Propositions \ref{prop:semihom-fmt} (1) and \eqref{eq:interval:case1:050} 
yield that if $i+m-1$ is even then 
$\Phi_{X_i \to X_{i+m}}^{\bfE_{i,m}}(\bbC_{x_i})
=\Phi_{X_{i+m-1} \to X_{i+m}}^{\bfE}
 \Phi_{X_i \to X_{i+m-1}}^{\bfE_{i,m-1}}(\bbC_{x_i})$ 
is a locally free sheaf. In the same manner if $i+m-1$ is odd, 
then 
$\Phi_{X_i \to X_{i+m}}^{\bfE_{i,m}}(\bbC_{x_i})
=\Phi_{X_{i+m-1} \to X_{i+m}}^{\bfE[2]}
 \Phi_{X_i \to X_{i+m-1}}^{\bfE_{i,m-1}}(\bbC_{x_i})$ 
is a locally free sheaf. 
Therefore $\bfE_{i,m}$ is also a locally free sheaf. 
Thus the conclusion holds.
\end{proof}

The next weak index theorem is 
the fundamental result for our computation of FMTs.

\begin{lem}\label{lem:A-B}
For $m \geq 0$, we have
\begin{align*}
\begin{array}{l l l l l}
&
&\Phi_{X_0 \to X_{m}}^{{\bfE}_{0,m}}(A_0),
&\Phi_{X_0 \to X_{m}}^{{\bfE}_{0,m}}(A_0^{\vee}[2])
&\in \Coh(X_{m}),
\\[6pt]
 \Phi_{X_0 \to X_{m}}^{{\bfE}_{0,m}}(B_0),
&\Phi_{X_0 \to X_{m}}^{{\bfE}_{0,m}}(B_0^{\vee}),
&\Phi_{X_0 \to X_{m}}^{{\bfE}_{0,m}^{\vee}[2]}(B_0),
&\Phi_{X_0 \to X_{m}}^{{\bfE}_{0,m}^{\vee}[2]}(B_0^{\vee})
&\in \Coh(X_{m}),
\\[6pt]
&
&\Phi_{X_1 \to X_{1+m}}^{{\bfE}_{1,m}}(A_1),
&\Phi_{X_1 \to X_{1+m}}^{{\bfE}_{1,m}}(A_1^{\vee}[2])
&\in \Coh(X_{1+m}),
\\[6pt]
 \Phi_{X_1 \to X_{1+m}}^{{\bfE}_{1,m}[2]}(B_1),
&\Phi_{X_1 \to X_{1+m}}^{{\bfE}_{1,m}[2]}(B_1^{\vee}),
&\Phi_{X_1 \to X_{1+m}}^{{\bfE}_{1,m}^{\vee}}(B_1),
&\Phi_{X_1 \to X_{1+m}}^{{\bfE}_{1,m}^{\vee}}(B_1^{\vee}) 
&\in \Coh(X_{1+m}).
\end{array}
\end{align*}
\end{lem}

\begin{proof}
We may assume that $m>0$.
Since $A_i, A_i^{\vee}[2]$, $i=0,1$ are 0-dimensional sheaves, 
the claims hold for $A_i, A_i^{\vee}[2]$, $i=0,1$.
Since $B_i, B_i^{\vee}$, $i=0,1$ are line bundles of degree 0,  
$\mu({\bfE}_{j,m}|_{X_j \times \{ x_{j+m} \}})=(-1)^j(b_m/a_m) \sqrt{n}$
implies that $\WIT_0$ (resp. $\WIT_2$) holds for $B_0, B_0^{\vee}$
with respect to $\Phi_{X_0 \to X_{m}}^{{\bfE}_{0,m}}$
(resp. $\Phi_{X_0 \to X_{m}}^{{\bfE}_{0,m}^{\vee}}$) and
$\WIT_2$ (resp. $\WIT_0$) holds for  
$B_1, B_1^{\vee}$
with respect to $\Phi_{X_1 \to X_{1+m}}^{{\bfE}_{1,m}}$
(resp. $\Phi_{X_1 \to X_{1+m}}^{{\bfE}_{1,m}^{\vee}}$).
\end{proof}

{}From \eqref{eq:039} and \eqref{eq:040} we have
\begin{align*}
\tag{$\#$}
0 \to \Phi_{X_0\to X_1}^\bfE(B_0) \to \Phi_{X_0\to X_1}^\bfE(A_0)
  \to B_1                          \to A_1                          
  \to 0.
\end{align*}
Applying $\Phi_{X_1\to X_0}^{\bfE^\vee[2]}$ to this sequence, we have
\begin{align*}
\tag{$\#^*$}
0 \to \Phi_{X_1\to X_0}^{\bfE^\vee}(B_1)
  \to \Phi_{X_0\to X_1}^{\bfE^\vee}(A_1)
  \to B_0 \to A_0
  \to 0.
\end{align*}

For later use, we prepare the next lemma.

\begin{lem}\label{lem:case1:012}
\begin{align*}
\begin{array}{l l l l l l l}
 \Phi_{X_0\to X_1}^{\bfE_{0,2 k +1}}  \calD_{X_0} \Phi_{X_1\to X_0}^{\bfE^\vee}
&=
&\Phi_{X_1\to X_1}^{\bfE_{1,2 k+2}[2]} \calD_{X_1},
&
&\Phi_{X_0\to X_1}^{\bfE_{0,2 k +1}}          \Phi_{X_1\to X_0}^{\bfE^\vee}
&=
&\Phi_{X_1\to X_1}^{\bfE_{1,2 k}},
\\[6pt]
 \Phi_{X_1\to X_1}^{\bfE_{1,2 k}[2]} \calD_{X_1} \Phi_{X_0\to X_1}^{\bfE}
&=
&\Phi_{X_0\to X_1}^{\bfE_{0,2 k-1}[2]} \calD_{X_0},
&
&\Phi_{X_1\to X_1}^{\bfE_{1,2 k}[2]}         \Phi_{X_0\to X_1}^{\bfE}
&=
&\Phi_{X_0\to X_1}^{\bfE_{0,2 k +1}}.
\end{array}
\end{align*}
\end{lem}
\begin{proof}
Recall Lemma~\ref{lem:Serre} and the definition of $\bfE_{i,m}$ (Definition \ref{defn:E_{i,m}}). One can compute that
\begin{align*}
&
\Phi_{X_0 \to X_1}^{{\bfE}_{0,2k+1}}\calD_{X_0}
\Phi_{X_1 \to X_0}^{{\bfE}^{\vee}}
=
\Phi_{X_0 \to X_1}^{{\bfE}_{0,2k+1}}
\Phi_{X_1 \to X_0}^{{\bfE}[2]}\calD_{X_0}=
\Phi_{X_1 \to X_1}^{{\bfE}_{1,2k+2}[2]}\calD_{X_0},
\\
&
\Phi_{X_0 \to X_1}^{{\bfE}_{0,2k+1}}
\Phi_{X_1 \to X_0}^{{\bfE}^{\vee}}
=
\Phi_{X_0 \to X_1}^{{\bfE}_{0,2k+1}}
(\Phi_{X_0 \to X_1}^{{\bfE}})^{-1}[-2]
=\Phi_{X_1 \to X_1}^{{\bfE}_{1,2k}[2]}[-2],
\\
&
\Phi_{X_1 \to X_1}^{{\bfE}_{1,2k}[2]} \calD_{X_1}
\Phi_{X_0 \to X_1}^{{\bfE}}
=
\Phi_{X_1 \to X_1}^{{\bfE}_{1,2k}[2]}
\Phi_{X_0 \to X_1}^{{\bfE}^{\vee}[2]} \calD_{X_0}
=
\Phi_{X_1 \to X_1}^{{\bfE}_{1,2k}[2]}
(\Phi_{X_1 \to X_0}^{{\bfE}[2]})^{-1}[2]\, \calD_{X_0}
=
\Phi_{X_0 \to X_1}^{{\bfE}_{0,2k-1}[2]}\, \calD_{X_0}.
\end{align*}
The last equality follows from the second one.
\end{proof}

\begin{dfn}\label{defn:F_m}
For  
\begin{align*}
F\seteq\Ker(B_1 \to A_1)
 \cong \Coker(\Phi_{X_0 \to X_1}^{{\bfE}}(B_0) 
              \to \Phi_{X_0 \to X_1}^{{\bfE}}(A_0)),
\end{align*}
we define the \emph{tame complex} $F_m^\bullet$ by

\begin{align*}
F_m^\bullet\seteq
\left\{
\begin{array}{l l l l}
\bigl(\Phi_{X_1 \to X_1}^{{\bfE}_{1,2}[2]}\bigr)^{(m+1)/2}(F^{\vee})
&=
\Phi_{X_1 \to X_1}^{{\bfE}_{1,2k}[2]}(F^{\vee})
&m=2k-1,
&k \geq 1,
\\[4pt]
\bigl(\Phi_{X_1 \to X_1}^{{\bfE}_{1,2}[2]}\bigr)^{m/2}(F)
&=\Phi_{X_1 \to X_1}^{{\bfE}_{1,2k}[2]}(F)
&m=2k,
&k \geq 0,
\\[4pt]
\bigl(\Phi_{X_1 \to X_1}^{{\bfE}_{1,2}[2]}\bigr)^{(m+1)/2}(F^{\vee})
&=\Phi_{X_1 \to X_1}^{{\bfE}_{1,2k}^{\vee}}(F^{\vee})
&m=-2k-1,
&k \geq 0,
\\[4pt]
\bigl(\Phi_{X_1 \to X_1}^{{\bfE}_{1,2}[2]}\bigr)^{m/2}(F)
&=\Phi_{X_1 \to X_1}^{{\bfE}_{1,2k}^{\vee}}(F)
&m=-2k,
&k \geq 1.
\end{array}
\right.
\end{align*}
\end{dfn}
$F_m^\bullet$ satisfies $v(F_m^\bullet)=(1,0,-\ell)$. 
We can also express them in terms of $A_0, B_0$.
\begin{lem}
\label{lem:another}
\begin{align*}
F_m^\bullet=
\left\{
\begin{array}{l l}
\Phi_{X_0 \to X_1}^{{\bfE}_{0,2k-1}[2]}\bigl([B_0 \to A_0]^{\vee}\bigr)
&m=2k-1,
\\[4pt]
\Phi_{X_0 \to X_1}^{{\bfE}_{0,2k+1}}\bigl([B_0 \to A_0]\bigr)
&m=2k,
\\[4pt]
\Phi_{X_0 \to X_1}^{{\bfE}_{1,2k+1}^{\vee}[2]}\bigl([B_0 \to A_0]^{\vee}\bigr)
&m=-2k-1,
\\[4pt]
\Phi_{X_0 \to X_1}^{{\bfE}_{1,2k-1}^{\vee}}\bigl([B_0 \to A_0]\bigr)
&m=-2k.
\end{array}
\right.
\end{align*}
\end{lem}

By Lemmas \ref{lem:A-B}, \ref{lem:another} and Definition \ref{defn:F_m},
we have the following description of $F_m^\bullet$, 
from which follows the first part of Proposition~\ref{prop:explicit} 
for the case $\ep=-1$.

\begin{prop}
\begin{enumerate}
\item
If $m=2k, k \geq 1$, then
$\Phi_{X_0 \to X_1}^{{\bfE}_{0,2k+1}}(B_{0})$ and 
$\Phi_{X_0 \to X_1}^{{\bfE}_{0,2k+1}}(A_0)$ are locally free sheaves,
and
\begin{align*}
&F_{2k}^\bullet
 =\bigl[\Phi_{X_0 \to X_1}^{{\bfE}_{0,2k+1}}(B_{0}) \to 
 \Phi_{X_0 \to X_1}^{{\bfE}_{0,2k+1}}(A_0)\bigr],
\\
&H^i(F_{2k}^\bullet)=
\left\{
 \begin{array}{l l}
 \Phi_{X_1 \to X_1}^{{\bfE}_{1,2k}}(A_1) & i=-1,
 \\[4pt]
 \Phi_{X_1 \to X_1}^{{\bfE}_{1,2k}[2]}(B_1) & i=0,
 \\[4pt]
 0 & i \not= -1,0.
 \end{array}
\right.
\end{align*}
In particular, we have an exact sequence
\begin{align*}
\tag{$\#_{2k}$}
0 \to \Phi_{X_1 \to X_1}^{{\bfE}_{1,2k}}(A_{1}) 
  \to \Phi_{X_0 \to X_1}^{{\bfE}_{0,2k+1}}(B_{0})
  \to \Phi_{X_0 \to X_1}^{{\bfE}_{0,2k+1}}(A_0)
  \to \Phi_{X_{1} \to X_1}^{{\bfE}_{1,2k}[2]}(B_{1})
  \to 0.
\end{align*}
The slopes of these four sheaves are 
\begin{align*}
  -\frac{b_{2k}       }{a_{2k}}  \sqrt{n}
 <-\frac{\ell a_{2k+1}}{b_{2k+1}}\sqrt{n}
 <-\frac{b_{2k+1}     }{a_{2k+1}}\sqrt{n}
 <-\frac{\ell a_{2k}  }{b_{2k}}  \sqrt{n}.
\end{align*}
\item
If $m=2k+1, k \geq 0$, then
$\Phi_{X_{1} \to X_1}^{{\bfE}_{1,2k+2}[2]}(A_{1}^{\vee})$ and 
$\Phi_{X_{1} \to X_1}^{{\bfE}_{1,2k+2}[2]}(B_{1}^{\vee})$ are
locally free sheaves, and 
\begin{align*}
&F_{2k+1}^\bullet=
 \bigl[
  \Phi_{X_{1} \to X_1}^{{\bfE}_{1,2k+2}[2]}(A_{1}^{\vee})\to 
  \Phi_{X_{1} \to X_1}^{{\bfE}_{1,2k+2}[2]}(B_{1}^{\vee})
 \bigr],
\\
&H^i(F_{2k+1}^\bullet)=
\left\{
\begin{array}{l l}
\Phi_{X_0 \to X_1}^{{\bfE}_{1,2k+1}}(B_0^{\vee})
& i=-1,
\\[4pt]
\Phi_{X_0 \to X_1}^{{\bfE}_{1,2k+1}}(A_0^{\vee}[2])
& i=0,
\\[4pt]
0
& i \neq -1,0.
\end{array}
\right.
\end{align*}
In particular, we have an exact sequence
\begin{align*}
\tag{$\#_{2k+1}$}
0 \to \Phi_{X_{0} \to X_1}^{{\bfE}_{0,2k+1}}(B_{0}^{\vee}) 
  \to \Phi_{X_{1} \to X_1}^{{\bfE}_{1,2k+2}[2]}(A_{1}^{\vee})
  \to \Phi_{X_{1} \to X_1}^{{\bfE}_{1,2k+2}[2]}(B_{1}^{\vee})
  \to \Phi_{X_{0} \to X_1}^{{\bfE}_{0,2k+1}}(A_{0}^{\vee}[2])
  \to 0.
\end{align*}
The slopes of these four sheaves are 
\begin{align*}
 -\frac{\ell a_{2k+1}}{b_{2k+1}}\sqrt{n}
<-\frac{b_{2k+2}}     {a_{2k+2}} \sqrt{n}
<-\frac{\ell a_{2k+2}}{b_{2k+2}}\sqrt{n}
<-\frac{b_{2k+1}}     {a_{2k+1}} \sqrt{n}.
\end{align*}
\item
If $m=-2k-1, k \geq 0$, then
$\Phi_{X_0 \to X_1}^{{\bfE}_{1,2k+1}^{\vee}[2]}(A_{0}^{\vee})$, 
$\Phi_{X_0 \to X_1}^{{\bfE}_{1,2k+1}^{\vee}[2]}(B_0^{\vee})$ 
are locally free sheaves, and
\begin{align*}
&F_{-2k-1}^\bullet[1]
=\bigl[
  \Phi_{X_0 \to X_1}^{{\bfE}_{1,2k+1}^{\vee}[2]}(A_{0}^{\vee}) \to 
  \Phi_{X_0 \to X_1}^{{\bfE}_{1,2k+1}^{\vee}[2]}(B_0^{\vee})
 \bigr],
\\
&H^i(F_{-2k-1}^\bullet[1])
=
\left\{
\begin{array}{l l}
\Phi_{X_1 \to X_1}^{{\bfE}_{1,2k}^{\vee}}(B_{1}^{\vee})& i=-1,
\\[4pt]
\Phi_{X_{1} \to X_1}^{{\bfE}_{1,2k}^{\vee}}(A_{1}^{\vee}[2])& i=0,
\\[4pt]
0 & i \ne -1,0.
\end{array}
\right.
\end{align*}
In particular, we have an exact sequence
\begin{align*}
\tag{$\#_{-2k-1}$}
0 \to \Phi_{X_1 \to X_1}^{{\bfE}_{1,2k}^{\vee}}(B_{1}^{\vee}) 
  \to \Phi_{X_0 \to X_1}^{{\bfE}_{1,2k+1}^{\vee}[2]}(A_{0}^{\vee}) 
  \to \Phi_{X_0 \to X_1}^{{\bfE}_{1,2k+1}^{\vee}[2]}(B_0^{\vee})
  \to \Phi_{X_{1} \to X_1}^{{\bfE}_{1,2k}^{\vee}}(A_{1}^{\vee}[2])
  \to 0.
\end{align*}
The slopes of these four sheaves are 
\begin{align*}
 \frac{\ell a_{2k}}  {b_{2k}}  \sqrt{n}
<\frac{b_{2k+1}}     {a_{2k+1}}\sqrt{n}
<\frac{\ell a_{2k+1}}{b_{2k+1}}\sqrt{n}
<\frac{b_{2k}}       {a_{2k}}  \sqrt{n}.
\end{align*}

\item
If $m=-2k-2, k \geq 0$, then 
$\Phi_{X_{1} \to X_1}^{{\bfE}_{1,2k+2}^{\vee}}(B_{1})$ and
$\Phi_{X_{1} \to X_1}^{{\bfE}_{1,2k+2}^{\vee}}(A_{1})$ are locally free
sheaves, and 
\begin{align*}
&F_{-2k-2}^\bullet=
\bigl[
 \Phi_{X_{1} \to X_1}^{{\bfE}_{1,2k+2}^{\vee}}(B_{1})
 \to 
 \Phi_{X_{1} \to X_1}^{{\bfE}_{1,2k+2}^{\vee}}(A_{1})
\bigr],
\\
&H^i(F_{-2k-2}^\bullet[1])=
\left\{
\begin{array}{l l}
\Phi_{X_{0} \to X_1}^{{\bfE}_{1,2k+1}^{\vee}}(A_{0})& i=-1,
\\[4pt]
\Phi_{X_{0} \to X_1}^{{\bfE}_{1,2k+1}^{\vee}[2]}(B_{0})& i=0,
\\[4pt]
0 & i \ne -1,0.
\end{array}
\right.
\end{align*}
In particular, we have an exact sequence
\begin{align*}
\tag{$\#_{-2k-2}$}
0 \to \Phi_{X_{0} \to X_1}^{{\bfE}_{1,2k+1}^{\vee}}(A_{0}) 
  \to \Phi_{X_{1} \to X_1}^{{\bfE}_{1,2k+2}^{\vee}}(B_{1})
  \to \Phi_{X_{1} \to X_1}^{{\bfE}_{1,2k+2}^{\vee}}(A_{1})
  \to \Phi_{X_{0} \to X_1}^{{\bfE}_{1,2k+1}^{\vee}[2]}(B_{0})
  \to 0.
\end{align*}
The slopes of these four sheaves are 
\begin{align*}
 \dfrac{b_{2k+1} }{a_{2k+1}}\sqrt{n}
<\dfrac{\ell a_{2k+2}}{b_{2k+2}}\sqrt{n}
<\dfrac{b_{2k+2} }{a_{2k+2}}\sqrt{n}
<\dfrac{\ell a_{2k+1}}{b_{2k+1}}\sqrt{n}.
\end{align*}
\end{enumerate}
\end{prop}

Now we show the second part of Proposition~\ref{prop:explicit} 
for the case $\ep=-1$.
By Definition~\ref{defn:F_m} and Lemma~\ref{lem:another}, 
the relations of $(\#_m)$, $m \in {\bbZ}$ are summarized as follows.
\begin{align}
\label{eq:relation-sharp}
\begin{array}{l l l l}
F_{2k}^\bullet&=&
 \Phi_{X_1 \to X_1}^{{\bfE}_{1,2k}[2]}\, \calD_{X_1}\,
 \bigl(\Phi_{X_1 \to X_1}^{{\bfE}_{1,2k}[2]}\bigr)^{-1}
 (F_{2k-1}^\bullet)
&k \in\bbZ_{\ge 1},
\\[4pt]
F_{2k+1}^\bullet&=&
 \Phi_{X_0 \to X_1}^{{\bfE}_{0,2k+1}[2]}\, \calD_{X_0}\,
 \bigl(\Phi_{X_0 \to X_1}^{{\bfE}_{0,2k+1}}\bigr)^{-1}(F_{2k}^\bullet)
&k \in\bbZ_{\ge 0},
\\[4pt]
F_{-2k-1}^\bullet&=&
 \Phi_{X_1 \to X_1}^{{\bfE}_{1,2k}^{\vee}}\, \calD_{X_1}\,
 \bigl(\Phi_{X_1 \to X_1}^{{\bfE}_{1,2k}^{\vee}}\bigr)^{-1}(F_{-2k}^\bullet)
&k \in\bbZ_{\ge 0},
\\[4pt]
F_{-2k-2}^\bullet&=&
 \Phi_{X_0 \to X_1}^{{\bfE}_{1,2k+1}^{\vee}}\, \calD_{X_0}\,
 \bigl(\Phi_{X_0 \to X_1}^{{\bfE}_{1,2k+1}^{\vee}[2]}\bigr)^{-1}
 (F_{-2k-1}^\bullet)
&k \in\bbZ_{\ge 0}.
\end{array}
\end{align}
These are nothing but the relations $\Psi(F_m^\bullet)=F_{m+1}^\bullet$ 
($m \in \bbZ$) indicated 
in the second statement of Proposition~\ref{prop:explicit}. Thus the proof of Proposition~\ref{prop:explicit} is completed.

\subsection{Application of the tame complex.}

For a FMT
$\Phi_{X_1 \to X'}^{{\bfG}}\colon{\bf D}(X_1) \to {\bf D}(X')$, 
we set $\lambda \seteq \mu({\bfG}_{x'})/\sqrt{n}$.
If $-\ell p / q<-\lambda<-q/p$, then
$\Phi_{X_1 \to X'}^{{\bfG}}(F_0^\bullet)$
is not a sheaf for all $F_0^\bullet$.

\begin{dfn}
We set
\begin{align*}
\begin{array}{l l l l l l}
I_0&\seteq&
 (0,\tfrac{b_1}{a_1}] 
 \cup 
 (\tfrac{\ell a_1}{b_1},\infty],
&
I_{-1}&\seteq&
 (-\infty,-\tfrac{\ell a_1}{b_1}] 
 \cup 
 (-\tfrac{b_1}{a_1},0],
\\[4pt]
I_{2k-1}&\seteq& 
 (\tfrac{b_{2k-1}}{a_{2k-1}},\tfrac{\ell a_{2k}}{b_{2k}}]
 \cup 
 (\tfrac{b_{2k}}{a_{2k}},\tfrac{\ell a_{2k-1}}{b_{2k-1}}],
&
I_{-2k-1}&\seteq&
 (-\tfrac{b_{2k}}{a_{2k}},-\tfrac{\ell a_{2k+1}}{b_{2k+1}}]
 \cup 
 (-\tfrac{b_{2k+1}}{a_{2k+1}},-\tfrac{\ell a_{2k}}{b_{2k}}],
\\[4pt]
I_{2k}&\seteq& 
 (\tfrac{\ell a_{2k}}{b_{2k}},\tfrac{b_{2k+1}}{a_{2k+1}}]
 \cup 
 (\tfrac{\ell a_{2k+1}}{b_{2k+1}},\tfrac{b_{2k}}{a_{2k}}],
&
I_{-2k}&\seteq&
 (-\tfrac{\ell a_{2k-1}}{b_{2k-1}},-\tfrac{b_{2k}}{a_{2k}}]
 \cup 
 (-\tfrac{\ell a_{2k}}{b_{2k}},-\tfrac{b_{2k-1}}{a_{2k-1}}].
\end{array}
\end{align*}
For $I=\coprod_i (s_i, t_i]$, we denote $I^*\seteq\coprod_i [s_i,t_i)$. 
\end{dfn}

By \eqref{eq:interval:case1:050}, we have decompositions 
${\bbP}^1({\bbR}) \setminus\{\pm\sqrt{\ell}\} =
\coprod_{m \in  {\bbZ}}I_m=
\coprod_{m \in {\bbZ}}I_m^*$.

\begin{thm}\label{thm:csv:1}
(1)
If $\lambda \in I_m$, then
$\Phi_{X_1 \to X'}^{{\bfG}}(F_m^\bullet)$ is a coherent sheaf.

(2)
If $\lambda \in I_m^*$, then
$\calD_{X'}\Phi_{X_1 \to X'}^{{\bfG}}(F_m^\bullet)=
\Phi_{X_1 \to X'}^{{\bfG}^{\vee}[2]}(F_m^{\bullet \vee})$ 
is a coherent sheaf.
\end{thm}

\begin{proof}
We only prove (1). Assume that $m=2k$. Then
$-\frac{b_{2k}}{a_{2k}} \leq -\lambda \leq -\frac{\ell a_{2k+1}}{b_{2k+1}}$
or  
$-\frac{b_{2k+1}}{a_{2k+1}} \leq -\lambda \leq -\frac{\ell a_{2k}}{b_{2k}}$.
By $(\#_{2k})$ and Lemma~\ref{lem:sheaf-criterion}, we get our claim.
If $m=2k-1$, then our claim also follows from 
$(\#_{2k-1})$ and Lemma~\ref{lem:sheaf-criterion}.
\end{proof}

\begin{rmk}\label{rmk:calF}
We will explain the operation $\calF$ in Fact \ref{fact:yoshioka-lem2-2}.
For simplicity, we assume that 
$\lambda \in (\tfrac{\ell a_{2k+1}}{b_{2k+1}},\tfrac{b_{2k}}{a_{2k}}]
=I_{2k} \cap (\sqrt{\ell},\infty)$.
Then we have 
\begin{align}
\label{eq:G-A-B}
\begin{array}{l l l}
\Phi_{X_1 \to X'}^{{\bfG}}
\Phi_{X_1 \to X_1}^{{\bfE}_{1,2i}[2]}(A_{1}),\;
&
\Phi_{X_1 \to X'}^{{\bfG}}
\Phi_{X_1 \to X_1}^{{\bfE}_{1,2i}[2]}(A_{1}^{\vee}[2]) 
&\in \Coh(X')\quad i \leq k,
\\[4pt]
\Phi_{X_1 \to X'}^{{\bfG}}
\Phi_{X_0 \to X_1}^{{\bfE}_{0,2i+1}[2]}(B_{0}),\;
&
\Phi_{X_1 \to X'}^{{\bfG}}
\Phi_{X_0 \to X_1}^{{\bfE}_{0,2i+1}[2]}(B_{0}^{\vee})
&\in \Coh(X')\quad i < k,
\\[4pt]
\Phi_{X_1 \to X'}^{{\bfG}}
\Phi_{X_0 \to X_1}^{{\bfE}_{0,2k+1}}(B_{0}),\;
&
\Phi_{X_1 \to X'}^{{\bfG}}
\Phi_{X_0 \to X_1}^{{\bfE}_{0,2k+1}}(B_{0}^{\vee})
&\in \Coh(X'),
\\[4pt]
\Phi_{X_1 \to X'}^{{\bfG}}
\Phi_{X_0 \to X_1}^{{\bfE}_{0,2i+1}}(A_{0}),\;
&
\Phi_{X_1 \to X'}^{{\bfG}}
\Phi_{X_0 \to X_1}^{{\bfE}_{0,2i+1}}(A_{0}^{\vee}[2])
&\in \Coh(X')\quad 0 \leq i, 
\\[4pt]
\Phi_{X_1 \to X'}^{{\bfG}}
\Phi_{X_1 \to X_1}^{{\bfE}_{1,2i}[2]}(B_{1}),\;
&
\Phi_{X_1 \to X'}^{{\bfG}}
\Phi_{X_1 \to X_1}^{{\bfE}_{1,2i}[2]}(B_{1}^{\vee}) 
&\in \Coh(X')\quad 0 \leq i.
\end{array}
\end{align}
By using \eqref{eq:G-A-B}, 
we have the following description of $F_m^\bullet$, $m=2j,2j+1$, $j<k$. 
\begin{itemize}
\item
$\Phi_{X_1 \to X'}^{{\bfG}}
 \Phi_{X_{1} \to X_1}^{{\bfE}_{1,2j}[2]}(B_{1})$ 
and  
$\Phi_{X_1 \to X'}^{{\bfG}}
 \Phi_{X_1 \to X_1}^{{\bfE}_{1,2j}[2]}(A_{1})$ 
are locally free sheaves, and
\begin{align*}
&
\Phi_{X_1 \to X'}^{{\bfG}}(F_{2j}^\bullet[1])
=
\bigl[\Phi_{X_1 \to X'}^{{\bfG}}
      \Phi_{X_{1} \to X_1}^{{\bfE}_{1,2j}[2]}(B_{1})
      \to 
      \Phi_{X_1 \to X'}^{{\bfG}}
      \Phi_{X_1 \to X_1}^{{\bfE}_{1,2j}[2]}(A_{1})
\bigr],
\\
&
H^i(\Phi_{X_1 \to X'}^{{\bfG}}(F_{2j}^\bullet[1]))=\left\{
\begin{array}{l l}
\Phi_{X_1 \to X'}^{{\bfG}}
\Phi_{X_0 \to X_1}^{{\bfE}_{0,2j+1}}(A_0)
& i=-1,
\\[4pt]
\Phi_{X_1 \to X'}^{{\bfG}}
\Phi_{X_0 \to X_1}^{{\bfE}_{0,2j+1}[2]}(B_0)
& i=0,
\\[4pt]
0, & i \not= -1,0.
\end{array}
\right.
\end{align*}

\item
$\Phi_{X_1 \to X'}^{{\bfG}}
 \Phi_{X_{0} \to X_1}^{{\bfE}_{0,2j+1}}(A_{0}^{\vee}[2])$
and 
$\Phi_{X_1 \to X'}^{{\bfG}}
 \Phi_{X_{0} \to X_1}^{{\bfE}_{0,2j+1}[2]}(B_{0}^{\vee})$
are locally free sheaves, and 
\begin{align*}
&
\Phi_{X_1 \to X'}^{{\bfG}}(F_{2j+1}^\bullet[1])
=\bigl[
 \Phi_{X_1 \to X'}^{{\bfG}}
 \Phi_{X_{0} \to X_1}^{{\bfE}_{0,2j+1}}(A_{0}^{\vee}[2])
 \to
 \Phi_{X_1 \to X'}^{{\bfG}}
 \Phi_{X_{0} \to X_1}^{{\bfE}_{0,2j+1}[2]}(B_{0}^{\vee})
\bigr],
\\
&
H^i(\Phi_{X_1 \to X'}^{{\bfG}}(F_{2j+1}^\bullet[1]))=\left\{
\begin{array}{l l}
\Phi_{X_1 \to X'}^{{\bfG}}
\Phi_{X_1 \to X_1}^{{\bfE}_{1,2j+2}[2]}(B_1^{\vee})
&i=-1,
\\[4pt]
\Phi_{X_1 \to X'}^{{\bfG}}
\Phi_{X_1 \to X_1}^{{\bfE}_{1,2j+2}[2]}(A_1^{\vee}[2])
& i=0,
\\[4pt]
0 & i \not= -1,0.
\end{array}
\right.
\end{align*}
\end{itemize}
In particular we have exact sequences
\begin{multline*}
\tag{$\natural_{2j}$}
0 \to 
\Phi_{X_1 \to X'}^{{\bfG}}
\Phi_{X_0 \to X_1}^{{\bfE}_{0,2j+1}}(A_0)
\to 
\Phi_{X_1 \to X'}^{{\bfG}}
\Phi_{X_{1} \to X_1}^{{\bfE}_{1,2j}[2]}(B_{1})
\\
\to 
\Phi_{X_1 \to X'}^{{\bfG}}
\Phi_{X_1 \to X_1}^{{\bfE}_{1,2j}[2]}(A_{1}) \to
\Phi_{X_1 \to X'}^{{\bfG}}
\Phi_{X_0 \to X_1}^{{\bfE}_{0,2j+1}[2]}(B_{0}) \to 
0,
\end{multline*}
\begin{multline*}
\tag{$\natural_{2j+1}$}
0 \to 
\Phi_{X_1 \to X'}^{{\bfG}}
\Phi_{X_{1} \to X_1}^{{\bfE}_{1,2j+2}[2]}(B_{1}^{\vee})
\to 
\Phi_{X_1 \to X'}^{{\bfG}}
\Phi_{X_{0} \to X_1}^{{\bfE}_{0,2j+1}}(A_{0}^{\vee}[2])
\\
\phantom{0\to \Phi_{X_1 \to X'}^{{\bfG}}\Phi_{X_{0} \to X_1}^{{\bfE}_{0,2j+1}[2]}(B_{0}^{\vee})}
\to 
\Phi_{X_1 \to X'}^{{\bfG}}
\Phi_{X_{0} \to X_1}^{{\bfE}_{0,2j+1}[2]}(B_{0}^{\vee}) \to
\Phi_{X_1 \to X'}^{{\bfG}}
\Phi_{X_{1} \to X_1}^{{\bfE}_{1,2j+2}[2]}(A_{1}^{\vee}[2]) \to
0.
\end{multline*}
We note that 
$\chi(H^1(\Phi_{X_1 \to X'}^{{\bfG}}(F_{m}^\bullet))(n H'))$, $n \gg 0$
are decreasing for $0 \leq m \leq 2k-1$. 
Finally for $m=2k$, 
$H^1(\Phi_{X_1 \to X'}^{{\bfG}}(F_{m}^\bullet))(n H'))$ becomes zero.
Thus 
$\Phi_{X_1 \to X'}^{{\bfG}}(F_{2k}^\bullet)$ is a coherent sheaf
and we have an exact sequence
\begin{multline*}
\tag{$\natural_{2k}$}
0 \to \Phi_{X_1 \to X'}^{{\bfG}}\Phi_{X_0 \to X_1}^{{\bfE}_{0,2k+1}}(B_{0}) 
  \to \Phi_{X_1 \to X'}^{{\bfG}}\Phi_{X_0 \to X_1}^{{\bfE}_{0,2k+1}}(A_0)
\\
  \to \Phi_{X_1 \to X'}^{{\bfG}}\Phi_{X_{1} \to X_1}^{{\bfE}_{1,2k}[2]}(B_{1})
  \to \Phi_{X_1 \to X'}^{{\bfG}}\Phi_{X_1 \to X_1}^{{\bfE}_{1,2k}[2]}(A_{1}) 
  \to 0.
\end{multline*}
Hence the complex $\Phi_{X_1 \to X'}^{{\bfG}}(F_{m}^\bullet)$ eventually becomes a coherent sheaf as $m$ increases.
By $(\natural_m)$, $(\natural_{m+1})$ and \eqref{eq:relation-sharp},
we see that the operation 
$\Phi_{X_1 \to X'}^{{\bfG}}(F_m^\bullet) 
\mapsto 
\Phi_{X_1 \to X'}^{{\bfG}}(F_{m+1}^\bullet)
$
is nothing but the operation $\calF$ in Fact \ref{fact:yoshioka-lem2-2}.
For example, if $m=2j$, then by
$\bigl(
\Phi_{X_1 \to X'}^{{\bfG}}\Phi_{X_1 \to X_1}^{{\bfE}_{0,2j}[2]}\bigr)^{-1}$,
$(\natural_{2j})$ transforms to $(\#)$.
Then applying 
$\Phi_{X_1 \to X'}^{{\bfG}} \Phi_{X_1 \to X_1}^{{\bfE}_{0,2j}[2]} \calD_{X_1}$,
we get $(\natural_{2j+1})$.
\end{rmk}

\subsection{Construction and the application of the tame complex: case II}

Recall the notation in \S\S\,\ref{subsect:solution}. 
Assume $\ep=q^2-\ell p^2=-1$.
We will prove Proposition~\ref{prop:explicit} for this case.

In the present case, the algebraic integers $a_m,b_m$ 
in Definition~\ref{dfn:ambm} satisfy 
\begin{align*}
0<\dfrac{b_m}{a_m} -\sqrt{\ell}<\dfrac{b_{m-1}}{a_{m-1}} -\sqrt{\ell}
\end{align*} 
for $m\in\bbZ_{>0}$.
We also have the following sequence of inequalities.
\begin{multline*}
 -\infty=\dfrac{b_0}{a_0}
<-\dfrac{q}{p}=\dfrac{b_{-1}}{a_{-1}}
<\dfrac{b_{-2}}{a_{-2}}<\cdots<-\sqrt{\ell}<\cdots
<\dfrac{\ell a_{-2}}{b_{-2}}<\dfrac{\ell a_{-1}}{b_{-1}}=-\dfrac{\ell p}{q}
\\
<\dfrac{\ell a_0}{b_1}=0
<\dfrac{\ell a_1}{b_1}=\dfrac{\ell p}{q}
<\dfrac{\ell a_2}{b_2}
<\cdots<\sqrt{\ell}
<\cdots<\dfrac{b_2}{a_2}<\dfrac{b_1}{a_1}=\dfrac{q}{p}
<\dfrac{b_0}{a_0}=\infty.
\end{multline*}
In this case we have a FMT $\Phi_{X_0 \to X_1}^{{\bfE}}$ such that
\begin{align}
\label{eq:2-51}
\theta(\Phi_{X_0 \to X_1}^{{\bfE}})=
\begin{bmatrix} q & \ell p\\ p & q \end{bmatrix}.
\end{align}
Then Lemma~\ref{lem:matrix} yields 
\begin{align}
\label{eq:2-52}
\theta(\Phi_{X_1 \to X_0}^{\bfE})
=\begin{bmatrix} q & \ell p \\ p & q
 \end{bmatrix},
\quad
 \theta(\Phi_{X_1 \to X_0}^{\bfE^{\vee}[2]})
=\theta(\Phi_{X_0 \to X_1}^{\bfE^{\vee}[2]})
=\begin{bmatrix} q & -\ell p\\ -p & q \end{bmatrix}.
\end{align}

\begin{dfn}\label{defn:E_{i,m}-2}
We define ${\bfE}_{i,m} \in {\bf D}(X_i \times X_{i+m})$ ($m \geq 0$) 
inductively by
\begin{align*}
{\bfE}_{i,0}={\calO}_{\Delta},
\quad
\Phi_{X_{i+m-1} \to X_{i+m}}^{{\bfE}^{\vee}[2]}\,
\Phi_{X_i \to X_{i+m-1}}^{{\bfE}_{i,m-1}^{\vee}}
=
\Phi_{X_i \to X_{i+m}}^{{\bfE}_{i,m}^{\vee}}.
\end{align*} 
\end{dfn}

\begin{lem}
${\bfE}_{i,m}$ is a locally free sheaf on $X_i \times X_{i+m}$
for $m >0$.
\end{lem}

\begin{proof}
Assume that ${\bfE}_{i,m-1}$ is a locally free sheaf.
We note that 
$\mu({\bfE}|_{X_0 \times \{ x_1 \}}^{\vee})
=\mu({\bfE}|_{\{ x_0 \} \times X_1}^{\vee})
=-\tfrac{q}{p}\sqrt{n}$.
Since 
$\mu(\Phi_{X_i \to X_{i+m-1}}^{{\bfE}_{i,m-1}^{\vee}}(\bbC_{x_i}))
=-\tfrac{b_m}{a_m}\sqrt{n}$, 
we find that 
$\Phi_{X_{i+m-1} \to X_{i+m}}^{{\bfE}_{i,m-1}^{\vee}[2]}
 \Phi_{X_i \to X_{i+m-1}}^{{\bfE}^{\vee}}(\bbC_{x_i})$ is
a sheaf. Hence ${\bfE}_{i,m}$ is also a locally free sheaf.
\end{proof}

\begin{lem}\label{lem:E_{i,m}-2}
\begin{enumerate}
\item
For $m\in\bbZ_{>0}$, we have the following. 
\begin{align*}
&
\Phi_{X_i \to X_{i+m}}^{{\bfE}_{i,m}^{\vee}[2]}=
\Phi_{X_{i+m-1} \to X_{i+m}}^{{\bfE}^{\vee}[2]}
\Phi_{X_{i+m-2} \to X_{i+m-1}}^{{\bfE}^{\vee}[2]}
\dotsm
\Phi_{X_i \to X_{i+1}}^{{\bfE}^{\vee}[2]},
\\
&
\Phi_{X_i \to X_{i+m}}^{{\bfE}_{i+m,m}}=
\Phi_{X_{i+m-1} \to X_{i+m}}^{{\bfE}}
\Phi_{X_{i+m-2} \to X_{i+m-1}}^{{\bfE}}
\dotsm 
\Phi_{X_i \to X_{i+1}}^{{\bfE}},
\\
&\Phi_{X_i \to X_{i+m}}^{{\bfE}_{i,m}}=
\Phi_{X_{i+m-1} \to X_{i+m}}^{{\bfE}}
\Phi_{X_{i+m-2} \to X_{i+m-1}}^{{\bfE}}
\dotsm 
\Phi_{X_i \to X_{i+1}}^{{\bfE}}.
\end{align*}
Under the identification $\xi\colon X_i \times X_{i+m} \to
X_{i+m} \times X_i$, we also have
\begin{align*}
{\bfE}_{i+m,m}={\bfE}_{i,m}.
\end{align*}

\item
For $m\in\bbZ_{>0}$, we have the following.
\begin{align*}
\theta(\Phi_{X_i \to X_{i+m}}^{{\bf E}_{i,m}^{\vee}[2]})
=\begin{bmatrix} b_m    & -\ell a_m  \\   -a_m & b_m    \end{bmatrix}
=\begin{bmatrix} b_{-m} & \ell a_{-m}\\ a_{-m} & b_{-m} \end{bmatrix},
\quad
\theta(\Phi_{X_i \to X_{i+m}}^{{\bf E}_{i,m}})
=\begin{bmatrix} b_m & \ell a_m \\ a_m & b_m \end{bmatrix}.
\end{align*}
\end{enumerate}
\end{lem}

\begin{proof}
(1)
The first claim is obvious from Definition~\ref{defn:E_{i,m}-2}.
For the second and third one, the computation proceed as follows.
\begin{align*}
  \Phi_{X_i \to X_{i+m}}^{{\bfE}_{i+m,m}}
&=\bigl(\Phi_{X_{i+m} \to X_{i+2m}}^{{\bfE}_{i+m,m}^{\vee}[2]}\bigr)^{-1}
\\
&=\bigl(\Phi_{X_{i+2m-1} \to X_{i+2m  }}^{{\bfE}^{\vee}[2]}\,
   \Phi_{X_{i+2m-2} \to X_{i+2m-1}}^{{\bfE}^{\vee}[2]}\,
   \dotsm \,
   \Phi_{X_{i+m}    \to X_{i+m+1 }}^{{\bfE}^{\vee}[2]}\bigr)^{-1}
\\
&=\Phi_{X_{i+m-1}   \to X_{i+m   }}^{{\bfE}}\,
  \Phi_{X_{i+m-2}   \to X_{i+m-1 }}^{{\bfE}}\,
  \dotsm\,
  \Phi_{X_i \to X_{i+1   }}^{{\bfE}},
\\
  \Phi_{X_i \to X_{i+m   }}^{{\bfE}_{i,m}}
&=\calD_{X_{i+m}}\,
  \Phi_{X_i \to X_{i+m   }}^{{\bfE}_{i,m}^{\vee}[2]}\,
  \calD_{X_i}
\\
&=\calD_{X_{i+m}}\,
  \Phi_{X_{i+m-1} \to X_{i+m  }}^{{\bfE}^{\vee}[2]}\,
  \Phi_{X_{i+m-2} \to X_{i+m-1}}^{{\bfE}^{\vee}[2]}\,
  \dotsm\,
  \Phi_{X_i       \to X_{i+1  }}^{{\bfE}^{\vee}[2]}\,
  \calD_{X_i}
\\
&= 
  \Phi_{X_{i+m-1} \to X_{i+m  }}^{{\bfE}}\,
  \Phi_{X_{i+m-2} \to X_{i+m-1}}^{{\bfE}}\,
  \dotsm \,
  \Phi_{X_i       \to X_{i+1  }}^{{\bfE}}.
\end{align*}
The last claim is a consequence of the rest ones.

(2)
This is a consequence of (1), \eqref{eq:2-51} and \eqref{eq:2-52}.
\end{proof}

In the present case $\ep=q^2-\ell p^2=1$, 
$F \in M_X^H(1,0,-\ell)$ has the cokernel presentation
\begin{align*}
0 \to \Phi_{X_0 \to X_1}^{{\bfE}}(A_0)^{\vee} 
  \to \Phi_{X_0 \to X_1}^{{\bfE}}(B_0)^{\vee} \to F \to 0.
\end{align*}
Thus we have an exact sequence
\begin{align}\label{eq:a}\tag{$\flat$}
0 \to \Phi_{X_0 \to X_1}^{{\bfE}}(A_0)^{\vee} 
  \to \Phi_{X_0 \to X_1}^{{\bfE}}(B_0)^{\vee} 
  \to B_1 \to A_1 \to 0.
\end{align}
Applying $\Phi_{X_1 \to X_0}^{{\bfE}^{\vee}}$ to the dual of
\eqref{eq:a}, we get
\begin{align*}
0 \to \Phi_{X_1 \to X_0}^{{\bfE}^{\vee}[2]}(A_1^{\vee}) 
  \to \Phi_{X_1 \to X_0}^{{\bfE}^{\vee}[2]}(B_1^{\vee}) 
  \to B_0 \to A_0 \to 0.
\end{align*}

Now we can define the tame complex for the case $\ep=1$.

\begin{dfn}\label{defn:F_m-2}
For $F\seteq\Ker(B_1 \to A_1)$, 
we set the \emph{tame complex} $F_m^\bullet$ by
\begin{align*}
F_m^\bullet\seteq
\left\{
\begin{array}{l l l l}
\bigl(\Phi_{X_1 \to X_1}^{{\bfE}_{1,2}^{\vee}[2]}\bigr)^{m/2}(F)
&=\Phi_{X_1 \to X_1}^{{\bfE}_{1,2k}^{\vee}[2]}(F)
&m=2k,&k\ge0,
\\[4pt]
\bigl(\Phi_{X_1 \to X_1}^{{\bfE}_{1,2}^{\vee}[2]}\bigr)^{(m+1)/2}(F^{\vee})
&=\Phi_{X_1 \to X_1}^{{\bfE}_{1,2k}^{\vee}[2]}(F^{\vee})
&m=2k-1,&k\ge1,
\\[4pt]
\bigl(\Phi_{X_1 \to X_1}^{{\bfE}_{1,2}^{\vee}[2]}\bigr)^{(m+1)/2}(F^{\vee})
&=\Phi_{X_1 \to X_1}^{{\bfE}_{1,2k}}(F^{\vee})
&m=-2k-1,&k\ge0,
\\[4pt]
\bigl(\Phi_{X_1 \to X_1}^{{\bfE}_{1,2}^{\vee}[2]}\bigr)^{m/2}(F)
&=\Phi_{X_1 \to X_1}^{{\bfE}_{1,2k}}(F)
&m=-2k,&k\ge1.
\end{array}
\right.
\end{align*}
\end{dfn}

We can rewrite $F_m^\bullet$ as follows.

\begin{lem}
\label{lem:another-2}
\begin{align*}
F_m^\bullet=
\left\{
\begin{array}{l l}
\Phi_{X_0 \to X_1}^{{\bfE}_{0,2k+1}^{\vee}[2]}
 \bigl([A_0^{\vee} \to B_0^{\vee}]\bigr)
&m=2k,
\\[4pt]
\Phi_{X_0 \to X_1}^{{\bfE}_{0,2k-1}^{\vee}[2]}
 \bigl([A_0^{\vee} \to B_0^{\vee}]^{\vee}\bigr)
&m=2k-1,
\\[4pt]
\Phi_{X_0 \to X_1}^{{\bfE}_{1,2k+1}}
 \bigl([A_0^{\vee} \to B_0^{\vee}]^{\vee}\bigr)
&m=-2k-1,
\\[4pt]
\Phi_{X_0 \to X_1}^{{\bfE}_{1,2k-1}}
 \bigl([A_0^{\vee} \to B_0^{\vee}]\bigr)
&m=-2k.
\end{array}
\right.
\end{align*}
\end{lem}
\begin{proof}
This is the consequence of the definition and Lemma \ref{lem:E*E-2} below.
\end{proof}

\begin{lem}\label{lem:E*E-2}
\begin{align*}
\begin{array} {l l l}
\Phi_{X_1 \to X_{1+m}}^{{\bfE}_{1,m}^{\vee}[2]} \calD_{X_1}
\Phi_{X_0 \to X_1}^{{\bfE}}
&=&\Phi_{X_0 \to X_{m+1}}^{{\bfE}_{0,m+1}^{\vee}[2]} \calD_{X_0},
\\[4pt]
\bigl(\Phi_{X_1 \to X_{1+m}}^{{\bfE}_{1,m}^{\vee}[2]} \calD_{X_1}\bigr)
\bigl(\calD_{X_1}\Phi_{X_0 \to X_1}^{{\bfE}}\bigr)
&=&\Phi_{X_0 \to X_{m-1}}^{{\bfE}_{0,m-1}^{\vee}[2]},
\\[4pt]
\Phi_{X_1 \to X_{1+m}}^{{\bfE}_{1,m}}\calD_{X_1}
\Phi_{X_0 \to X_1}^{{\bfE}}
&=&\Phi_{X_0 \to X_{m-1}}^{{\bfE}_{0,m-1}}\calD_{X_0},
\\[4pt]
\bigl(\Phi_{X_1 \to X_{1+m}}^{{\bfE}_{1,m}}\calD_{X_1}\bigr)
\bigl(\calD_{X_1}\Phi_{X_0 \to X_1}^{{\bfE}}\bigr)
&=&\Phi_{X_0 \to X_{1+m}}^{{\bfE}_{0,1+m}}.
\end{array}
\end{align*}
\end{lem}

\begin{proof}
The second and the last equality follows from the definitions of $E_{i,m}$. For the remained equations, we have 
\begin{align*}
&\Phi_{X_1 \to X_{1+m}}^{{\bfE}_{1,m}^{\vee}[2]} \calD_{X_1}
\Phi_{X_0 \to X_1}^{{\bfE}}=
\Phi_{X_1 \to X_{1+m}}^{{\bfE}_{1,m}^{\vee}[2]}
\Phi_{X_0 \to X_1}^{{\bfE}^{\vee}[2]} \calD_{X_0}
=\Phi_{X_0 \to X_{m+1}}^{{\bfE}_{0,m+1}^{\vee}[2]} \calD_{X_0},
\\[4pt]
&\Phi_{X_1 \to X_{1+m}}^{{\bfE}_{1,m}} \calD_{X_1}
\Phi_{X_0 \to X_1}^{{\bfE}}=
\Phi_{X_1 \to X_{1+m}}^{{\bfE}_{1,m}}
\Phi_{X_0 \to X_1}^{{\bfE}^{\vee}[2]} \calD_{X_0}
=\Phi_{X_0 \to X_{m-1}}^{{\bfE}_{0,m-1}} \calD_{X_0}.
\end{align*}
\end{proof}

For the proof of Proposition~\ref{prop:explicit}, we prepare some calculations.

\begin{lem}\label{lem:A-B-2}
\begin{align*}
\Phi_{X_i \to X_{i+m}}^{\bfE_{i,m}^{\vee}[2]}(B_i),\,
\Phi_{X_i \to X_{i+m}}^{\bfE_{i,m}^{\vee}[2]}(B_i^{\vee}),\,
\Phi_{X_i \to X_{i+m}}^{\bfE_{i,m}}(B_i),\,
\Phi_{X_i \to X_{i+m}}^{\bfE_{i,m}}(B_i^{\vee})\,
\in \Coh(X_{i+m}).
\end{align*}
\end{lem}

\begin{proof}
Since $\mu({\bfE}_{i,m}|_{\{ x_i \} \times X_{i+m}})=
\mu({\bfE}_{i,m}|_{X_i \times \{ x_{i+m} \}})=
b_m/a_m>0$ for $m \geq 1$ and $\deg(B_i)=0$, we get our claim.
\end{proof}

\begin{proof}[{Proof of Proposition~\ref{prop:explicit}}]

By Lemma \ref{lem:A-B-2}, we have the following exact sequences
for $k \geq 0$.
\begin{align*}
\tag{$\flat_{2k+2}$}
&
0 \to \Phi_{X_1 \to X_1}^{{\bfE}_{1,2k+2}^{\vee}}(A_1)
  \to \Phi_{X_0 \to X_1}^{{\bfE}_{0,2k+3}^{\vee}}(A_0^{\vee}[2])
  \to \Phi_{X_0 \to X_1}^{{\bfE}_{0,2k+3}^{\vee}[2]}(B_0^{\vee})
  \to \Phi_{X_1 \to X_1}^{{\bfE}_{1,2k+2}^{\vee}[2]}(B_1)
  \to 0,
\\[4pt]
\tag{$\flat_{2k+1}$}
&
0 \to \Phi_{X_0 \to X_1}^{{\bfE}_{0,2k+1}^{\vee}}(A_0)
  \to \Phi_{X_1 \to X_1}^{{\bfE}_{1,2k+2}^{\vee}}(A_1^{\vee}[2])
  \to \Phi_{X_1 \to X_1}^{{\bfE}_{1,2k+2}^{\vee}[2]}(B_1^{\vee})
  \to \Phi_{X_0 \to X_1}^{{\bfE}_{0,2k+1}^{\vee}[2]}(B_0)
  \to 0,
\\[4pt]
\tag{$\flat_{-2k-1}$}
&
0 \to \Phi_{X_1 \to X_1}^{{\bfE}_{1,2k}}(B_1^{\vee})
  \to \Phi_{X_0 \to X_1}^{{\bfE}_{0,2k+1}}(B_0)
  \to \Phi_{X_0 \to X_1}^{{\bfE}_{0,2k+1}}(A_0)
  \to \Phi_{X_1 \to X_1}^{{\bfE}_{1,2k}}(A_1^{\vee}[2])
  \to 0,
\\[4pt]
\tag{$\flat_{-2k-2}$}
&
0 \to \Phi_{X_0 \to X_1}^{{\bfE}_{0,2k+1}}(B_0^{\vee})
  \to \Phi_{X_1 \to X_1}^{{\bfE}_{1,2k+2}}(B_1)
  \to \Phi_{X_1 \to X_1}^{{\bfE}_{1,2k+2}}(A_1)
  \to \Phi_{X_0 \to X_1}^{{\bfE}_{0,2k+1}}(A_0^{\vee}[2])
  \to 0.
\end{align*}

These exact sequences and Lemma~\ref{lem:another-2} 
give the following relations of FMTs.
\begin{align*}
\begin{array}{l l l}
F_{2k}^\bullet&=& 
 \Phi_{X_1 \to X_1}^{{\bfE}_{1,2k}^{\vee}}\,
 \calD_{X_1}\,
 \bigl(\Phi_{X_1 \to X_1}^{{\bfE}_{1,2k}^{\vee}}\bigr)^{-1}
  (F_{2k-1}^\bullet),
\\[4pt]
F_{2k+1}^\bullet&=& 
 \Phi_{X_0 \to X_1}^{{\bfE}_{0,2k+1}^{\vee}}\,
 \calD_{X_0}\,
 \bigl(\Phi_{X_0 \to X_1}^{{\bfE}_{0,2k+1}^{\vee}}\bigr)^{-1}
  (F_{-2k-1}^\bullet),
\\[4pt]
F_{-2k-1}^\bullet&=& 
 \Phi_{X_1 \to X_1}^{{\bfE}_{1,2k}}\,
 \calD_{X_1}\,
 \bigl(\Phi_{X_1 \to X_1}^{{\bfE}_{1,2k}}\bigr)^{-1}(F_{-2k}^\bullet),
\\[4pt]
F_{-2k-2}^\bullet&=&
 \Phi_{X_0 \to X_1}^{{\bfE}_{0,2k+1}}\,
 \calD_{X_0}\,
 \bigl(\Phi_{X_0 \to X_1}^{{\bfE}_{0,2k+1}}\bigr)^{-1}(F_{-2k-1}^\bullet).
\end{array}
\end{align*}

{}From these relation, we can easily get the claim.
\end{proof}

\begin{dfn}
We set
\begin{align*}
\begin{array}{l l l l l l}
I_0     &\seteq&(0,\tfrac{\ell a_1}{b_1}] \cup (\tfrac{b_1}{a_1},\infty],
&
I_{-1}  &\seteq&(-\infty,-\tfrac{b_1}{a_1}] \cup (-\tfrac{\ell a_1}{b_1},0],
\\[4pt]
I_{m}   &\seteq&(\tfrac{\ell a_m}{b_m},\tfrac{\ell a_{m+1}}{b_{m+1}}]
                \cup
                (\tfrac{b_{m+1}}{a_{m+1}},\tfrac{b_{m} }{a_{m} }],
&
I_{-m-1}&\seteq&(-\tfrac{b_m}{a_m},-\tfrac{b_{m+1}}{a_{m+1}}]
                \cup
               (-\tfrac{\ell a_{m+1}}{b_{m+1}},-\tfrac{\ell a_{m}}{b_{m}}]\quad
                m \geq 1.
\end{array}
\end{align*}
\end{dfn}
Then we have
${\bbP}^1({\bbR}) \setminus \{\pm\sqrt{\ell}\}=
\coprod_{m \in \bbZ} I_m =\coprod_{m \in \bbZ} I_m^*$.

\begin{thm}\label{thm:csv:2}
For a FMT $\Phi_{X_1 \to X'}^{{\bfG}}\colon{\bf D}(X_1) \to {\bf D}(X')$,
we set $\lambda \seteq \mu({\bfG}_{x'})/\sqrt{n}$.

(1)
If $\lambda \in I_m$, then
$\Phi_{X_1 \to X'}^{\bfG}(F_m^\bullet)$ is a coherent sheaf.

(2)
If $\lambda \in I_m^*$, then
$\calD_{X'}\, \Phi_{X_1 \to X'}^{\bfG}(F_m^\bullet)=
\Phi_{X_1 \to X'}^{{\bfG}^{\vee}[2]}\, \calD_{X_1}(F_m^\bullet)$ 
is a coherent sheaf.
\end{thm}

The proof is based on the calculation in this subsection and is similar to that of Theorem \ref{thm:csv:1}. We omit the detail.

\begin{rmk}
In \cite[Theorem 3]{Mukai:1979} Mukai proved the case $\ell=2$, $n=1$ of our Theorem \ref{thm:csv:1}. He also claims that the operation $\calF$  in Remark \ref{rmk:calF} gives an isomorphism between the moduli and the Hilbert scheme. We will return to this point in our future work.
\end{rmk}

\section{Birational morphism of the moduli}
\label{sect:birational}

In this section we give an explicit description of the birational map between the moduli of stable sheaves and the Hilbert scheme of points: $M_X^H(v)\cdots\to X\times\Hilb{\mpr{v^2}/2}{X}$. Fix an abelian surface $X$ with $\NS(X)=\bbZ H$ and set $n\seteq(H^2)/2$.

\subsection{Quadratic forms and the description of the birational map}
\label{subsec:quad}

Recall the matrix description introduced in \S\S\,\ref{subsec:matrix}.
We begin with the next simple remark.

\begin{lem}
\label{lem:quad}
For an isotropic Mukai vector $w=(p^2,\tfrac{p q}{\sqrt{n}}H,q^2)$ and a Mukai vector $v=(r,d H,a)$,
\begin{align*}
\langle w,v \rangle=Q_v(q,-p).
\end{align*}
Here the symbol $Q_v$ is the quadratic form defined to be
\begin{align*}
Q_v(x,y)
\seteq-(r x^2+2 d\sqrt{n}x y+a y^2)
=-
\begin{bmatrix}x &y\end{bmatrix}
\begin{bmatrix}r & d\sqrt{n}\\d\sqrt{n} & a\end{bmatrix}
\begin{bmatrix}x \\y\end{bmatrix}.
\end{align*} 
\end{lem}

\begin{proof}
\begin{align*}
\langle w,v \rangle 
&=\phantom{m}
B \left(
\begin{bmatrix}p^2 & p q \\ p q & q^2 \end{bmatrix},
\begin{bmatrix}r & d\sqrt{n}\\ d\sqrt{n} & a \end{bmatrix} 
\right)
\phantom{mn}
=-\tr \left(
\begin{bmatrix}q^2 & -p q \\ -p q & p^2 \end{bmatrix}
\begin{bmatrix}r & d\sqrt{n}\\ d\sqrt{n} & a \end{bmatrix} 
\right)\\
&=-\tr \left(
\begin{bmatrix}q \\ -p \end{bmatrix}
\begin{bmatrix}q &-p \end{bmatrix}
\begin{bmatrix}r & d\sqrt{n}\\ d\sqrt{n} & a \end{bmatrix} 
\right)
= 
-\tr \left(
\begin{bmatrix}q &-p\end{bmatrix}
\begin{bmatrix}r & d\sqrt{n}\\d\sqrt{n} & a\end{bmatrix}
\begin{bmatrix}q \\-p\end{bmatrix}
\right)
\\
&=Q_v(q,-p).
\end{align*}
\end{proof}

\begin{prop}\label{prop:biraional}
Fix a positive Mukai vector $v=(r,d H,a)$ and assume $\ell\seteq \mpr{v^2}/2$ is positive. Assume further that there exists a solution $(p_1,q_1)$ of the quadratic indefinite equation
\begin{align}
\label{eq:quad}
Q_v(q_1,-p_1)=-a p_1^2+2d\sqrt{n}p_1 q_1-r q_1^2=\ep,\quad \ep=\pm1,
\end{align}
with the additional condition  
\begin{align}
\label{eq:p_1,q_1}
p_1^2,\ p_1q_1/\sqrt{n},\ q_1^2\ \in \bbZ.
\end{align}
We also set
\begin{align*}
p_2\seteq \ep(d \sqrt{n} p_1-r q_1),\quad
q_2\seteq \ep(-d\sqrt{n} q_1+a p_1).
\end{align*}
\begin{enumerate}
\item
Take $(p_1,q_1)$ satisfying \eqref{eq:p_1,q_1}
such that $|p_1|$ is minimum among them. 
Then a general member of $M_X^H(v)$ has the next semi-homogeneous presentation.
\begin{align*}
\begin{array}{l l l}
\text{(a)}
&\Ker(E_1\to E_2)
&\text{ if } - p_2/p_1 \ge 0\ \text{ and }\ \ep=+1, 
\\
\text{(b)}
&\Ker(E_2\to E_1)
&\text{ if } - p_2/p_1  <  0\ \text{ and }\ \ep=-1, 
\\
\text{(c)}
&\Coker(E_1\to E_2)
&\text{ if } - p_2/p_1 \ge 0\ \text{ and }\ \ep=-1, 
\\
\text{(d)}
&\Coker(E_2\to E_1)
&\text{ if } - p_2/p_1 <   0\ \text{ and }\ \ep=+1. 
\end{array}
\end{align*}
Here $E_i$ are semi-homogeneous sheaves with $v(E_1)=\ell (p_1^2, \tfrac{p_1 q_1}{\sqrt{n}} H, q_1^2)$ and $v(E_2)=(p_2^2, \tfrac{p_2 q_2}{\sqrt{n}} H, q_2^2)$.
\item
The element
\begin{align*}
g \seteq
\begin{bmatrix}
q_1 & -q_2 \\-p_1 & p_2
\end{bmatrix}
\end{align*}
belongs to $G$, and we can take a universal family $\bfF$ on $X\times X'$ such that $\theta(\Phi_{X \to X'}^{{\bfF}})=\pm g\in\PSL(2,\bbR)$. Then the following functor induces the birational map  $M_X^H(r,d H,a) \cdots \to M_{X'}^{H'}(1,0,-\ell)$ for each case (a)-(d).
\begin{align*}
\text{(a)}\; \calD_{X'}\Phi_{X \to X'}^{{\bfF}[1]},\quad
\text{(b)}\; \Phi_{X \to X'}^{{\bfF}[-2]},\quad
\text{(c)}\; \calD_{X'}\Phi_{X \to X'}^{{\bfF}},\quad
\text{(d)}\; \Phi_{X \to X'}^{{\bfF}[-1]}.
\end{align*}
By this functor, the semi-homogeneous sheaves $E_1$, $E_2$ appearing in (1) are transformed in to $L\otimes \calO_Z$ and $L$, where $Z$ is a 0-dimensional subscheme of $X$ with length $\ell$ and $L$ is a line bundle on $X$.
\end{enumerate}
\end{prop}

\begin{proof}
The condition \eqref{eq:quad} and \eqref{eq:p_1,q_1} immediately show $g\in G$, so that $\Phi_{X\to X'}^\bfF$ does exist. Since the rest statements of (2) follow from those of (1), we give the proof of (1).

Lemma~\ref{lem:quad} shows that the equation \eqref{eq:quad} is equivalent to the numerical equation for $v$. By Theorem~\ref{thm:mukai-conj1'}, we know that a general member of $M_X^H(v)$ has a semi-homogeneous presentation which associates to the numerical solution of minimal rank. Then Proposition~\ref{prop:sp} means that the statements follows from the WIT with respect to $\Phi_{X\to X'}^\bfF$. We study the WIT with respect to the inverse transform $\Phi_{X'\to X}^{\bfF^{\vee}}$, since it is more convenient.

A simple calculation shows that
\begin{align*}
&
{}^t g^{-1} \begin{bmatrix}0 & 0 \\ 0 & 1\end{bmatrix} g^{-1}
=\begin{bmatrix}p_1^2 & p_1 q_1\\ p_1 q_1 & q_1^2\end{bmatrix},
\quad
{}^t g^{-1} \begin{bmatrix}1 & 0 \\ 0 & 0\end{bmatrix} g^{-1}
=\begin{bmatrix}p_2^2 & p_2 q_2\\ p_2 q_2 & q_2^2\end{bmatrix},
\\
&
{}^t g^{-1}\begin{bmatrix}1 & 0 \\ 0 & -\ell\end{bmatrix} g^{-1}
=-\epsilon \begin{bmatrix}r & d\sqrt{n}\\d\sqrt{n} & a\end{bmatrix}.
\end{align*}
Then Proposition~\ref{prop:matrix} and the equality $\theta(\Phi_{X' \to X}^{\bfF^{\vee}}) = \pm g^{-1}=\pm\begin{bmatrix}p_2 &q_2 \\ p_1 &q_1\end{bmatrix}$ yield
\begin{align}
\nonumber
&v(\Phi_{X' \to X}^{\bfF^{\vee}}(\bbC_{x'}))
=(p_1^2, \tfrac{p_1 q_1}{\sqrt{n}} H,q_1^2),
\quad
v(\Phi_{X' \to X}^{\bfF^{\vee}}(\calO_{X'}))
=(p_2^2, \tfrac{p_2 q_2}{\sqrt{n}} H,q_2^2),
\\
\nonumber
&v(\Phi_{X' \to X}^{\bfF^{\vee}}(I_Z \otimes L))=v(E_2)-v(E_1),
\end{align}
where we used symbols $x'\in X'$ and $L \in \Piczero{X'}$ 
and the symbol $Z$ means a 0-dimensional subscheme of length $\ell$. 
We also find that $\WIT_0$ holds for $\calO_{X'}$ 
with respect to $\Phi_{X' \to X}^{\bfF^\vee}$ if and only if $-p_2/p_1<0$, 
since $\mu({\bfF}|_{\{x \} \times X'})=-p_2/p_1$.

Hence if $\mu(\bfF|_{\{x\} \times X'})<0$, 
then both $\Phi_{X' \to X}^{\bfF^{\vee}}(\bbC_{x'})$ 
and $\Phi_{X' \to X}^{\bfF^{\vee}}(\calO_{X'})$ are semi-homogeneous sheaves. 
By the minimality of $|p_1|$, 
$F^\bullet\seteq \Phi_{X' \to X}^{\bfF^{\vee}}(I_Z \otimes L)$ 
is a coherent sheaf up to shift for a general $I_Z \otimes L$. 
The shift is determined by the sign of $\rk(F^\bullet)=\ell p_1^2-p_2^2$, 
which equals to $\ep r$ by the definition of $p_2$ and \eqref{eq:quad}. 
Thus we have the statements for the cases (b) and (d).

If $\mu(\bfF|_{\{x\} \times X'})\ge 0$, 
then both $\Phi_{X' \to X}^{\bfF^{\vee}[2]}\,\calD_{X'}(\bbC_{x'})$ 
and $\Phi_{X' \to X}^{\bfF^{\vee}[2]}\,\calD_{X'}(\calO_{X'})$ 
are semi-homogeneous sheaves. 
Similar arguments yield the statements for (a) and (c).
\end{proof}

\subsection{Structures of $G$}

As the preliminary of the next subsection, 
we study the structure of the arithmetic group $G$ 
introduced in Definition~\ref{dfn:G}.

\begin{dfn}
Let $n=p_1^{e_1}p_2^{e_2}\dotsm p_N^{e_N}$ be the prime decomposition 
of the integer $n \seteq (H^2)/2$.
We define a map $\widetilde{\phi}_i \colon \bbZ \to \{ 0,1 \}$ sending $m$ to
\begin{align*}
\widetilde{\phi}_i(m)=
\begin{cases}
1\; & p_i \,|\, m, \\
0\; & p_i \not |\, m.
\end{cases} 
\end{align*}
We identify $\{ 0,1 \}$ with $\bbZ/2 \bbZ$ in a natural way,
and regard it as a group. Then we also define a map:
\begin{align*}
\begin{array}{c c c c c}
\widetilde{\phi}\colon 
 & \bbZ & \to     & (\bbZ/2 \bbZ)^{\oplus N}\\
 & m    & \mapsto 
 & (\widetilde{\phi}_1(m),\widetilde{\phi}_2(m),\ldots,\widetilde{\phi}_N(m)).
\end{array}
\end{align*}
Clearly this $\widetilde{\phi}$ is a group homomorphism.
Then for $g=
\begin{bmatrix} 
a\sqrt{r} & b\sqrt{s} \\ c\sqrt{s} & d\sqrt{r}
\end{bmatrix}
\in G$, 
we set $\phi(g) \seteq \widetilde{\phi}(r)$.
By Lemma \ref{lem:group}, $r$ is uniquely determined by $g$.  
Thus $\phi(g)$ is well-defined.
\end{dfn}

\begin{lem}
$\phi(g g')=\phi(g)+\phi(g')$.
Thus $\phi:G \to (\bbZ/2\bbZ)^{\oplus N}$
is a homomorphism. 
\end{lem}

\begin{proof}
We use the notation in \eqref{eq:multiplication}.
We set $r=p r_1$, $r'=p r_1'$, $\gcd(r_1,r_1')=1$
and $s=q s_1$, $s'=q s_1'$, $\gcd(s_1,s_1')=1$.
Then $n=p q r_1 s_1=p q r_1' s_1'$.
Hence $r_1'=s_1$ and $s_1'=r_1$.
Thus 
\begin{align*}
\begin{bmatrix}
a\sqrt{r} & b\sqrt{s} \\
c\sqrt{s} & d\sqrt{r}
\end{bmatrix}
\begin{bmatrix}
a'\sqrt{r'} & b'\sqrt{s'} \\
c'\sqrt{s'} & d'\sqrt{r'}
\end{bmatrix}
=
\begin{bmatrix}
(a a' p+b c' q)\sqrt{r_1 s_1} & (a b' r_1+b d' s_1)\sqrt{p q} \\
(c a' s_1+d c' r_1)\sqrt{p q} & (c b' q+d d' p)\sqrt{r_1 s_1}
\end{bmatrix}.
\end{align*}
By the uniqueness of the description (Lemma \ref{lem:group}), 
$\phi(g g')=\widetilde{\phi}(r_1 s_1)=\widetilde{\phi}(r_1 r_1')
=\widetilde{\phi}(p r_1)+\widetilde{\phi}(p r_1')
=\phi(g)+\phi(g')$.
\end{proof}

\begin{prop}
\begin{enumerate}
\item
$\phi$ is surjective and $\ker \phi$ is isomorphic 
to the modular group $\Gamma_0(n)$.
%
\item
For $g=
\begin{bmatrix}
a\sqrt{r} & b\sqrt{s} \\
c\sqrt{s} & d\sqrt{r}
\end{bmatrix} \in G$ and $Y \in \FM(X)$, we set
$Y_{g^{-1}} \seteq M_Y^{\widehat{H}}(c^2 s, c d \widehat{H},d^2 r)$,
where $\widehat{H}$ is the ample generator of $\NS(Y)$.
Then $(Y_{g^{-1}})_{{g'}^{-1}} \cong Y_{(g' g)^{-1}}$ and 
$Y_{g^{-1}}$ is determined by $\phi(g)$.
\item
Assume that $\rk \End(X)=1$. Then
$Y_{g^{-1}} \cong Y$ if and only if $g \in \ker \phi$.
In particular, $\ker \phi=\theta(\mathrm{Eq}_0({\bf D}(Y),{\bf D}(Y)))$.
\end{enumerate}
\end{prop}

\begin{proof}
(1)
The surjectivity is easy. 
We note that
\begin{align*}
\ker \phi=
\left\{
\left.
\begin{bmatrix}
a & b\sqrt{n} \\
c\sqrt{n} & d
\end{bmatrix} 
\, \right|\; a,b,c,d \in \bbZ,\;
a d-b c n=1
\right\}.
\end{align*}
It is a conjugate of $\Gamma_0(n)$. 
Thus the second claim also holds. 

(2)
We have a FMT $\Phi_g:{\bfD}(Y_{g^{-1}}) \to {\bfD}(Y)$ such that
$\theta(\Phi_g)=g$.
Since the composition $\Phi_g \Phi_{g'}:{\bfD}(
(Y_{g^{-1}})_{{g'}^{-1}}) \to
{\bf D}(Y_{g^{-1}}) \to {\bfD}(Y)$ satisfies
$\theta(\Phi_g \Phi_{g'})=g' g=\theta(\Phi_{g' g})$, 
we have $(Y_{g^{-1}})_{{g'}^{-1}} \cong Y_{(g' g)^{-1}}$.
The second claim follows from the following Lemma~\ref{lem:FMpartner}.

(3) is due to Orlov \cite[Example~4.16]{Orlov:2002}.
See also the argument in the next Lemma~\ref{lem:FMpartner}.
\end{proof}

\begin{lem}
\label{lem:FMpartner}
Let $p,q\in\bbZ$ and $r,s\in\bbZ_{>0}$ with $r s=n\seteq(H^2)/2$.
\begin{enumerate}
\item
If $\gcd(p r,q s)=1$, then we have 
$M_X^H(p^2 r,p q H,q^2 s) \cong M_X^H(r,H,s)$.
\item
Moreover if $\rk \End(X)=1$, then
$M_X^H(r,H,s) \cong X$ if and only if $s=1$.
\end{enumerate}
\end{lem}

\begin{proof}
(1)
First assume that $p q\neq 0$.
Let $E$ be a semi-homogeneous vector bundle
with $v(E)=(p^2 r,p q H,q^2 s)$.
Then we have 
$M_X^H(p^2 r,p q H,q^2 s) \cong X/K(E)$, where 
$K(E)\seteq \{ x \in X \mid T_x^*(E) \cong E \}$.
By \cite[Cor.~7.8]{Mukai:1978},
$K(E)=p^2 r K(p q H)$,
where
$K(H) \seteq \{x \in X \mid T_x^*({\calO}_X(H)) \cong {\calO}_X(H) \}$.
We write $X=\bbC^2/\Lambda$, where
$\Lambda$ is a lattice of rank 4.
Then there is a decomposition $\Lambda=L_1 \oplus L_2$
such that $L_i \cong \bbZ^{\oplus 2}$  ($i=1,2$) and
$K(H)=\frac{1}{n}L_1/L_1=\frac{1}{r s}L_1/L_1$.
Since
$K(p q H)=\frac{1}{p q r s}L_1/L_1
\oplus \frac{1}{p q}L_2/L_2$,
we have
$p^2 r K(p q H)=p(\frac{1}{q s}L_1/L_1)
\oplus p r(\frac{1}{q}L_2/L_2)$.
Since $\gcd(p,q s)=1$ and $\gcd(p r,q)=1$, we have
$p(\frac{1}{q s}L_1/L_1)
\oplus p r(\frac{1}{q}L_2/L_2)=\frac{1}{q s}L_1/L_1
\oplus \frac{1}{q}L_2/L_2$.
Thus $p^2 r K(p q H)$ contains the subgroup $X_q$ of $q$-torsion
points of $X$.  
Since the kernel of $q 1_X:X \to X$ is $X_q$, 
$q 1_X$ induces an isomorphism
$p^2 r K(p q H)/X_q \cong r K(H)$.
Hence we have isomorphisms
$M_X^H(p^2 r,p q H,q^2 s) \cong X/(p^2 r K(p q H)) \cong
X/(r K(H)) \cong M_X^H(r,H,s)$.

Next we suppose $p q=0$. 
Then by the assumption we have $(p, q, r, s)=(1,0,1,n)$ or $(0,1,n,1)$.
For the first case, we have $M_X^H(1,0,0)=\widehat{X}\cong M_X^H(1,H,n)$.
For the second case, $M_X^H(0,0,1)=X\cong X/(n K(H))\cong M_X^H(n,H,1)$.
Thus the conclusion holds.

(2)
Assume that $\rk \mathrm{End}(X)=1$.
If $X/(r K(H)) \cong X$, then
we have a homomorphism $\psi:X \to X/(r K(H)) \to X$.
By our assumption,
$m_1 \psi = m_2 1_X$ for $m_1\in \bbZ_{>0}$ and $m_2 \in \bbZ$.
Then $r K(H)=\ker(\psi)=m_1 X_{m_2}$.
If $s\ne 1$, this is impossible.
Hence $X/(r K(H)) \not \cong X$ if $s\ne 1$.
If $s=1$, then clearly we have $X/(r K(H)) \cong X$.
\end{proof}

We conclude this subsection with the next lemma,
although it will not be used afterwards.

\begin{lem}
Assume that $r,a\in\bbZ_{>0}$ satisfies $r a=n\seteq(H^2)/2$.
Then the dual variety of $M_X^H(r,H,a)$ is isomorphic to $M_X^H(a,H,r)$.
\end{lem}

\begin{proof}
Set $Y\seteq M_X^H(r,H,a)$ and $Z\seteq M_X^H(a,H,r)$. 
We construct an isomorphism $Z \to \widehat{Y}$ using FMT. 
Since $\gcd(r,a)=1$, 
there exists a pair $(p,q)\in\bbZ^2$ such that $p r-q a=1$. 
Consider the moduli $\widetilde{Y}\seteq M_X^H(p^2 r,p q H,q^2 a)$. 
Since $\mpr{(p^2r,p q H,q^2 a),(a,H,r)}=-(p r-q a)^2=-1$, 
there is a universal family $\bfE$ on $Z\times \widetilde{Y}$. 
We have $Y\cong\widetilde{Y}$ by (1). 
Thus there exits a universal family $\bfF$ on $Z\times Y$.  
Then choosing an appropriate integer $k$, 
we have a FMT $\Phi_{Z\to Y}^{\bfF[k]}$ 
which sends $\bbC_z$ ($z\in Z$) to a line bundle on $Y$. 
Retaking $\bfF$ as in Proposition~\ref{prop:semihom-fmt}~(4), 
we have $Z\cong M_Y^{\widehat{H}}(1,0,0)$.
\end{proof}

\subsection{Principally polarized case}
\label{subsec:pp} 

In the case $n=1$, namely when $X$ is principally polarized, 
we have a simple description of the birational map.

First we note the next proposition.

\begin{prop}
Assume that $X$ is an abelian surface with $\NS(X)=\bbZ H$ and 
$(H^2)=2$.
Then $\FM(X)=\{X\}$.
\end{prop}

\begin{proof}
If $\NS(X)=\bbZ H$, then a two-dimensional fine moduli of the stable sheaves 
on $X$ can be written as $M_X^H(p^2 r, p q H,q^2 a)$ with $p,q\in\bbZ$ and 
$r,a\in\bbZ_{>0}$ satisfying $\gcd(p r,q a)=1$ and $r a=(H^2)/2$.
Then by Lemma~\ref{lem:FMpartner} we have 
$\FM(X)=\{M_X^H(r,H,a)\mid \gcd(r,a)=1, r a=(H^2)/2 \}/\sim$.
If moreover $(H^2)=2$, then $\FM(X)=\{M_X^H(1,H,1)\}$. 
Since $M_X^H(1,H,1)\cong\widehat{X}\cong X$, we have the consequence.
\end{proof}

Next we recall the $G$ introduced in Definition \ref{dfn:G}. 
In the case $n=1$, $G$ coincides with $\SL(2,\bbZ)$. 
Then the study of numerical solutions reduces 
to the theory of integral quadratic forms.

\begin{dfn}
Two quadratic forms $f_1(x,y)=ax^2+2b x y+c y^2$ and $f_2(x,y)=a' x^2+2b' x y+c' y^2$ defined over $\bbZ$ (that is, $a,b,c,a',b',c'\in\bbZ$) are called \emph{equivalent} if there exists a matrix $A\in\GL(2,\bbZ)$ such that ${}^t A\begin{bmatrix}a&b\\b&c\end{bmatrix}A = \begin{bmatrix}a'&b'\\b'&c'\end{bmatrix}$. We denote it by $f_1(x,y)\equiv f_2(x,y)$.

The discriminant of the quadratic form $f(x,y)=a x^2+2b x y+c y^2$ is defined as $\det f\seteq b^2-a c$.

The class number of quadratic forms of discriminant $D$ is the number of the equivalent classes of quadratic forms whose discriminant are $D$.
\end{dfn}

\begin{thm}\label{thm:class-number}
Let $X$ be a principally polarized abelian surface with $\NS(X)=\bbZ H$.
Let $v=(r,d H,a)$ be a Mukai vector satisfying the following condition.
\begin{align}
\label{cond:cls}
\begin{split}
&\ell\seteq\mpr{v^2}/2\in\bbZ_{>0} \mbox{ is not a square number, and }
\\
&\mbox{the class number of quadratic forms with discriminant $\ell$ is 1.}
\end{split}
\end{align}
Then the birational morphism $M_X^H(v)\cdots\to X\times{\rm Hilb}^{\ell}X$ is given by the following description.
\begin{enumerate}
\item
Take a Mukai vector $v_1=(p_1^2,p_1q_1H,q_1^2)$ which is of minimum rank among those satisfying $2p_1q_1d-p_1^2a-q_1^2r=\epsilon$, where $\epsilon=1$ or $-1$.  %
\item
Then the matrix 
\begin{align*}
\gamma\seteq\pm
\begin{bmatrix}
 q_1&\epsilon(d q_1-a p_1)\\-p_1&\epsilon(d p_1-r q_1)
\end{bmatrix}
\in\PSL(2,\bbZ)
\end{align*}
diagonalizes the matrix $Q_v=\begin{bmatrix}r&d\\d&a\end{bmatrix}$, that is, 
${}^t\gamma Q_v\gamma=-\epsilon\begin{bmatrix}1&0\\0&-\ell\end{bmatrix}$.

\item
Set $v_2\seteq(p_2^2,p_2q_2H,q_2^2)$, where $q_2=-\epsilon(d q_1-a p_1)$ and $p_2=\epsilon(d p_1-r q_1)$. Then a general member of $M_X^H(v)$ has a semi-homogeneous presentation with numerical solution $v=\pm(\ell v_1-v_2)$. 
Moreover the FMT $\Phi\seteq\Phi_{X\to X}^{\bfE }$ 
such that $\theta(\Phi)=\gamma$ or the composition $\calD_X\circ\Phi$ 
gives the birational correspondence $M_X^H(v)\cdots\to X\times\Hilb{\ell}{X}$ 
up to shift.
\end{enumerate}
\end{thm}

\begin{proof}
The condition \eqref{cond:cls} implies that there exists at least one numerical solution for $v$. Then the conclusion follows from Proposition\,\ref{prop:biraional}.
\end{proof}

\begin{rmk}
It is conjectured that there are infinite number of real quadratic forms whose class number is one.
\end{rmk}

We close this paper with the comment 
on the birational types of lower dimensional moduli. 
By the above theorem, one is lead to 
the classification of quadratic forms 
$Q(x,y)=r x^2+2d x y+a y^2$ of small discriminant $\ell=d^2-r a$.
We search in the range $0<\ell<11$.

If $r=0$, then $Q=2 d x y+a y^2=2 d(x+\lambda y)y+(a-2 d \lambda)y^2$ 
for any $\lambda\in\bbZ$. 
Hence we assume $|a|\le d$. 
If $a\neq 0$, then replacing $x$ and $y$ we reduce to the case $r>0$. 
If $r>0$, then we may assume that $0\le d\le r/2$. 
If $a\ge 0$, then from $0\le a =(d/r)d-\ell/r<d/2\le r/4$ 
we have $Q\equiv r' x^2+2d' x y+a' y^2$ with $0\le r'<r$. 

Thus by the action of $\GL(2,\bbZ)$, we may assume
\\
\qquad
(i) $r=a=0$,\quad or
\qquad
(ii) $r\neq 0$ and $0\le 2 d\le r\le -a$.

In the first case $r=a=0$, 
the quadratic form corresponding to a primitive Mukai vector 
is uniquely determined to be $Q=2x y$.

In the second case, $a\le 0$ yields $d\le\sqrt{\ell}$. 
The above discussion allows us to assume $d\le3$. 
If $r=1$, then $d=0$ and we have $(r,d,a)=(1,0,-\ell)$.
If $|a|\ge r>1$, then the possible values of $(r,d,a)$ are 
\begin{align*}
\begin{array}{l l l}
(2,d,-2) \quad d=0,1,2; &
(2,d,-3) \quad d=0,1,2; &
(2,d,-4) \quad d=0,1;   \\
(2,d,-5) \quad d=0;     &
(3,d,-3) \quad d=0,1. 
\end{array}
\end{align*}
We can remove the cases $(r,d,a)=(2,0,-2),(2,0,-4),(3,0,-3)$ 
since these correspond to non-primitive Mukai vectors.
Moreover we can find that $2x^2-3 y^2\equiv 2x^2-4x y-y^2\equiv x^2-6 y$, 
$2x^2+2x y-3 y^2\equiv 2x^2+6x y+y^2\equiv x^2-7y ^2$ and 
$2x^2-5 y^2 \equiv 3x^2+2x y-3y^2$.

We summarize the result in the table below.
\begin{align*}
\begin{tabular}{l|l l || l | l}
$\ell$ & $Q_v$ & & $\ell$ & $Q_v$ \\
\hline
1 & $2x y$, $x^2-y^2$ & &6 & $x^2-6 y ^2$\\
2 & $x^2-2 y^2$       & &7 & $x^2-7 y ^2$\\
3 & $x^2-3 y^2$       & &8 & $x^2-8 y^2$\\
4 & $x^2-4 y^2$       & &9 & $2 x^2+2 x y-4 y^2$, $x^2-9 y^2$\\
5 & $2 x^2+2x y-2y^2$, $x^2-5 y^2$
&&10& $3 x^2+2 x y-3 y^2$, $x^2-10 y^2$\\
\end{tabular}
\end{align*}

Note that if $Q_v$ is not primitive and $\ell>1$, then $M_X^H(v)$ is not birationally equivalent to $X\times\Hilb{\ell}{X}$. For example, if $r,a \in 2\bbZ$, then $Q_v\not\equiv x^2-\ell y^2$. 
\\

The following proposition describes the birational classes of the moduli of dimension $\le 20$.

\begin{prop}
\begin{enumerate}
\item
If $\ell=1,2,3,4,6,7,8$, then $M_X^H(v)$ is birationally equivalent to $X\times \Hilb{\ell}{X}$. 

\item
If $\ell=5$, then the birational type of $M_X^H(v)$ is either $M_X^H(2,H,-2)$ or $X\times \Hilb{\ell}{X}$.

\item
If $\ell=9$, then the birational type of $M_X^H(v)$ is either $M_X^H(2,H,-4)$ or $X\times \Hilb{\ell}{X}$.
\end{enumerate}
\end{prop}
\begin{proof}
We only need to show that $M_X^H(2,H,-2)$, $M_X^H(2,H,-4)$ and $X\times \Hilb{\ell}{X}$ are birationally different. The case $M_X^H(2,H,-2)$ is shown in \cite[Example 4.1]{Yoshioka:2001:ann}. The case $M_X^H(2,H,-4)$ can be shown similarly. 
\end{proof}

\begin{rmk}
In the case $\ell=10$, we cannot distinguish $M_X^H(3,H,-3)$ and $X\times \Hilb{10}{X}$ by the method of \cite[Example 4.1]{Yoshioka:2001:ann}.
\end{rmk}


\begin{thebibliography}{MM}

\bibitem[Hu]{Huybrechts:book}
D.~Huybrechts,
\emph{Fourier-Mukai transforms in algebraic geometry},
Oxford University Press, 2006.

\bibitem[HL]{HuybrechtsLehn:book}
D.~Huybrechts, M.~Lehn,
\emph{The geometry of moduli spaces of sheaves},
Aspects of mathematics, Vieweg, 1997.

\bibitem[M1]{Mukai:1978}
S.~Mukai,
\emph{Semi-homogeneous vector bundles on an abelian variety},
J. Math. Kyoto Univ. {\bf 18} (1978), no. 2, 239--272. 

\bibitem[M2]{Mukai:1979}
S.~Mukai,
\emph{On Fourier functors and their applications to vector bundles on abelian surfaces} (in Japanese),
in the proceeding of "Daisuu-kikagaku Symposium" (T\^ohoku University, June 1979), 76--93. 

\bibitem[M3]{Mukai:1980}
S.~Mukai,
\emph{On classification of vector bundles on abelian surfaces} (in Japanese),
RIMS K\^{o}ky\^{u}roku {\bf 409} (1980), 103--127. 

\bibitem[M4]{Mukai:1981}
S.~Mukai,
\emph{Duality between $D(X)$ and $D(\widehat{X})$ with its application to Picard sheaves},
Nagoya Math. J. {\bf 81} (1981), 153--175. 

\bibitem[M5]{Mukai:1984}
S.~Mukai,
\emph{Symplectic structure of the moduli space of sheaves on an abelian and K3 surface},
Invent. Math. {\bf 77} (1984), 101--116. 

\bibitem[M6]{Mukai:1987:tata}
S.~Mukai,
\emph{On the moduli space of bundles on $K3$ surfaces. I}, in \emph{Vector bundles on algebraic varieties}, 341--413, Tata Inst. Fund. Res., 1987.

\bibitem[M7]{Mukai:1987:ASPM}
S.~Mukai,
\emph{Fourier Functor and its Application to the Moduli of Bundles on an Abelian Variety}, Advanced Studies in Pure Mathematics, {\bf 10} (1987), 515--550, .


\bibitem[M8]{Mukai:1998}
S.~Mukai,
\emph{Abelian variety and spin representation}, preprint, 1998.

\bibitem[Or]{Orlov:2002}
D.~Orlov,
\emph{Derived categories of coherent sheaves on abelian varieties and equivalences between them},
Izv. Math. {\bf 66} (2002), no. 3, 569--594

\bibitem[Y1]{Yoshioka:2001:ann}
K.~Yoshioka,
\emph{Moduli spaces of stable sheaves on abelian surfaces},
Math. Ann. {\bf 321} (2001), 817--884.

\bibitem[Y2]{Yoshioka:2003:twistI}
K.~Yoshioka,
\emph{Twisted stability and Fourier-Mukai transforms I},
Compositio Math. {\bf 138} (2003), 261--288.

\bibitem[Y3]{Yoshioka:2009:stabilityII}
K.~Yoshioka,
\emph{Stability and the Fourier-Mukai transform II},
Compositio Math. {\bf 145} (2009), 112--142.

\bibitem[Y4]{Yoshioka:2009:ann}
K.~Yoshioka,
\emph{Fourier-Mukai transform on abelian surfaces},
to appear in Math. Ann.,
online version available at
http://www.springerlink.com/content/100442/.
\end{thebibliography}
\end{document}